\documentclass[11pt,oneside]{amsart}

\usepackage{bm}

\usepackage{newtxtext}
\usepackage{newtxmath}
\usepackage{fix-cm}
\DeclareMathAlphabet{\mathcal}{OMS}{cmsy}{m}{n}

\usepackage{tikz-cd}

\usepackage{tcolorbox}

\usepackage[arrow,curve,matrix,tips,2cell]{xy}

\usepackage[hmargin=2cm,vmargin=1.5cm]{geometry}

\usepackage{fancyhdr}

\usepackage{amsmath}
\usepackage{amscd}
\usepackage{amssymb}
\usepackage{latexsym}
\usepackage{url}

\usepackage{graphicx}

\newtheorem{theorem}{Theorem}[section]
\newtheorem*{theorem*}{Theorem}
\newtheorem{lemma}[theorem]{Lemma}
\newtheorem*{lemma*}{Lemma}
\newtheorem{corollary}[theorem]{Corollary}
\newtheorem{proposition}[theorem]{Proposition}

\newtheorem{remark}[theorem]{Remark}
\newtheorem{definition}[theorem]{Definition}
\newtheorem*{definition*}{Definition}

\newtheorem{question}[theorem]{Question}
\newtheorem*{question*}{Question}

\newtheorem{example}[theorem]{Example}
\newtheorem{examples}[theorem]{Examples}

\makeatletter
\def\revddots{\mathinner{\mkern1mu\raise\p@
\vbox{\kern7\p@\hbox{.}}\mkern2mu
\raise4\p@\hbox{.}\mkern2mu\raise7\p@\hbox{.}\mkern1mu}}
\makeatother 
\newcommand{\bgl}{\begin{equation}} 
\newcommand{\egl}{\end{equation}}
\newcommand{\bgloz}{\begin{equation*}} 
\newcommand{\egloz}{\end{equation*}}
\newcommand{\bgln}{\begin{eqnarray}} 
\newcommand{\egln}{\end{eqnarray}}
\newcommand{\bglnoz}{\begin{eqnarray*}} 
\newcommand{\eglnoz}{\end{eqnarray*}}
\newcommand{\btheo}{\begin{theorem}}
\newcommand{\etheo}{\end{theorem}}
\newcommand{\btheooz}{\begin{theorem*}}
\newcommand{\etheooz}{\end{theorem*}}

\newcommand{\blemma}{\begin{lemma}}
\newcommand{\elemma}{\end{lemma}}
\newcommand{\blemmaoz}{\begin{lemma*}}
\newcommand{\elemmaoz}{\end{lemma*}}
\newcommand{\bproof}{\begin{proof}}
\newcommand{\eproof}{\end{proof}}
\newcommand{\bbew}{\begin{beweis}}
\newcommand{\ebew}{\end{beweis}}
\newcommand{\bremark}{\begin{remark}\em}
\newcommand{\eremark}{\end{remark}}
\newcommand{\bdefin}{\begin{definition}}
\newcommand{\edefin}{\end{definition}}
\newcommand{\bdefinoz}{\begin{definition*}}
\newcommand{\edefinoz}{\end{definition*}}
\newcommand{\bex}{\begin{example}}
\newcommand{\eex}{\end{example}}
\newcommand{\bexs}{\begin{examples}}
\newcommand{\eexs}{\end{examples}}

\newcommand{\bprop}{\begin{proposition}}
\newcommand{\eprop}{\end{proposition}}
\newcommand{\bcor}{\begin{corollary}}
\newcommand{\ecor}{\end{corollary}}
\newcommand{\bfa}{\begin{cases}} 
\newcommand{\efa}{\end{cases}}

\newcommand{\bquestion}{\begin{question}}
\newcommand{\equestion}{\end{question}}
\newcommand{\bquestionoz}{\begin{question*}}
\newcommand{\equestionoz}{\end{question*}}
\newcommand{\cE}{\mathcal E}

\newcommand{\cG}{\mathcal G}

\newcommand{\cJ}{\mathcal J}
\newcommand{\cK}{\mathcal K}
\newcommand{\cL}{\mathcal L}

\newcommand{\cO}{\mathcal O}

\newcommand{\cR}{\mathcal R}

\newcommand{\cT}{\mathcal T}

\newcommand{\cW}{\mathcal W}
\newcommand{\cX}{\mathcal X}

\def\Cz{\mathbb{C}}
\def\Fz{\mathbb{F}}

\def\Tz{\mathbb{T}}
\def\Zz{\mathbb{Z}}

\def\1z{\mathbb{1}}
\newcommand{\an}[1]{``#1''} 
\newcommand{\ti}{\tilde}

\newcommand{\lori}{\longrightarrow}

\newcommand{\ma}{\mapsto} 
\newcommand{\onto}{\twoheadrightarrow} 
\newcommand{\into}{\hookrightarrow} 

\newcommand{\ve}{\varepsilon}

\def\SEMI{\mbox{$\times\kern-2pt\vrule height5pt width.6pt \kern3pt $}}

\newcommand{\End}{{\rm End}\,}

\newcommand{\Spec}{{\rm Spec\,}} 

\newcommand{\id}{{\rm id}}

\newcommand{\img}{{\rm im\,}}
\renewcommand{\ker}{{\rm ker}\,}

\newcommand{\lcm}{{\rm lcm}} 
\newcommand{\reg}{^\times} 
\newcommand{\lspan}{{\rm span}} 
\newcommand{\clspan}{\overline{\lspan}} 
\newcommand{\abs}[1]{\lvert#1\rvert} 
\newcommand{\defeq}{\mathrel{:=}} 
\newcommand{\eqdef}{\mathrel{=:}} 

\newcommand{\dop}{\text{: }} 
\newcommand{\ilim}{\varinjlim} 
\newcommand{\lge}{\left\{} 
\newcommand{\rge}{\right\}} 
\newcommand{\lru}{\left(} 
\newcommand{\rru}{\right)} 
\newcommand{\lsp}{\left\langle} 
\newcommand{\rsp}{\right\rangle} 
\newcommand{\rukl}[1]{\lru #1 \rru} 
\newcommand{\gekl}[1]{\lge #1 \rge} 
\newcommand{\spkl}[1]{\lsp #1 \rsp} 
\newcommand{\menge}[2]{\gekl{ #1 \dop #2 }} 
\begin{document}

\title{C*-algebras of\\right LCM one-relator monoids and Artin-Tits monoids of finite type}

\thispagestyle{fancy}

\author{Xin Li}
\address{School of Mathematics and Statistics, University of Glasgow, University Place, Glasgow G12 8QQ, United Kingdom}
\email{Xin.Li@glasgow.ac.uk}

\author{Tron Omland}
\address{Department of Mathematics, University of Oslo, NO-0316 Oslo, Norway \\ and \newline Department of Computer Science, Oslo Metropolitan University, NO-0130 Oslo, Norway}
\email{trono@math.uio.no}

\author{Jack Spielberg}
\address{School of Mathematical and Statistical Sciences, Arizona State University, Tempe, AZ 85287-1804, USA}
\email{jack.spielberg@asu.edu}


\thanks{The second author is funded by the Research Council of Norway through FRINATEK, project no.~240913.}

\begin{abstract}
We study C*-algebras generated by left regular representations of right LCM one-relator monoids and Artin-Tits monoids of finite type. We obtain structural results concerning nuclearity, ideal structure and pure infiniteness. Moreover, we compute K-theory. Based on our K-theory results, we develop a new way of computing K-theory for certain group C*-algebras and crossed products.
\end{abstract}

\maketitle

\setcounter{tocdepth}{3}

\setlength{\parindent}{0cm} \setlength{\parskip}{0.5cm}

\section{Introduction}

Ever since the work in \cite{Co,Cun77}, C*-algebras generated by isometries have provided interesting examples in the theory of general C*-algebras. In particular, C*-algebras generated by left regular representations of left-cancellative semigroups, also called semigroup C*-algebras, form a natural class of C*-algebras which have been studied for several types of semigroups (see \cite{CELY} and the references therein). This construction has attracted attention because of connections to several topics such as index theory, representation theory, amenability, and non-selfadjoint operator algebras. Recently, it was discovered in \cite{CDL,Li14,Li16a,Li16b} that there is even a connection to number theory. This observation has led to a renewed interest in the subject. Apart from these examples from number theory, there is another natural class of examples coming from group theory, for instance semigroups given by generators and relations.

This present paper makes a contribution towards a better understanding of the C*-algebras attached to this second class of monoids. More precisely, we study semigroup C*-algebras of two classes of such monoids: one-relator monoids and Artin-Tits monoids of finite type. Our motivating example is the classical Braid monoid $B_n^+$ which sits inside the Braid group $B_n$. Geometrically speaking, elements of $B_n^+$ can be thought of as braids in which the crossings are allowed to go in one direction only. The (reduced) semigroup C*-algebra of $B_n^+$ turns out to be an extension of the reduced group C*-algebra of $B_n$ by another C*-algebra which is very similar to a graph C*-algebra \cite{Rae} in the case $n=3$ and which is not nuclear in the case $n \geq 4$. Our search for other classes of monoids whose C*-algebras have a similar structure as the one of $B_3^+$ has led us to one-relator monoids, while monoids whose C*-algebras behave like the ones of $B_n^+$, for $n \geq 4$, are given by more general Artin-Tits monoids (see \cite{BS,Del,Sal}). In this paper, we mainly focus on structural results and K-theory computations for these classes of semigroup C*-algebras. The latter turn out to be interesting because they lead to a new way of computing K-theory for certain group C*-algebras and crossed products.

Let us now describe our main results in more detail. For a class of cancellative right LCM one-relator monoids as specified in \S~\ref{sec:nonOre} and \S~\ref{sec:Ore}, we show that there is a big difference between the reversible and the non-reversible cases. In the non-reversible case, we show that the boundary quotient (in the sense of \cite[\S~5.7]{CELY}) of our semigroup C*-algebra is purely infinite simple (see Corollary~\ref{Cor:bdq:pis_gen}, where this is proven in even greater generality). It follows that the semigroup C*-algebra is the extension of a purely infinite simple C*-algebra by the algebra of compact operators if the underlying monoid is finitely generated, and it is purely infinite simple itself if the monoid is not finitely generated (see Corollary~\ref{Cor:OneRelM-Omega=P+bd}). Moreover, this boundary quotient can be identified with a graph C*-algebra \cite{Rae} in special cases, so that we can further deduce that our semigroup C*-algebra is nuclear (see Corollary~\ref{Cor:bd=graph_nonOre}). 

In the reversible case, the boundary quotient is given by the reduced C*-algebra of the one-relator group given by the same presentation as our one-relator monoid. This means that the boundary quotient is almost never nuclear since most one-relator groups are not amenable. However, we show that in special cases, the kernel of the canonical projection map onto the boundary quotient is a Toeplitz graph C*-algebra, in particular nuclear (see Theorem~\ref{Thm:OreOneRel:graph}).

The main advantage of our new approach is that we are able to identify general conditions on the presentations defining our monoids which allow us to deduce the above-mentioned structural results. This is in contrast to previous work, which focused on case-by-case studies of classes of monoids defined by generators and relations, where the methods had to be adapted to the defining relations.

Our analysis of one-relator monoids prompts the question whether similar results still hold in the case of more relations. We show that in general, this is not the case. For C*-algebras attached to Artin-Tits monoids of finite type, which are not one-relator monoids, not only does the boundary quotient fail to be nuclear, but even the kernel of the canonical projection onto the boundary quotient is typically not nuclear (see Proposition~\ref{prop:nonamenable}). However, we show that this kernel is still purely infinite simple for irreducible Artin-Tits monoids with at least three generators (see Theorem~\ref{Thm:simplepi}). Our results complement the work in \cite{CL02,CL07,ELR} on C*-algebras of right-angled Artin-Tits monoids.

Finally, we use our analysis of the ideal corresponding to the boundary quotient to provide general tools for K-theory computations. We focus on dihedral Artin-Tits groups and torus knot groups. In each of these cases, we construct two six-term exact sequences in K-theory. The first one allows to compute K-theory for the ideal corresponding to the boundary quotient, while the second one requires K-theory of the ideal as an input and determines K-theory of the reduced group C*-algebras of dihedral Artin-Tits groups and torus knot groups. All this works for arbitrary coefficients and therefore allows for K-theory computations for crossed products (see Theorem~\ref{THM:K_examples}). We then present explicit K-theory computations for a variety of examples, including group C*-algebras of dihedral Artin-Tits groups, torus knot groups, the braid group $B_4$ and a semidirect product group arising from Artin's representation of braid groups (see \S~\ref{ss:Ex-K_*}).

The main ingredient for our computations is a K-theory formula for crossed products in the presence of the independence condition (see \cite{CEL1,CEL2}). Therefore, our K-theory results point towards the possibility of using K-theory formulas for semigroup C*-algebras and closely related crossed products as in \cite{CEL1,CEL2} to compute K-theory for group C*-algebras and crossed products. Moreover, since we want to allow arbitrary coefficients, we introduce graph C*-algebras twisted by coefficient algebras (a special case of this construction appears in \cite{Cun81}). This construction might be of independent interest.

We would like to thank Nathan Brownlowe and Dave Robertson for inviting us to the workshop \an{Interactions Between Semigroups and Operator Algebras} in Newcastle (Australia), where this project has been initiated. Moreover, we thank the referees for many helpful comments which improved our paper considerably.

\section{Preliminaries}

Let us start by reviewing some basics concerning semigroups and their C*-algebras.
\setlength{\parindent}{0cm} \setlength{\parskip}{0cm}

\subsection{Monoids}

Semigroups with an identity are called monoids in this paper. Although it is often not strictly necessary, we will focus on monoids throughout, just to keep things simple and because all our examples do have an identity.

\subsubsection{Monoids defined by presentations}
\label{sss:MonoidsPres}

Let $\Sigma$ be a countable set, $\Sigma^*$ the free monoid generated by $\Sigma$ (also viewed as the set of finite words in $\Sigma$) and $R$ a subset of $\Sigma^* \times \Sigma^*$. We call $(\Sigma,R)$ a presentation, where $\Sigma$ is the set of generators and $R$ is the set of relations. Words $u \in \Sigma^*$ for which $(u,v)$ or $(v,u)$ lies in $R$ are called relators. Given a presentation $(\Sigma,R)$, we form the monoid $P = \spkl{\Sigma \, \vert \, R}^+$ generated by $\Sigma$ subject to the relations $u = v$ for all $(u,v) \in R$ (see \cite[\S~1.12]{CP1}). Note that countability of $\Sigma$ is not essential; this assumption is included only to guarantee that our semigroup C*-algebras are separable.
\setlength{\parindent}{0cm} \setlength{\parskip}{0.5cm}

Given $x \in \Sigma^*$, we also denote by $x$ the element in $P$ represented by $x$. Given $x, y \in \Sigma^*$, we write $x \equiv y$ if $x$ and $y$ coincide in $\Sigma^*$, i.e., $x$ is the same word as $y$, and $x = y$ if $x$ and $y$ represent the same element in $P$. We write $\ell^*$ for the length of a word in $\Sigma^*$ and $\ell(x) \defeq \min \menge{\ell^*(w)}{w \in \Sigma^*, \, w = x}$ for the word length of $x \in P$ with respect to $\Sigma$. Given $\sigma \in \Sigma$ and $x \in \Sigma^*$, $\ell_{\sigma}(x)$ counts how many times $\sigma$ appears in $x$. We write $\varepsilon$ for the empty word in $\Sigma^*$.

Recall how $P = \spkl{\Sigma \, \vert \, R}^+$ is constructed: For two words $w, w' \in \Sigma^*$, we write $w \equiv_R w'$ if $w \equiv w'$ or there exist $a, u, v, z \in \Sigma^*$ with $(u,v) \in R$ or $(v,u) \in R$ such that $w \equiv a u z$ and $w' \equiv a v z$. Two words $x, y \in \Sigma^*$ represent the same element in $P$ if and only if there is a finite sequence $w_1, \dotsc, w_s \in \Sigma^*$ with $w_1 \equiv x$, $w_s \equiv y$, and for all $1 \leq i \leq s-1$, we have $w_i \equiv_R w_{i+1}$.
\setlength{\parindent}{0cm} \setlength{\parskip}{0cm}

\subsubsection{Standing assumptions on our presentations}
\label{sss:Standing}

All our presentations will satisfy that for all $(u,v) \in R$, we have $u \neq \varepsilon \neq v$. It then immediately follows from the construction of $P = \spkl{\Sigma \, \vert \, R}^+$ that $P$ does not contain any left- or right-invertible elements other than the identity. In particular, the group of units $P^*$ in $P$ is trivial, i.e., $P^* = \gekl{e}$, where $e$ is the identity in $P$. 
\setlength{\parindent}{0cm} \setlength{\parskip}{0.5cm}

In the following, we will also always exclude the (degenerate) case that a generator $\sigma \in \Sigma$ is redundant, i.e., we always assume that for every $\sigma \in \Sigma$ and $w \in (\Sigma \setminus \gekl{\sigma})^*$ that $\sigma \neq w$. Otherwise, if $\sigma \in \Sigma$ satisfies $\sigma = w$ for some $w \in (\Sigma \setminus \gekl{\sigma})^*$, then we could always give another presentation for $P$, $P \cong \spkl{\Sigma \setminus \gekl{\sigma} \, \vert \, R'}^+$, where $R'$ is obtained from $R$ by replacing every $(u,v) \in R$ by $(u',v')$, where we obtain $u'$ from $u$ and $v'$ from $v$ by replacing every $\sigma$ in $u$ and $v$ by $w$.

Moreover, we will always assume (without loss of generality) that $R$ does not contain $(u,u)$ for any $u \in \Sigma^*$.

Note that under these standing assumptions, $P = \spkl{\Sigma \, \vert \, R}^+$ is finitely generated if and only if $\abs{\Sigma} < \infty$.
\setlength{\parindent}{0cm} \setlength{\parskip}{0cm}

\subsubsection{One-relator monoids}
\label{sss:OneRelM}

We will mainly focus on the case of one-relator monoids, i.e., monoids of the form $P = \spkl{\Sigma \, \vert \, R}^+$ where $R$ consists of a single pair $(u,v) \in \Sigma^* \times \Sigma^*$. This is (arguably) the most basic case. Moreover, in examples, we actually see very different phenomena when we drop this assumption (see \S~\ref{ss:ATM_FiniteType}). 
\setlength{\parindent}{0cm} \setlength{\parskip}{0.5cm}

We often write $(\Sigma, u=v)$ for $(\Sigma, R)$ and $\spkl{\Sigma \, \vert \, u = v}^+$ for $\spkl{\Sigma \, \vert \, R}^+$. Let $G = \spkl{\Sigma \, \vert \, u = v}$ be the group generated by $\Sigma$ subject to the relation $u = v$. Let us assume that the first letter of $u$ is not equal to the first letter of $v$ and that the last letter of $u$ is not equal to the last letter of $v$. Then it was proven in \cite{Ad} that the canonical map $P \to G$ induced by the identity on $\Sigma$ is injective. In particular, $P$ is cancellative (i.e., left and right cancellative). Conversely, since we assume that $u \not\equiv v$, if $P$ is cancellative, then the first and last letters of $u$ and $v$ must differ.

Therefore, as we do want cancellation to hold, whenever we consider one-relator monoids of the form $\spkl{\Sigma \, \vert \, u = v}^+$, we will assume that the first letter of $u$ is not equal to the first letter of $v$ and that the last letter of $u$ is not equal to the last letter of $v$.

\bremark
$G = \spkl{\Sigma \, \vert \, u = v}$ for $u, v \in \Sigma^*$ is a one-relator group. If the first and last letters of $u$ and $v$ differ, then it follows from \cite[Theorem~4.12]{MKS} that $G$ is torsion-free. Moreover, $G$ is amenable if and only if $G$ is cyclic or $G \cong \spkl{a,b \, \vert \, ab = b^ka}$ for some $0 \neq k \in \Zz$ (see for instance \cite{CG}).
\eremark
\setlength{\parindent}{0cm} \setlength{\parskip}{0cm}

\subsubsection{The LCM property}
\label{sss:LCM}

We will assume throughout \S~\ref{sec:nonOre} -- \ref{sec:exseq-K} that all our monoids $P$ are right LCM, i.e., given $p, q \in P$, we have $pP \cap qP = \emptyset$ or there exists $r \in P$ with $pP \cap qP = rP$. 
\setlength{\parindent}{0cm} \setlength{\parskip}{0.5cm}

Applying the notion of right reversing in the sense of \cite{Deh03} to our presentation $(\Sigma, u=v)$, let us present a sufficient condition for the right LCM property: 
\setlength{\parindent}{0.5cm} \setlength{\parskip}{0cm}

Let us introduce another copy of $\Sigma$ by forming $\Sigma \cup \Sigma^{-1} = \menge{\sigma}{\sigma \in \Sigma} \cup \menge{\sigma^{-1}}{\sigma \in \Sigma}$. For $x, y \in (\Sigma \cup \Sigma^{-1})^*$, we write $x \curvearrowright_r y$ if there exist $w_1, w_2, \dotsc, w_j \in (\Sigma \cup \Sigma^{-1})^*$ with $x \equiv w_1$, $y \equiv w_j$, and for all $1 \leq i \leq j-1$, we have $w_i \equiv \alpha \sigma^{-1} \tau \zeta$ for $\sigma, \tau \in \Sigma$ and $\alpha, \zeta \in (\Sigma \cup \Sigma^{-1})^*$, and $w_{i+1} \equiv \alpha \zeta$ if $\sigma \equiv \tau$, and $w_{i+1} \equiv \alpha s t^{-1} \zeta$ if $\sigma s \equiv u$ and $\tau t \equiv v$. Then, following \cite[Definition~3.1]{Deh03}, we say that $(\Sigma, u=v)$ satisfies the strong $\boldsymbol{r}$-cube condition on $\Sigma$ if for all $\sigma, \tau, \upsilon \in \Sigma$, $\sigma^{-1} \tau \tau^{-1} \upsilon \curvearrowright_r xy^{-1}$ for $x, y \in \Sigma^*$ implies that $(\sigma x)^{-1} (\upsilon y) \curvearrowright_r \varepsilon$. Note that for a word $w = \sigma_1 \sigma_2 \dotsm \sigma_{l-1} \sigma_l \in \Sigma^*$, we write $w^{-1} = \sigma_l^{-1} \sigma_{l-1}^{-1} \dotsm \sigma_2^{-1} \sigma_1^{-1} \in (\Sigma^{-1})^*$. 

Now, in our case of the presentation $(\Sigma, u=v)$, it is easy to verify the strong $\boldsymbol{r}$-cube condition on $\Sigma$: If $\sigma \not\equiv \tau$ and $\sigma, \tau$ are not the first letters of $u$ and $v$, then we cannot have $\sigma^{-1} \tau \tau^{-1} \upsilon \curvearrowright_r xy^{-1}$ for some $x, y \in \Sigma^*$. Similarly, if $\tau \not\equiv \upsilon$ and $\tau, \upsilon$ are not the first letters of $u$ and $v$, then we cannot have $\sigma^{-1} \tau \tau^{-1} \upsilon \curvearrowright_r xy^{-1}$ for some $x, y \in \Sigma^*$. If $\sigma \equiv \tau \equiv \upsilon$, then $x \equiv y \equiv \varepsilon$ and the strong $\boldsymbol{r}$-cube condition on $\Sigma$ is obviously satisfied. If $\sigma \equiv \tau$ and $u \equiv \tau x$, $v \equiv \upsilon y$, then $\sigma^{-1} \tau \tau^{-1} \upsilon \equiv \tau^{-1} \tau \tau^{-1} \upsilon \curvearrowright_r xy^{-1}$, and $(\tau x)^{-1} (\upsilon y) \curvearrowright_r \varepsilon$. The case $u \equiv \sigma x$, $v \equiv \tau y$ and $\tau \equiv \upsilon$ is analogous. Finally, if $\sigma \equiv \upsilon$ and $u \equiv \sigma s \equiv \upsilon s$, $v \equiv \tau t$, then $\sigma^{-1} \tau \tau^{-1} \upsilon \curvearrowright_r s t^{-1} t s^{-1} \curvearrowright_r ss^{-1}$, so that $x \equiv s$, $y \equiv s$, and $(\sigma s)^{-1} (\upsilon s) \equiv (\sigma s)^{-1} (\sigma s) \curvearrowright_r \varepsilon$. This shows the strong $\boldsymbol{r}$-cube condition on $\Sigma$. 

In order to deduce the strong $\boldsymbol{r}$-cube condition (on $\Sigma^*$), we need to check whether $(\Sigma, u=v)$ is $\boldsymbol{r}$-homogeneous in the sense of \cite[Definition~4.1]{Deh03}. This is for instance the case if $\ell^*(u) = \ell^*(v)$ (for we can take $\lambda = \ell^*$ in \cite[Definition~4.1]{Deh03}). $(\Sigma, u=v)$ is also $\boldsymbol{r}$-homogeneous if $\ell^*(u) < \ell^*(v)$ and there exists $\sigma \in \Sigma$ with $\ell_{\sigma}(u) > \ell_{\sigma}(v)$. In that case, let $\delta \defeq \ell^*(v) - \ell_{\sigma}(v) - (\ell^*(u) - \ell_{\sigma}(u))$. Then $0 < \ell_{\sigma}(u) - \ell_{\sigma}(v) < \delta$. Let $\eta$ and $\zeta$ be positive integers with $\delta \eta = \lcm(\delta, \ell_{\sigma}(u) - \ell_{\sigma}(v)) = (\ell_{\sigma}(u) - \ell_{\sigma}(v)) \zeta$. Then define $\lambda$ as the monoid homomorphism from $\Sigma^*$ to the non-negative integers by setting $\lambda(\sigma) = \zeta$ and $\lambda(\tau) = \eta$ for all $\tau \in \Sigma \setminus \gekl{\sigma}$. We then have $\lambda(u) = \zeta \ell_{\sigma}(u) + \eta (\ell^*(u) - \ell_{\sigma}(u)) = \delta \eta + \zeta \ell_{\sigma}(v) + \eta (\ell^*(u) - \ell_{\sigma}(u)) = \zeta \ell_{\sigma}(v) + \eta (\delta + \ell^*(u) - \ell_{\sigma}(u)) = \zeta \ell_{\sigma}(v) + \eta (\ell^*(v) - \ell_{\sigma}(v)) = \lambda(v)$. Hence $(\Sigma, u=v)$ is $\boldsymbol{r}$-homogeneous. Once we know that $(\Sigma, u=v)$ is $\boldsymbol{r}$-homogeneous, \cite[Proposition~4.4]{Deh03} tells us that $(\Sigma, u=v)$ is $\boldsymbol{r}$-complete in the sense of \cite[Definition~2.1]{Deh03}. This in turn implies that $P = \spkl{\Sigma, u=v}^+$ is right LCM by \cite[Proposition~6.10]{Deh03}.
\setlength{\parindent}{0cm} \setlength{\parskip}{0.5cm}

\bremark
\label{REM:stillLCM}
Note that there are examples of presentations $(\Sigma, u=v)$ which are not $\boldsymbol{r}$-homogeneous but where $\spkl{\Sigma, u=v}^+$ is still right LCM. For instance, let $\Sigma = \gekl{a,b}$ and $u = b^d a b^c$ for some positive integers $c, d$, and $v = a$. If we define $x_i \defeq a b^{ic}$ for $i = 0, 1, \dotsc$, then $x_i = b^d x_{i+1}$ for all $i$. Thus $\spkl{\Sigma \, \vert \, u=v}^+$ is not $\boldsymbol{r}$-Noetherian in the sense of \cite[Definition~2.6]{Deh03}, so that $(\Sigma, u=v)$ cannot be $\boldsymbol{r}$-homogeneous by \cite[Proposition~4.3]{Deh03}. However, $\spkl{\Sigma \, \vert \, u=v}^+$ is right LCM by \cite[Proposition~2.10]{Sp12}.
\eremark
\setlength{\parindent}{0cm} \setlength{\parskip}{0cm}

\subsubsection{Reversibility}
\label{sss:Rev}

Let us discuss the question when $P = \spkl{\Sigma \, \vert \, u=v}^+$ is left reversible, i.e., for all $p, q \in P$, we have $pP \cap qP \neq \emptyset$. We first observe that if $(\Sigma, u=v)$ is $\boldsymbol{r}$-homogeneous, then by \cite[Proposition~6.7]{Deh03}, $P = \spkl{\Sigma \, \vert \, u=v}^+$ is left reversible if and only if there exists $\Sigma \subseteq \Sigma' \subseteq \Sigma^*$ such that for all $w, x \in \Sigma'$, there exist $y, z \in \Sigma'$ with $(xz)^{-1}(wy) \curvearrowright_r \varepsilon$. In particular, $P = \spkl{\Sigma \, \vert \, u=v}^+$ is left reversible if there exists $\Sigma \subseteq \Sigma' \subseteq \Sigma^*$ such that for all $w, x \in \Sigma'$, there exist $y, z \in \Sigma'$ with $x^{-1}w \curvearrowright_r zy^{-1}$ (see \cite[Remark~6.9]{Deh03}).
\setlength{\parindent}{0cm} \setlength{\parskip}{0.5cm}

If $\abs{\Sigma} \geq 3$, then $P$ cannot be left reversible: Take $\sigma \in \Sigma$ such that $\sigma$ is not the first letter of $u$ and also not the first letter of $v$. Then for every $\tau \in \Sigma \setminus \gekl{\sigma}$, we must have $\sigma P \cap \tau P = \emptyset$ because we cannot find $t, w \in \Sigma^*$ with $tuw \equiv tvw$ such that the first letter of $tuw$ is $\sigma$ and the first letter of $tvw$ is $\tau$ (or vice versa).

Moreover, let us present examples of presentations $(\gekl{a,b}, u=v)$ which do not give rise to left reversible monoids. Define ${\rm OVL}(v) = \menge{x \in \gekl{a,b}^*}{v \equiv xy \equiv wx \ {\rm for} \ {\rm some} \ \varepsilon \neq w, y \in \gekl{a,b}^*}$. Assume that ${\rm OVL}(v) = \gekl{\varepsilon}$. Further suppose that $\ell_a(u) < \ell_a(v)$ or $\ell_b(u) < \ell_b (v)$. The first condition implies that $(\gekl{a,b}, u=v)$ satisfies the Church-Rosser condition (by \cite{Tru}), while the second condition ensures that $(\gekl{a,b}, u=v)$ is noetherian. (The reader may consult \cite[\S~1.1]{BO} for explanations of these notions.) This implies that given $w \in \gekl{a,b}^*$, we can find a unique irreducible word $x \in \gekl{a,b}^*$ (i.e., $x$ does not contain $v$ as a subword) such that $w = x$ in $P$, and two irreducible words $x, y \in \gekl{a,b}^*$ represent the same element in $P$ if and only if $x \equiv y$ (see \cite[Theorem~1.1.12]{BO}). Now assume that ${\rm OVL}(v) = \gekl{\varepsilon}$, and that $\ell_a(u) < \ell_a(v)$ or $\ell_b(u) < \ell_b (v)$. If the first letter of $v$ is $a$, but $v \not\equiv ab^k$ for some $k \geq 1$, then $P$ is not left reversible: Let $p = b^m$, $q = a b^n$ for some $m, n \geq \ell^*(v)$. Take $x, y \in \gekl{a,b}^*$ irreducible. Then $px$ and $qy$ are irreducible, but $px \neq qy$ as $px \not\equiv qy$. Hence $pP \cap qP = \emptyset$.
\setlength{\parindent}{0cm} \setlength{\parskip}{0cm}

\subsection{Semigroup C*-algebras}
\label{SgpC}

Let $P$ be a left cancellative semigroup and $\lambda$ its left regular representation on $\ell^2 P$. This means that for $p \in P$, $\lambda(p)$ is the isometry $\ell^2 P \to \ell^2 P$ determined by $(\lambda(p)\xi)(px) = \xi(x)$ for all $x \in P$ and $(\lambda(p)\xi)(y) = 0$ if $y \notin pP$. We denote by $C^*_{\lambda}(P)$ the C*-algebra generated by $\lambda(P)$, and call it the (reduced) semigroup C*-algebra of $P$. Let $S$ be the inverse semigroup of partial isometries on $\ell^2 P$ generated by $\lambda(P)$. Let $D_{\lambda}(P) \defeq C^*(\menge{s^*s}{s \in S}$. If $P$ embeds into a group $G$, then we can also describe $D_{\lambda}(P)$ by $D_{\lambda}(P) = C^*_{\lambda}(P) \cap \ell^{\infty}(P)$. $D_{\lambda}(P)$ is always a commutative C*-algebra, and we define $\Omega_P \defeq \Spec \rukl{D_{\lambda}(P)}$. We refer the reader to \cite{Li12,Li13} as well as \cite[\S~5]{CELY} for more details about general semigroup C*-algebras.

\subsubsection{Description as groupoid C*-algebras}

If $P$ is a right LCM monoid, then $D_{\lambda}(P) = C^*(\menge{1_{pP}}{p \in P}) \subseteq \ell^{\infty}(P)$. Let $\cJ_P \defeq \menge{pP}{p \in P} \cup \gekl{\emptyset}$. $\cJ_P$ is a semilattice under intersection because $P$ is right LCM. Points in $\Omega_P$, i.e., characters on $D_{\lambda}(P)$, are in bijection to multiplicative surjective maps $\chi: \: \cJ_P \to \gekl{0,1}$ with $\chi(\emptyset) = 0$ (such a map $\chi$ corresponds to the character $\omega: \: D_{\lambda}(P) \to \Cz$ determined by $\omega(1_{pP}) = \chi(pP)$). The topology on $\Omega_P$ corresponds to the topology of point-wise convergence for maps $\cJ_P \to \gekl{0,1}$. In the following, we identify $\Omega_P$ with multiplicative surjective maps $\chi: \: \cJ_P \to \gekl{0,1}$ as explained above.
\setlength{\parindent}{0cm} \setlength{\parskip}{0.5cm}

If $P$ embeds into a group $G$, then we can fix such an embedding and view $P$ as a submonoid of $G$, and consider the partial action $G \curvearrowright \Omega_P$ determined by the partial homeomorphisms $U_{g^{-1}} \to U_g, \, \chi \ma g.\chi$, where $U_{g^{-1}}$ is the subspace of all $\chi \in \Omega_P$ which satisfy $\chi(qP) = 1$ for some $q \in P$ such that $g = pq^{-1}$ for some $p \in P$. For such $\chi$, $g.\chi$ is determined by $(g.\chi)(xP) = \chi(qyP)$ if $xP \cap pP = pyP$ and $(g.\chi)(xP) = 0$ if $xP \cap pP = \emptyset$. Note that $U_{g^{-1}}$ is non-empty if and only if $g = pq^{-1}$ for some $p, q \in P$. As explained in \cite[Theorem~5.6.41]{CELY}, the semigroup C*-algebra of $P$ is canonically isomorphic to the reduced crossed product attached to this partial dynamical system $G \curvearrowright \Omega_P$ as well as isomorphic to the reduced C*-algebra of the corresponding partial transformation groupoid, $C^*_{\lambda}(P) \cong C^*_r(G \ltimes \Omega_P)$.
\setlength{\parindent}{0cm} \setlength{\parskip}{0cm}

\subsubsection{Distinguished subspaces of the character space}
\label{sss:DistSubspOmega}

For brevity, let us drop indices and write $\cJ \defeq \cJ_P$, $\cJ\reg = \cJ \setminus \gekl{\emptyset}$ and $\Omega \defeq \Omega_P$. Let us describe distinguished elements and subspaces of $\Omega$. First of all, every $pP \in \cJ$ determines a point $\chi_{pP} \in \Omega$ given by $\chi_{pP}(xP) = 1$ if $pP \subseteq xP$ and $\chi_{pP}(xP) = 0$ if $pP \nsubseteq xP$. This allows us to identify $\cJ\reg$ with a subset of $\Omega$. It is easy to see that $\cJ\reg$ is dense in $\Omega$. We define $\Omega_{\infty} = \Omega \setminus \cJ\reg$. Clearly, $\cJ\reg$ and $\Omega_{\infty}$ are $G$-invariant subspaces. In the case where $P^* = \gekl{e}$, we have a bijection $P \cong \cJ\reg, \, p \ma pP$, so that we will not distinguish between $P$ and $\cJ\reg$, and hence we have $\Omega_{\infty} = \Omega \setminus P$. Among the points in $\Omega_{\infty}$, we single out those $\chi \in \Omega$ for which $\chi^{-1}(1)$ is maximal, i.e., whenever $\omega \in \Omega$ satisfies $\omega(xP) = 1$ for all $xP \in \cJ$ with $\chi(xP) = 1$, then we must have $\omega = \chi$. We set $\Omega_{\max} \defeq \menge{\chi \in \Omega}{\chi^{-1}(1) \ {\rm is} \ {\rm maximal}}$. Note that $\Omega_{\max} \subseteq \Omega_{\infty}$. Moreover, we define $\partial \Omega \defeq \overline{\Omega_{\max}}$. Let us now collect a few facts about $\partial \Omega$, which are obtained in \cite[\S~5.7]{CELY} in greater generality than needed here. $\partial \Omega$ is the minimal non-empty closed $G$-invariant subspace of $\Omega$. Moreover, $\partial \Omega$ reduces to a single point (namely $\chi \in \Omega$ given by $\chi(xP) = 1$ for all $x \in P$; we usually denote this $\chi$ by $\infty$) if and only if $P$ is left reversible. If $\partial \Omega$ is not a point, then $G \curvearrowright \partial \Omega$ is purely infinite. Since $\partial \Omega$ is always a closed $G$-invariant subspace of $\Omega$, we can define a quotient of $C^*_{\lambda}(P)$ by setting $\partial C^*_{\lambda}(P) \defeq C^*_r(G \ltimes \partial \Omega)$. $\partial C^*_{\lambda}(P)$ is called the boundary quotient of $C^*_{\lambda}(P)$.

\subsubsection{Characters attached to words}
\label{sss:Char_w}

Now let us assume that our monoid $P$ is given by a presentation $(\Sigma, R)$ as in \S~\ref{sss:MonoidsPres} satisfying the standing assumptions in \S~\ref{sss:Standing}. Assume that $P$ is right LCM. Let $G$ be the group given by the same presentation as $P$, and suppose that $P \to G$ induced by $\id_{\Sigma}$ is an embedding (this must be the case if $P$ embeds into a group). Since $P^* = \gekl{e}$ we may write $P$ for the subset $\cJ$ of $\Omega$. In this setting, let us describe $\Omega_{\infty}$ by infinite words (with letters) in $\Sigma$. We write $\Sigma^{\infty}$ for the set of all these infinite words. Given $w \in \Sigma^{\infty}$, define $w_{[1,i]}$ to be the word consisting of the first $i$ letters of $w$ (for $i \geq 1$). Let us write $w_i \defeq w_{[1,i]}$. Then define $\chi_w \in \Omega$ by setting $\chi_w(zP) = 1$ if and only if $w_i \in zP$ for some $i$. It is easy to see that $\chi_w \in \Omega_{\infty}$. Conversely, we have
\blemma
\label{LEM:chi=chi_w}
Every $\chi \in \Omega_{\infty}$ is of the form $\chi_w$ for some $w \in \Sigma^{\infty}$.
\elemma
\setlength{\parindent}{0cm} \setlength{\parskip}{0cm}

\bproof
Given $\chi \in \Omega_{\infty}$, let us enumerate $zP \in \cJ$ with $\chi(zP) = 1$, so that $\menge{zP}{\chi(zP) = 1} = \gekl{z_1 P, z_2 P, \dotsc}$. Then choose $y_i \in \Sigma^*$ with $y_i P = z_1 P \cap \dotso \cap z_i P$ and $y_{i+1} \in y_i \Sigma^*$. Note that all $y_i P$ themselves lie in $\gekl{z_1 P, z_2 P, \dotsc}$. Let $l_i \defeq \ell^*(y_i)$. Since $\chi \in \Omega_{\infty}$, we must have $\lim_{i \to \infty} l_i = \infty$. So there is a unique infinite word $w$ in $\Sigma$ such that $w_{l_i} = y_i$. We claim that $\chi = \chi_w$. Indeed, given $x \in P$, we have $\chi(xP) = 1$ if and only if $xP \in \gekl{z_1P, z_2P, \dotsc}$ if and only if $y_i \in xP$ for some $i$ if and only if $w_j \in xP$ for some $j$ if and only if $\chi_w(xP) = 1$.
\eproof

In this picture with infinite words, the $G$-action on $\Omega_{\infty}$ is given as follows: $\chi \in \Omega_{\infty}$ lies in $U_{g^{-1}}$ if and only if $\chi = \chi_w$ for some $w \in \Sigma^{\infty}$ which starts with $q \in \Sigma^*$, i.e., of the form $w \equiv qx$ for some $x \in \Sigma^{\infty}$, and $g = pq^{-1}$ for some $p \in \Sigma^*$, and in this case we have $g.\chi = g.\chi_w = \chi_{px}$.
\setlength{\parindent}{0cm} \setlength{\parskip}{0.5cm}

The topology on $\Omega_{\infty}$ is the subspace topology from $\Omega$, which in turn is determined by the basic open subspaces $U(pP; p_1P, \dotsc, p_k P) \defeq \menge{\chi \in \Omega}{\chi(pP) = 1, \, \chi(p_1 P) = \dotso = \chi(p_k P) = 0}$, where $p_j P \subsetneq pP$ for all $1 \leq j \leq k$. So basic open subspaces of $\Omega_{\infty}$ are given by $\Omega_{\infty} \cap U(pP; p_1P, \dotsc, p_k P)$. 

In general, two different infinite words can give rise to the same element in $\Omega_{\infty}$. However, this cannot happen for infinite words which contain no relator as a finite subword:
\blemma
\label{LEM:w=w}
For $w, \bar{w} \in \Sigma^{\infty}$, let $w_i \defeq w_{[1,i]}$ and $\bar{w}_j \defeq \bar{w}_{[1,j]}$ for all $i, j \geq 1$. Then $\chi_w = \chi_{\bar{w}}$ if and only if for all $i$, there exists $j$ with $\bar{w}_j \in w_i P$, and for all $j$, there exists $i$ with $w_i \in \bar{w}_j P$. 
\setlength{\parindent}{0.5cm} \setlength{\parskip}{0cm}

Assume that the infinite word $w$ does not contain any relator as a subword. Then $\chi_w = \chi_{\bar{w}}$ if and only if $w \equiv \bar{w}$.
\elemma
\setlength{\parindent}{0cm} \setlength{\parskip}{0cm}

\bproof
The first part follows immediately from the construction of characters from infinite words. For the second part, $\chi_w = \chi_{\bar{w}}$ implies that given $j$, there exists $i$ satisfying without loss of generality that $i \geq j$ such that $w_i = \bar{w}_j x$ for some $x \in \Sigma^*$. Since $w_i$ contains no relator as a subword, we must have $w_i \equiv \bar{w}_j x$, and hence $w_j \equiv \bar{w}_j$.
\eproof
\setlength{\parindent}{0cm} \setlength{\parskip}{0.5cm}

\section{The non-reversible case}
\label{sec:nonOre}

Throughout this section, let $P$ be a right LCM monoid given by a presentation $(\Sigma, R)$ as in \S~\ref{sss:MonoidsPres} satisfying our standing assumptions in \S~\ref{sss:Standing}. Let $G$ be the group given by the same presentation as $P$, and suppose that $P \to G$ induced by $\id_{\Sigma}$ is an embedding.

Let us start with the following general result which tells us that if our presentation does not have too many relations, then typically the boundary $\partial \Omega$ is as large as it can be.

\btheo
\label{THM1:bd=}
Suppose that $\abs{\Sigma} \geq 2$. Assume that there is $\varepsilon \neq z \in \Sigma^*$ such that $z$ is not a subword of a relator, no relator is a subword of $z$, no prefix of $z$ is a suffix of a relator, and no suffix of $z$ is a prefix of a relator.
\setlength{\parindent}{0cm} \setlength{\parskip}{0cm}

\begin{enumerate}
\item[(i)] If there exist $a_1, \dotsc, a_n \in P$ such that $P \setminus \gekl{e} = \bigcup_{i=1}^n a_iP$, then we have $\overline{G.\chi} = \Omega_{\infty}$ for every $\chi \in \Omega_{\infty}$.
\item[(ii)] If we cannot write $P \setminus \gekl{e} = \bigcup_{i=1}^n a_iP$ for some $a_1, \dotsc, a_n \in P$, then we have $\overline{G.\chi} = \Omega$ for every $\chi \in \Omega_{\infty}$.
\end{enumerate}
Moreover, in either case there exists $\omega \in \Omega_{\infty}$ with trivial stabilizer group, $G_{\omega} = \gekl{e}$.
\etheo 
\setlength{\parindent}{0cm} \setlength{\parskip}{0cm}

\bproof
Consider a basic open subset
$$U = U(pP; p_1 P, \dotsc, p_k P) = \menge{\chi \in \Omega}{\chi(pP) = 1, \, \chi(p_1 P) = \dotso = \chi(p_k P) = 0},$$
with $p_j P \subsetneq pP$ for all $1 \leq j \leq k$. In case (i), we may assume that $\Omega_{\infty} \cap U \neq \emptyset$ and need to find $g \in G$ with $g.\chi \in U$ ($g.\chi$ will then automatically lie in $\Omega_{\infty} \cap U$ because $\chi \in \Omega_{\infty}$ and $\Omega_{\infty}$ is $G$-invariant). In case (ii), we only know that $U \neq \emptyset$, and we need to find $g \in G$ with $g.\chi \in U$.
\setlength{\parindent}{0cm} \setlength{\parskip}{0.5cm}

Let us first of all show that in case (ii), we actually also have $\Omega_{\infty} \cap U \neq \emptyset$. Choose generators $\sigma_j \in \Sigma$ such that $p_j P \subseteq p \sigma_j P$. As $P \setminus \gekl{e} \supsetneq \bigcup_j \sigma_jP$, there exists $\sigma \in \Sigma$ with $\sigma \notin \bigcup_j \sigma_j P$. Consider the infinite word $w = p \sigma z z z \dotsm$, and the corresponding element $\chi_w$. Obviously, $\chi_w(pP) = 1$. If $\chi_w(p_j P) = 1$, then $p \sigma z^m \in p_j P$ for some $m$, which implies that $\sigma z^m \in \sigma_j P$. So there exists $x \in \Sigma^*$ with $\sigma z^m = \sigma_j x$. Thus there exists a sequence $w_1, \dotsc, w_s \in \Sigma^*$ such that $w_1 \equiv \sigma z^m$, $w_s \equiv \sigma_j x$, and for all $1 \leq i \leq s-1$, we have $w_i \equiv_R w_{i+1}$. By our assumptions on $z$, we must be able to write $w_r \equiv v_r z^m$ with $\ell^*(v_r) \geq 1$ and $\sigma = v_r$ for all $1 \leq r \leq s$. Hence $v_s z^m \equiv \sigma_j x$, which implies that the first letter of $v_s$ must be $\sigma_j$, so that $v_s \in \sigma_j P$. However, we know that $\sigma = v_s$. This is a contradiction as $\sigma \notin \bigcup_j \sigma_j P$ by our choice of $\sigma$.

Now choose $\check{p}_j \in \Sigma^*$ with $\check{p}_j = p_j$ in $P$. Let $\check{\Sigma} \subseteq \Sigma$ be a finite subset such that $\check{p}_j \in \check{\Sigma}^*$ for all $j$. Define $\check{\ell}^*$ to be the word length in $\check{\Sigma}^*$, and for $x \in P$, set $\check{\ell}(x) \defeq \inf \menge{L \in \Zz_{\geq 0}}{x = \sigma_1 \dotsm \sigma_L, \, \sigma_1, \dotsc, \sigma_L \in \check{\Sigma}}$. Note that if $x$ cannot be represented by a word in $\check{\Sigma}^*$, then $\check{\ell}(x) = \infty$. We claim that in both cases (i) and (ii), we can always find $x \in pP$ with $x \notin p_j P$ for all $j$ and $\check{\ell}(x) \geq \check{\ell}^*(\check{p}_j) = \ell^*(\check{p}_j)$ for all $j$. In case (ii), this follows from our argument above and because $\lim_{m \to \infty} \ell(p \sigma z^m) = \infty$. In case (i), we can find $t \in \Sigma^{\infty}$ with $\chi_t \in U$ because $\Omega_{\infty} \cap U \neq \emptyset$. Let $t_i \defeq t_{[1,i]}$ be the word consisting of the first $i$ letters of $t$. Then $t_1 P \supsetneq t_2 P \supsetneq t_3 P \supsetneq \dotso$, which implies that $\sup_i \check{\ell}(t_i) = \infty$. Otherwise $\menge{t_i}{i = 1, 2, \dotsc}$ would be a finite set since $\check{\Sigma}$ is finite, but then we could not have $t_1 P \supsetneq t_2 P \supsetneq t_3 P \supsetneq \dotso$. Now choose $i$ with $t_i \in pP$ and $\check{\ell}(t_i) \geq \max_j \check{\ell}^*(\check{p}_j)$, and set $x \defeq t_i$. 

So let $x \in pP$ satisfy $x \notin p_j P$ for all $j$ and $\check{\ell}(x) \geq \check{\ell}^*(\check{p}_j)$ for all $j$. Set $g = xz$. We have to show that $g.\chi \in U$. Take $s \in \Sigma^{\infty}$ with $\chi = \chi_s$, and let $s_i \defeq s_{[1,i]}$ be the word consisting of the first $i$ letters of $s$. Then $g.\chi = g.\chi_s = \chi_{xzs}$. Since $x \in pP$, we have $(g.\chi)(pP) = 1$. Assume that $(g.\chi)(p_j P) = 1$ for some $j$. Then there exists $n$ with $x z s_n = \check{p}_j r$ for some $r \in \Sigma^*$. This means that we can find $y_1, \dotsc, y_m \in \Sigma^*$ with $y_1 \equiv x z s_n$, $y_m \equiv \check{p}_j r$ and for all $1 \leq l \leq m-1$, there exist $\alpha, u, v, \zeta \in \Sigma^*$ with $(u,v) \in R$ or $(v,u) \in R$ or $(u,v) = (\varepsilon, \varepsilon)$ such that $y_l \equiv \alpha u \zeta$ and $y_{l+1} \equiv \alpha v \zeta$. By our assumptions on $z$, the subword $z$ of $x z s_n$ will not be changed. Hence for every $j$, we must have $y_l \equiv x_l z s_{n,l}$, with $x_l = x$ and $s_{n,l} = s_n$ in $P$. In particular, we have $x_m z s_{n,m} \equiv \check{p}_j r$. Now we claim that $\ell^*(x_m) \geq \ell^*(\check{p}_j)$. If not, then $\ell^*(x_m) < \ell^*(\check{p}_j)$, so that $x_m \in \check{\Sigma}^*$, and thus $\ell^*(x_m) = \check{\ell}^*(x_m) \geq \check{\ell}(x) \geq \check{\ell}^*(\check{p}_j) = \ell^*(\check{p}_j)$, which is a contradiction. So we have $\ell^*(x_m) \geq \ell^*(\check{p}_j)$, which implies $x_m \equiv \check{p}_j q$ for some $q \in \Sigma^*$. Hence it follows that $x \in \check{p}_j P = p_j P$, which is a contradiction. This shows $g.\chi \in U$.

Now let us find $\omega \in \Omega_{\infty}$ with trivial stabilizer group. First assume that $z$ contains two different letters $\alpha$ and $\beta$. By assumption, $\alpha$ and $\beta$ are not relators. Choose a sequence $\sigma_1, \sigma_2, \dotsc$ in $\gekl{\alpha,\beta}$ such that the infinite word $\sigma_1 \sigma_2 \dotsm$ has no period. Now form the infinite word $\bar{w} \defeq z \sigma_1 z \sigma_2 z \dotsm$. It has no period, either, because we cannot shift by a multiple of $\ell^*(z) + 1$ as $\sigma_1 \sigma_2 \dotsm$ has no period, and if shifting $\bar{w}$ by a different number of letters produces $\bar{w}$ again, then all the $\sigma_i$ would have to coincide with a fixed letter of $z$, which again is not possible as $\sigma_1 \sigma_2 \dotsm$ has no period. By our assumption on $z$, $\bar{w}$ does not contain a relator as a finite subword. Now assume that $h \in G$ satisfies $h.\chi_{\bar{w}} = \chi_{\bar{w}}$. Then we must have $h = pq^{-1}$, and $\chi_{\bar{w}}(qP) = 1$. The latter implies that $z \sigma_1 z \sigma_2 \dotsm z \sigma_n \in qP$ for some $n$, so that $z \sigma_1 z \sigma_2 \dotsm z \sigma_n = qr$ in $P$. But since $z \sigma_1 z \sigma_2 \dotsm z \sigma_n$ contains no relator as a subword, we must have $z \sigma_1 z \sigma_2 \dotsm z \sigma_n \equiv qr$. In particular, $\bar{w} \equiv q \ti{w}$ for some $\ti{w} \in \Sigma^{\infty}$, and $\ti{w}$ also contains no relator as a finite subword and has no period. Now we have $\chi_{\bar{w}} = h.\chi_{\bar{w}} = h.\chi_{q \ti{w}} = \chi_{p \ti{w}}$. Since $\bar{w}$ contains no relator as a finite subword, Lemma~\ref{LEM:w=w} implies that $q \ti{w} \equiv \bar{w} \equiv p \ti{w}$. Since $\ti{w}$ has no period, we must have $\ell^*(p) = \ell^*(q)$, which implies $p \equiv q$. Hence $h = pq^{-1} = e$.

Now assume that $z$ is a word in a single letter. Then, as we assume $\abs{\Sigma} \geq 2$, there must be another letter $\sigma \in \Sigma$ which does not appear in $z$. Form the infinite word $\bar{w} \defeq z \sigma z^2 \sigma z^3 \sigma \dotsm$. Assume that $h \in G$ satisfies $h.\chi_{\bar{w}} = \chi_{\bar{w}}$. Then we must have $h = pq^{-1}$, and $\chi_{\bar{w}}(qP) = 1$. Then there must exist $n \in \Zz_{\geq 1}$ such that $z \sigma z^2 \sigma \dotsm z^n \sigma = qr$ for some $r \in P$. Our assumptions on $z$ imply that $q = z \sigma z^2 \sigma \dotsm \sigma \zeta$ for some subword $\zeta$ of a power of $z$. So $\bar{w} = q \ti{w}$ with $\ti{w} = \eta \sigma z^k \sigma z^{k+1} \sigma \dotsm$, where $\eta$ is a subword of some power of $z$. Now write $p = x_1 a_1 \dotsm x_m a_m$, where $x_i$ are subwords of powers of $z$ and $a_i$ are words not containing the letter appearing in $z$. As we have seen above, $q$ is of the same form, i.e., $q = y_1 b_1 \dotsm y_n b_n$, where $y_j$ are subwords of powers of $z$ and $b_j$ are words not containing the letter appearing in $z$. We have $p \ti{w} = x_1 a_1 \dotsm x_m a_m \eta \sigma z^k \sigma z^{k+1} \sigma \dotsm$ and $q \ti{w} = y_1 b_1 \dotsm y_n b_n \eta \sigma z^k \sigma z^{k+1} \sigma \dotsm$. Now our assumptions on $z$ imply that $m=n$ and $x_1 = y_1$, $a_1 = b_1$, ..., $x_m = y_m$ and $a_m = b_m$, and thus $p = q$, i.e., $h = p q^{-1} = e$.
\eproof

\setlength{\parindent}{0cm} \setlength{\parskip}{0cm}

The following is an immediate consequence:
\bcor
\label{Cor:bdq:pis_gen}
Let us keep the same assumptions as in Theorem~\ref{THM1:bd=}. If we are in case (i), i.e., there exist elements $a_1, \dotsc, a_n$ in our semigroup $P$ such that $P \setminus \gekl{e} = \bigcup_{i=1}^n a_iP$, then we have $\partial \Omega = \Omega \setminus P$ and, if $G$ is exact, a short exact sequence $0 \to \cK(\ell^2 P) \to C^*_{\lambda}(P) \to \partial C^*_{\lambda}(P) \to 0$.
\setlength{\parindent}{0.5cm} \setlength{\parskip}{0cm}

If we are in case (ii), i.e., we cannot write $P \setminus \gekl{e} = \bigcup_{i=1}^n a_iP$ for some $a_1, \dotsc, a_n \in P$, then we have $\partial \Omega = \Omega$ and $C^*_{\lambda}(P) = \partial C^*_{\lambda}(P)$. 

In both cases $\partial C^*_{\lambda}(P)$ is a purely infinite simple C*-algebra (we do not need exactness of $G$ for this).
\ecor
\setlength{\parindent}{0cm} \setlength{\parskip}{0cm}

Here we used \cite[Theorem~22.9]{Ex}. The last part of the corollary follows from Theorem~\ref{THM1:bd=} together with \cite[Corollary~5.7.17]{CELY}.
\setlength{\parindent}{0cm} \setlength{\parskip}{0.5cm}

\bremark
Let us assume that our presentation $(\Sigma,R)$ is $r$-homogeneous in the sense of \cite[Definition~4.1]{Deh03}, or equivalently (see \cite[Proposition~4.3]{Deh03}), that it is $r$-Noetherian in the sense of \cite[Definition~2.6]{Deh03}. Then there exist elements $a_1, \dotsc, a_n$ in our semigroup $P = \spkl{\Sigma \, \vert \, R}^+$ such that $P \setminus \gekl{e} = \bigcup_{i=1}^n a_iP$ if and only if $\abs{\Sigma} < \infty$. The \an{if}-part is clear. For the \an{only if}-part, assume $\abs{\Sigma} = \infty$ and that $P \setminus \gekl{e} = \bigcup_{i=1}^n a_iP$ for some $a_i \in P$. Then $P$ does not coincide with the monoid $\spkl{a_1, \dotsc, a_n}^+$ generated by $a_1, \dotsc, a_n$ (see \S~\ref{sss:Standing}). Take $x \in P \setminus \spkl{a_1, \dotsc, a_n}^+$. Define $x_0 \defeq x$, and given $x_k \in P \setminus \spkl{a_1, \dotsc, a_n}^+$, we have $x_k \in \bigcup_{i=1}^n a_i P$ so that we can find $1 \leq i \leq n$ and $x_{k+1} \in P$ such that $x_k = a_i x_{k+1}$. Since $x_k \in P \setminus \spkl{a_1, \dotsc, a_n}^+$, we must have $x_{k+1} \in P \setminus \spkl{a_1, \dotsc, a_n}^+$. So we obtain a sequence $(x_k)$ such that $Px_0 \subsetneq Px_1 \subsetneq Px_2 \subsetneq \dotso$. This contradicts our assumption that $(\Sigma,R)$ is $r$-Noetherian.
\eremark

Let us now apply our findings to the one-relator case. From now on until the rest of this section, let $P = \spkl{\Sigma \, \vert \, u = v}^+$ be a right LCM one-relator monoid, with $(\Sigma, u = v)$ satisfying our standing assumptions in \S~\ref{sss:Standing} and the assumptions in \S~\ref{sss:OneRelM}. 
\blemma
\label{Lem:OneRel-zex}
Let $P = \spkl{\Sigma \, \vert \, u = v}^+$ with $\abs{\Sigma} \geq 3$. Then there exists $\varepsilon \neq z \in \Sigma^*$ such that $z$ is not a subword of a relator, no relator is a subword of $z$, no prefix of $z$ is a suffix of a relator, and no suffix of $z$ is a prefix of a relator.
\elemma
\setlength{\parindent}{0cm} \setlength{\parskip}{0cm}

\bproof
Let $c \in \Sigma$ be a generator which is not the first letter of $u$ or $v$, and let $a \in \Sigma$ be a generator which is not the last letter of $u$ or $v$. Consider $a^i c^j$ for $i, j > \max(\ell^*(u), \ell^*(v))$. Then no prefix of $a^i c^j$ is a suffix of a relator, and no suffix of $a^i c^j$ is a prefix of a relator. Moreover, $a^i c^j$ is not a subword of a relator. The only problem is that $a^i c^j$ could contain a relator as a subword. Certainly, $a^i c^j$ can contain at most one relator as a subword, say $u$. Then $u \equiv a^k c^l$ with $k, l \geq 1$. By our assumptions, this means that $v$ does not start with $a$ and does not end with $c$. Choose $b \in \Sigma$ with $b \not\equiv a$, $b \not\equiv c$. Then $a^i b c^j$ has all the desired properties.
\eproof

\bcor
\label{Cor:OneRelM-Omega=P+bd}
Let $P = \spkl{\Sigma \, \vert \, u = v}^+$ with $\abs{\Sigma} \geq 3$. Then if $\abs{\Sigma} < \infty$, we have $\partial \Omega = \Omega \setminus P$ and a short exact sequence $0 \to \cK(\ell^2 P) \to C^*_{\lambda}(P) \to \partial C^*_{\lambda}(P) \to 0$. If $\abs{\Sigma} = \infty$, we have $\partial \Omega = \Omega$ and $C^*_{\lambda}(P) = \partial C^*_{\lambda}(P)$. In both cases, $\partial C^*_{\lambda}(P)$ is a purely infinite simple C*-algebra.
\ecor
\bproof
This follows from Lemma~\ref{Lem:OneRel-zex} and Corollary~\ref{Cor:bdq:pis_gen} once we show that for $P = \spkl{\Sigma \, \vert \, u = v}^+$, there exist elements $a_1, \dotsc, a_n \in P$ such that $P \setminus \gekl{e} = \bigcup_{i=1}^n a_iP$ if and only if $\abs{\Sigma} < \infty$. The \an{if}-part is clear. For the \an{only if}-part, assume $\abs{\Sigma} = \infty$. Choose $\check{a}_i \in \Sigma^*$ with $a_i = \check{a}_i$, let $\Sigma_A$ be the set of letters appearing in one of the $\check{a}_i$, and let $\Sigma_R$ be the set of letters which appear in $u$ or $v$. Take $\tau \in \Sigma \setminus (\Sigma_A \cup \Sigma_R)$. Then $\tau \neq e$ and $\tau \notin \bigcup_{i=1}^n a_iP$ since otherwise, $\tau = \check{a}_i x$ (for some $i$ and some $x \in P$) would imply that we can find a sequence $w_1, \dotsc, w_s$ in $\Sigma^*$ with $w_1 \equiv \tau$ and $w_s \equiv \check{a}_i x$ such that $w_r \equiv_R w_{r+1}$ for all $1 \leq r \leq s-1$. But $\tau \equiv w_1 \equiv_R w_2$ implies $w_2 \equiv \tau$, and we conclude inductively that $w_s \equiv \tau$, hence $\check{a}_i x \equiv \tau$. This, however, is impossible. Hence $P \setminus \gekl{e} \neq \bigcup_{i=1}^n a_iP$.

Note that exactness of the one-relator group $\spkl{\Sigma \, \vert \, u = v}$ follows from \cite{Gue}.
\eproof
\setlength{\parindent}{0cm} \setlength{\parskip}{0.5cm}

\bremark
\label{Rem:Toe=qlo-K}
Let $P = \spkl{\Sigma \, \vert \, u = v}^+$ with $\abs{\Sigma} \geq 3$, and let $G = \spkl{\Sigma \, \vert \, u = v}$ be the group given by the same presentation as $P$. $P \subseteq G$ is Toeplitz in the sense of \cite[\S~5.8]{CELY} if and only if $P \subseteq G$ is quasi-lattice ordered in Nica's sense \cite{Ni}, i.e., for every $g \in G$, $P \cap gP = \emptyset$ or there exists $p \in P$ with $P \cap gP = pP$.
\setlength{\parindent}{0.5cm} \setlength{\parskip}{0cm}

If $P \subseteq G$ is Toeplitz, then it follows from \cite[Theorem~5.2]{CEL2} that $K_*(C^*_{\lambda}(P)) \cong K_*(\Cz)$, and more precisely,
$$K_*(C^*_{\lambda}(P)) = (\Zz[1],\gekl{0}).$$
Here we use that one-relator groups satisfy the Baum-Connes conjecture with coefficients \cite{BBV}. If in addition $3 \leq \abs{\Sigma} < \infty$, the short exact sequence $0 \to \cK(\ell^2 P) \to C^*_{\lambda}(P) \to \partial C^*_{\lambda}(P) \to 0$ induces in K-theory the exact sequence 
$$
  0 \to K_1(\partial C^*_{\lambda}(P)) \to \Zz \to \Zz \to K_0(\partial C^*_{\lambda}(P)) \to 0,
$$
where the homomorphism $\Zz \to \Zz$ is given by multiplication with $1 - \abs{\Sigma} + 1 = 2 - \abs{\Sigma}$ since the minimal projection $1_{\gekl{e}} \in \cK(\ell^2 P)$ is given by $1_{\gekl{e}} = 1_{P} - \sum_{\sigma \in \Sigma} 1_{\sigma P} + 1_{uP}$, and all the projections $1_{P}$, $1_{\sigma P}$ and $1_{uP}$ are Murray-von Neumann equivalent to $1$ in $C^*_{\lambda}(P)$. Hence
$$(K_0(\partial C^*_{\lambda}(P)), [1], K_1(\partial C^*_{\lambda}(P))) \cong (\Zz / (\abs{\Sigma} - 2) \Zz, 1, \gekl{0}).$$
If $\abs{\Sigma} = \infty$, we have $C^*_{\lambda}(P) = \partial C^*_{\lambda}(P)$, so that
$$(K_0(\partial C^*_{\lambda}(P)), [1], K_1(\partial C^*_{\lambda}(P))) \cong (\Zz, 1, \gekl{0}).$$
\eremark

Let us discuss another natural topic, namely nuclearity of $\partial C^*_{\lambda}(P)$, where $P = \spkl{\Sigma \, \vert \, u = v}^+$, with $\abs{\Sigma} \geq 3$. Note that there are examples of such monoids for which the boundary quotient is not nuclear. Take for instance $P = \spkl{a, b, c \, \vert \, aba = bab}^+$. Then $P = B_3^+ * \spkl{c}^+$, where $B_3^+$ is the Braid monoid $B_3^+ = \spkl{a, b \, \vert \, aba = bab}^+$ and $\spkl{c}^+$ is the free monoid on one generator $c$. The Braid group $B_3 = \spkl{a, b \, \vert \, aba = bab}$ appears as a stabilizer group for the partial action $G \curvearrowright \partial \Omega$, where $G = \spkl{a, b, c \, \vert \, aba = bab}$. Hence $G \ltimes \partial \Omega$ is not amenable, and $\partial C^*_{\lambda}(P)$ is not nuclear.

Let us now present conditions on the presentation $(\Sigma, u=v)$ which ensure that $\partial C^*_{\lambda}(P)$ is nuclear. For this purpose, we introduce the following
\bdefin
\label{DEF:GraphModel}
Let $X \subseteq \Omega$ be a closed subset. The partial action $G \curvearrowright \Omega$ restricts to a partial dynamical system $G \curvearrowright X$, so that we can form the transformation groupoid $G \ltimes X$. Let $\cE$ be a graph with vertices $\cE^0$ and edges $\cE^1$. Let $\bm{\sigma}: \: \cE^1 \to \Sigma$ be a map. $\bm{\sigma}$ extends to maps on the finite path space $\cE^* \to \Sigma^*, \, \mu_1 \dotsm \mu_n \ma \bm{\sigma}(\mu_1) \dotsm \bm{\sigma}(\mu_n)$ and the infinite path space $\cE^{\infty} \to \Sigma^{\infty}, \, \mu_1 \mu_2 \dotsm \ma \bm{\sigma}(\mu_1) \bm{\sigma}(\mu_2) \dotsm $, both of which are again denoted by $\bm{\sigma}$. Consider the composition $\varphi: \: \partial \cE \to \Sigma^* \cup \Sigma^{\infty} \to \Omega$, where the first map is the restriction of $\bm{\sigma}$ to $\partial \cE \subseteq \cE^* \cup \cE^{\infty}$, and the second map is given by $\Sigma^* \cup \Sigma^{\infty} \to \Omega, \, w \ma \chi_w$. $(\cE,\bm{\sigma})$ is called a graph model for $G \curvearrowright X$ if $\Phi: \: \cG_{\cE} \to G \ltimes X, \, (\lambda \zeta, l(\lambda) - l(\zeta), \mu \zeta) \ma (\bm{\sigma}(\lambda) \bm{\sigma}(\mu)^{-1}, \varphi(\mu \zeta))$ is an isomorphism of topological groupoids. Here $l(\cdot)$ denotes the path length in $\cE^*$.
\edefin
\setlength{\parindent}{0cm} \setlength{\parskip}{0cm}

We refer the reader to \cite[\S~3.2]{Li17} and the references therein for more information on how $\partial \cE$ and $\cG_{\cE}$ are constructed. Note that we use the same convention as in \cite{Rae,Li17}, which is different from the one in \cite{Sp02}. Moreover, the character $\chi_w$ attached to a word $w$ has been constructed in \S~\ref{sss:DistSubspOmega} and \S~\ref{sss:Char_w}.

\bremark
Once $(\cE,\bm{\sigma})$ has been constructed, an important step in showing that $(\cE,\bm{\sigma})$ is a graph model is to show that $\varphi$ is a homeomorphism onto $X$. Note that, however, we cannot expect $\bm{\sigma}$ to be injective on the set $\cE^*$ of finite paths. The reason is that $\bm{\sigma}$ will not remember the source of our finite paths. This problem will disappear once we go over to $\partial \cE$, as $\partial \cE$ only contains those finite paths in $\cE^*$ whose source is a vertex $z$ with $\abs{r^{-1}(z)} = \infty$ or $0$ (where $r(\cdot)$ is the range map of $\cE$). The point is that in our graph models, there will be at most one vertex $z$ with $\abs{r^{-1}(z)} = \infty$, and for paths in $\cE^{\infty}$, we can reconstruct sources by looking at subsequent parts of our infinite path.
\eremark

Now let $P = \spkl{\Sigma \, \vert \, u = v}^+$ and $G = \spkl{\Sigma \, \vert \, u = v}$. Assume that
\begin{eqnarray}
\label{OVL}
&& {\rm OVL}(v) = \menge{x \in \Sigma^*}{v \equiv xy \equiv wx \ {\rm for} \ {\rm some} \ \varepsilon \neq w, y \in \Sigma^*} = \gekl{\varepsilon},\\
\label{l<l}
 && \ell_{\sigma}(u) < \ell_{\sigma}(v) \ {\rm for} \ {\rm some} \ \sigma \in \Sigma,\\
\label{uNotvNotu}
&& u \text{ is not a subword of } v, \text{ no prefix of } u \text{ is a suffix of } v, \text{ and no suffix of } u \text{ is a prefix of } v.
\end{eqnarray}
We now construct graph models for $G \curvearrowright \partial \Omega$. 
\btheo
\label{THM:GraphModel-nonOre}
Assume that \eqref{OVL}, \eqref{l<l} and \eqref{uNotvNotu} are valid. Set $L \defeq \ell^*(v) - 1$.
\begin{enumerate}
\item[(i)] Suppose that $3 \leq \abs{\Sigma} < \infty$. Define a graph as follows: Set $\cE^0 \defeq \menge{x \in \Sigma^*}{\ell^*(x) = L}$ as the set of vertices, and let $\cE^1 \defeq \menge{(y,x)}{y, x \in \cE^0; \; \exists \, \sigma, \tau \in \Sigma: \: y \tau \equiv \sigma x \not\equiv v}$. Define range and source maps by $r(y,x) = y$ and $s(y,x) = x$. Moreover, define $\bm{\sigma}: \: \cE^1 \to \Sigma$ by setting $\bm{\sigma}(y,x)$ as the first letter of $y$. Then $(\cE,\bm{\sigma})$ is a graph model for $G \curvearrowright \partial \Omega$.
\item[(ii)] Suppose that $\abs{\Sigma} = \infty$. Define a graph as follows: Let $\Sigma_R$ be the finite subset of $\Sigma$ containing all the letters that appear in the relators $u$ and $v$. Set $\cE^0 \defeq \menge{x \in \Sigma_R^*}{\ell^*(x) \leq L}$ as the set of vertices, and let 
\begin{align*}
  \cE^1 \defeq &\menge{(y,x)}{y, x \in \cE^0; \; \ell^*(y) = \ell^*(x) = L; \; \exists \, \sigma, \tau \in \Sigma_R: \: y \tau \equiv \sigma x \not\equiv v}\\
  \cup &\menge{(y,x)}{y, x \in \cE^0; \; \ell^*(y) = \ell^*(x) + 1; \; \exists \, \sigma \in \Sigma_R: \: y \equiv \sigma x}\\
  \cup &\menge{\mu_{\tau,x}}{x \in \cE^0, \, \tau \in \Sigma \setminus \Sigma_R}
\end{align*}
Define range and source maps by $r(y,x) = y$, $s(y,x) = x$, and $r(\mu_{\tau,x}) = \ve$, $s(\mu_{\tau,x}) = x$. Moreover, define $\bm{\sigma}: \: \cE^1 \to \Sigma$ by setting $\bm{\sigma}(y,x)$ as the first letter of $y$, and $\bm{\sigma}(\mu_{\tau,x}) = \tau$. Then $(\cE,\bm{\sigma})$ is a graph model for $G \curvearrowright \Omega$.
\end{enumerate}
\etheo
\setlength{\parindent}{0cm} \setlength{\parskip}{0cm}

Note that $3 \leq \abs{\Sigma}$ ensures that our monoid is not left reversible (see \S~\ref{sss:Rev}).

\bproof
(i): First of all, Corollary~\ref{Cor:OneRelM-Omega=P+bd} implies that in this case, we have $\partial \Omega = \Omega_{\infty}$. 
As explained in \S~\ref{sss:Rev}, our conditions \eqref{OVL} and \eqref{l<l} imply that for every $z \in \Sigma^*$, we can find a unique irreducible word $x \in \Sigma^*$ (i.e., $x$ does not contain $v$ as a subword) such that $z = x$ in $P$, and two irreducible words $x, y \in \Sigma^*$ represent the same element in $P$ if and only if $x \equiv y$. Given $z \in \Sigma^*$, we write $\rho(z)$ for the unique irreducible word in $\Sigma^*$ with $\rho(z) = z$ in $P$. Given $\chi \in \Omega_{\infty}$, Lemma~\ref{LEM:chi=chi_w} implies that we can find $x \in \Sigma^{\infty}$ such that $\chi = \chi_x$. Then $x$ can be chosen to be of the form $x \equiv y_1 v y_2 v y_3 v \dotsm$, where $y_i \in \Sigma^*$ do not contain $v$ as a subword. Set $w \defeq y_1 u y_2 u y_3 u \dotsm$. Then $w$ does not contain $v$ as a finite subword by our assumptions on $u$ and $v$, and we have $\chi = \chi_x = \chi_w$ by Lemma~\ref{LEM:w=w} (since $x_i = w_i$ for infinitely many $i$, with $w_i = w_{[1,i]}$ being the word consisting of the first $i$ letters of $w$). Moreover, given infinite words $\ti{w}, \bar{w} \in \Sigma^{\infty}$ which do not contain $v$ as finite subword such that $\chi_{\ti{w}} = \chi_{\bar{w}}$, then Lemma~\ref{LEM:w=w} implies that for all $j$, there exists $i$ such that $\ti{w}_i \in \bar{w}_j P$, where $\ti{w}_i = \ti{w}_{[1,i]}$ and $\bar{w}_j = \bar{w}_{[1,j]}$ (the words consisting of the first $i$ or $j$ letters of $\ti{w}$ or $\bar{w}$, respectively). In particular, for $j \geq l \defeq \ell^*(v)$, because of our assumption on $u$ and $v$, we conclude from $\ti{w}_i \in \bar{w}_j P$ that $i \geq j-l$ and $\ti{w}_{j-l} \equiv \bar{w}_{j-l}$. Hence $\ti{w} \equiv \bar{w}$. This shows that we have a bijection 
$$
  \menge{w \in \Sigma^{\infty}}{w \text{ is } v \text{-free}} \cong \Omega_{\infty}, \, w \ma \chi_w.
$$
Here $w \in \Sigma^{\infty}$ is called $v$-free if $w$ does not contain $v$ as finite subword. To see that $\bm{\sigma}$ gives rise to a bijection $\cE^{\infty} \cong \menge{w \in \Sigma^{\infty}}{w  \text{ is } v \text{-free}}$, we construct an inverse as follows: Given a $v$-free $w \in \Sigma^{\infty}$, say $w = \sigma_1 \sigma_2 \sigma_3 \dotsm$ with $\sigma_i \in \Sigma$, let $x_i \defeq \sigma_i \dotsm \sigma_{i+L-1}$, set $\mu_i \defeq (x_i,x_{i+1})$, and define $\bm{\mu}(w) \defeq \mu_1 \mu_2 \mu_3 \dotsm$. In this way, we obtain a map $\bm{\mu}: \: \menge{w \in \Sigma^{\infty}}{w \text{ is } v \text{-free}} \to \cE^{\infty}$ which is the inverse of $\cE^{\infty} \to \menge{w \in \Sigma^{\infty}}{w \text{ is } v \text{-free}}, \mu \ma \bm{\sigma}(\mu)$. Therefore, $\varphi$ is bijective as it is a composition of two bijections. Moreover, it is easy to see that $\varphi$ is a homeomorphism.
\setlength{\parindent}{0.5cm} \setlength{\parskip}{0cm}

Clearly, $\Phi$ is a continuous groupoid homomorphism. To construct its inverse, define $\Psi: \: G \ltimes \partial \Omega \to \cG_{\cE}, \, (g,\chi_w) \ma (\varphi^{-1}(g.\chi_w), \lim_{n \to \infty} (\ell^*(\rho(g w_n)) - \ell^*(w_n)), \varphi^{-1}(\chi_w))$, where $w$ is a $v$-free word, and $w_n = w_{[1,n]}$ (the word consisting of the first $n$ letters of $w$). To show that $\Psi$ is a well-defined continuous groupoid homomorphism, let us show that $\ell^*(\rho(g w_n)) - \ell^*(w_n)$ is well-defined and eventually constant (i.e., it does not depend on $n$ as long as $n$ is big enough). First of all, $\chi_w \in U_{g^{-1}}$ implies that there exists $M \in \Zz_{\geq 0}$ such that $g w_M \in P$. Moreover, our assumptions imply that there exists an integer $N \geq M$ ($N = M + L$ works) such that for all integers $n \geq N$, we have $\rho(g w_n) = \rho(g w_N) w_{[N+1,n]}$, where $w_{[N+1,n]}$ is the subword of $w$ starting with the $(N+1)$-th letter and ending with the $n$-th letter. Thus we have for all $n \geq N$ that $\ell^*(\rho(g w_n)) - \ell^*(w_n) = \ell^*(\rho(g w_N) w_{[N+1,n]}) - \ell^*(w_N w_{[N+1,n]}) = \ell^*(\rho(g w_N)) - \ell^*(w_N)$. So, indeed, it does not depend on $n$. This shows that $\Psi$ is a well-defined continuous groupoid homomorphism. It is straightforward to see that $\Psi$ is the inverse of $\Phi$.
\setlength{\parindent}{0cm} \setlength{\parskip}{0.5cm}

(ii): We have $\Omega = P \cup \Omega_{\infty}$. As in (i), we have a bijection 
$$
  \menge{w \in \Sigma^* \cup \Sigma^{\infty}}{w \text{ is } v \text{-free}} \cong \Omega, \, w \ma \chi_w.
$$
We claim that $\bm{\sigma}$ gives rise to a bijection $\partial \cE \cong \menge{w \in \Sigma^* \cup \Sigma^{\infty}}{w \text{ is } v \text{-free}}$. First of all, note that $\partial \cE$ consists of the finite paths in $\cE$ which start in $\ve$ (i.e., which have source equal to $\ve$) and the infinite paths in $\cE$. We now set out to construct an inverse of $\bm{\sigma}$. Given a $v$-free $w \in \Sigma_R^* \cup\Sigma_R^{\infty}$, say $w = \sigma_1 \sigma_2 \sigma_3 \dotsm$ with $\sigma_i \in \Sigma_R$, let $x_i \defeq \sigma_i \dotsm \sigma_{\max(i+L-1, \, l(w))}$, set $\mu_i \defeq (x_i,x_{i+1})$, and define $\bm{\mu}(w) \defeq \mu_1 \mu_2 \mu_3 \dotsm$ ($\bm{\mu}(w) \defeq \mu_1 \mu_2 \mu_3 \dotsm \mu_{l(w)}$ if $\ell^*(w) < \infty$). A general $v$-free word $w \in \Sigma^* \cup \Sigma^{\infty}$ is of the form $w = w_1 \tau_1 w_2 \tau_2 \dotsm$ or $w = w_1 \tau_1 w_2 \tau_2 \dotsm z$, where $w_i$ are $v$-free words in $\Sigma_R^*$, $\tau_i \in \Sigma \setminus \Sigma_R$, and $z$ is a $v$-free word in $\Sigma_R^{\infty}$. In the first case, set $\bm{\mu}(w) = \bm{\mu}(w_1) \mu_{\tau_1, \, r(\bm{\mu}(w_2))} \bm{\mu}(w_2) \mu_{\tau_2, \, r(\bm{\mu}(w_3))} \dotsm$, and in the second case, set $\bm{\mu}(w) = \bm{\mu}(w_1) \mu_{\tau_1, \, r(\bm{\mu}(w_2))} \bm{\mu}(w_2) \mu_{\tau_2, \, r(\bm{\mu}(w_3))} \dotsm \bm{\mu}(z)$. In this way, we obtain a map $\bm{\mu}: \: \menge{w \in \Sigma^* \cup \Sigma^{\infty}}{w \text{ is } v \text{-free}} \to \partial \cE$ which is the inverse of $\partial \cE \to \menge{w \in \Sigma^* \cup \Sigma^{\infty}}{w \text{ is } v \text{-free}}, \mu \ma \bm{\sigma}(\mu)$. Therefore, $\varphi$ is bijective as it is a composition of two bijections. Moreover, it is easy to see that $\varphi$ is a homeomorphism.
\setlength{\parindent}{0.5cm} \setlength{\parskip}{0cm}

Again, $\Phi$ is a continuous groupoid homomorphism, so we just need to construct its inverse $\Psi: \: G \ltimes \partial \Omega \to \cG_{\cE}$. For a $v$-free word $w \in \Sigma^{\infty}$, we define $\Psi(g,\chi_w)$ as in (i). For a $v$-free word $w \in \Sigma^*$, define $\Psi(g,\chi_w) \defeq (\varphi^{-1}(g.\chi_w), \ell^*(\rho(gw)) - \ell^*(w), \varphi^{-1}(\chi_w))$. Then $\Psi$ is a well-defined continuous groupoid homomorphism which is the inverse of $\Phi$.
\eproof
\setlength{\parindent}{0cm} \setlength{\parskip}{0.5cm}

Combining Theorem~\ref{THM:GraphModel-nonOre} with Corollary~\ref{Cor:OneRelM-Omega=P+bd}, we obtain
\bcor
\label{Cor:bd=graph_nonOre}
Let $P = \spkl{\Sigma \, \vert \, u = v}^+$ with $\abs{\Sigma} \geq 3$. Assume that \eqref{OVL}, \eqref{l<l} and \eqref{uNotvNotu} hold. Let $\cE$ be the graph model constructed in Theorem~\ref{THM:GraphModel-nonOre}. Then the boundary quotient of $C^*_{\lambda}(P)$ is isomorphic to the graph C*-algebra of $\cE$, $\partial C^*_{\lambda}(P) \cong C^*(\cE)$, and it is a UCT Kirchberg algebra. In particular, $C^*_{\lambda}(P)$ and $\partial C^*_{\lambda}(P)$ are nuclear, and $G \curvearrowright \Omega$ is amenable (as is $G \ltimes \Omega$).
\ecor

\bex
\label{ex:A*B}
Here is a concrete class of examples where Corollary~\ref{Cor:bd=graph_nonOre} applies: Let $A$ and $B$ be sets with $\abs{A} + \abs{B} \geq 3$. Choose $u \in A^*$ arbitrary and $v \in B^*$ with ${\rm OVL}(v) = \gekl{\varepsilon}$. Then the presentation $(A \cup B, u = v)$ satisfies all the assumptions in Theorem~\ref{THM:GraphModel-nonOre} and Corollary~\ref{Cor:bd=graph_nonOre}, so that Corollary~\ref{Cor:bd=graph_nonOre} applies to the monoid $P = \spkl{A \cup B \, \vert \, u = v}^+$. 
\eex

\section{The reversible case}
\label{sec:Ore}

Throughout this section, let $P$ be left reversible. In that case, $\partial \Omega$ reduces to a point, given by the map $\cJ \to \gekl{0,1}$ sending every non-empty set in $\cJ$ to $1$. We denote this point by $\infty$. Let us also assume that $P$ is right LCM throughout this section. Recall that $\Omega_{\infty} = \Omega \setminus P$. In the following, we write $\ti{\Omega} \defeq \Omega \setminus \gekl{\infty}$ and $\ti{\Omega}_{\infty} \defeq \Omega_{\infty} \setminus \gekl{\infty}$.

\subsection{Reversible one-relator monoids}
\label{ss:OreOneRel}

Let us present one-relator monoids $P = \spkl{\Sigma \, \vert \, u=v}^+$ which are left reversible and for which we can describe $G \ltimes \ti{\Omega}_{\infty}$. Let $(\Sigma, u = v)$ satisfy the assumptions in \S~\ref{sss:Standing} and \S~\ref{sss:OneRelM}.

As explained in \S~\ref{sss:Rev}, $P = \spkl{\Sigma \, \vert \, u=v}^+$ can only be left reversible if $\abs{\Sigma} \leq 2$. Since our assumptions force $\abs{\Sigma} \geq 2$, we must have $\abs{\Sigma} = 2$, say $\Sigma = \gekl{a,b}$. The next lemma gives a special condition which enforces left reversibility.
\blemma
\label{LEM:w}
If there is $w \in P$ with $w = a x = x x' = b y = y y'$ for some $x, y, x', y' \in P$, then we have
\setlength{\parindent}{0cm} \setlength{\parskip}{0cm}
\begin{itemize}
\item $pw \in wP$ for all $p \in P$;
\item for every $p \in P$, there exists $i \geq 0$ with $w^i \in pP$.
\end{itemize}
\elemma
\setlength{\parindent}{0cm} \setlength{\parskip}{0cm}

\bproof
To prove the first claim, take $p \in \gekl{a,b}^*$. We proceed inductively on $\ell^*(p)$. Our claim certainly holds for $p = \varepsilon$. If $pw \in wP$, then $paw = p a x x' = pw x' \in wP$. Similarly, $pbw \in wP$.
\setlength{\parindent}{0cm} \setlength{\parskip}{0.5cm}

To prove the second claim, again take $p \in \gekl{a,b}^*$ and proceed inductively on $\ell^*(p)$. Our claim certainly holds for $p = \varepsilon$. Assume $w^i \in pP$, say $w^i = pq$ for some $q \in P$. By (i), we know that $qw = wr$ for some $r \in P$. Hence $p a x r = p w r = p q w = w^{i+1}$. This shows that $w^{i+1} \in paP$. Similarly, $w^{i+1} \in pbP$.
\eproof
For instance, if $P$ is a Garside monoid in the sense of \cite{DP, Deh02, Deh15}, then we may take for $w$ a Garside element in $P$.
\setlength{\parindent}{0cm} \setlength{\parskip}{0.5cm}

The following is now an immediate consequence of Lemma~\ref{LEM:w}.
\bcor
If there is $w \in P$ with $w = a \alpha = \alpha \gamma = b \beta = \beta \delta$ for some $\alpha, \beta, \gamma, \delta \in P$, then $\chi \in \Omega$ satisfies $\chi(pP) = 1$ for all $p \in P$ if and only if $\chi(w^i P) = 1$ for all $i \geq 0$. Moreover, let $X_i \defeq \menge{\chi \in \Omega}{\chi(w^i P) = 1, \, \chi(w^{i+1} P) = 0}$. Then $\Omega$ admits the following decomposition as a (set-theoretic) disjoint union
$$
  \Omega = \Big( \coprod_{i=0}^{\infty} X_i \Big) \amalg \gekl{\infty}.
$$
\ecor
\setlength{\parindent}{0cm} \setlength{\parskip}{0cm}

We write $X \defeq X_0$. Our goal is to describe $G \ltimes \ti{\Omega}_{\infty}$. Since $Y \defeq \Omega_{\infty} \cap X$ is obviously a compact open subspace of $\Omega_{\infty}$ which meets every $G$-orbit of $\ti{\Omega}_{\infty}$, we have $G \ltimes \ti{\Omega}_{\infty} \sim_M (G \ltimes \ti{\Omega}_{\infty}) \, \vert \, Y$. Hence it suffices to describe the restriction $(G \ltimes \ti{\Omega}_{\infty}) \, \vert \, Y$ of $G \ltimes \ti{\Omega}_{\infty}$ to $Y$. Let us now -- in the same spirit as Theorem~\ref{THM:GraphModel-nonOre} -- present two cases where we can provide a graph model for $G \curvearrowright Y$. 

\btheo
\label{Thm:OreOneRel:graph}
$ $
\begin{enumerate}
\item[I.] Assume that in Lemma~\ref{LEM:w}, we may take $w = u = v$. Let $L \defeq \max(\ell^*(u), \ell^*(v)) - 1$. Define a graph as follows: As $\cE^0$, take the (finite) collection of all finite subwords of length $L$ in $u, v$-free infinite words. Let $\cE^1 \defeq \menge{(y,x)}{y, x \in \cE^0; \; \exists \, \sigma, \tau \in \gekl{a,b}: \: y \tau \equiv \sigma x \text{ is } u, v \text{-free}}$. Define range and source maps by $r(y,x) = y$ and $s(y,x) = x$. Moreover, define $\bm{\sigma}: \: \cE^1 \to \gekl{a,b}$ by setting $\bm{\sigma}(y,x)$ as the first letter of $y$. Then $(\cE,\bm{\sigma})$ is a graph model for $G \curvearrowright Y$.
\item[II.] Now we do not assume that $w$ in Lemma~\ref{LEM:w} satisfies $w = u = v$. Suppose the following hold:
\begin{enumerate}
\item[II.1.] Let $\cX \defeq \menge{z \in \gekl{a,b}^* \ v \text{-free}}{z \in wP}$. Assume that there exists a finite subset $\cW \subseteq \cX$ such that for all $z \in \cX$ there is $x \in \cW$ such that $x$ is a subword of $z$.
\item[II.2.] $u$ is not a subword of $v$, no prefix of $u$ is a suffix of $v$, no suffix of $u$ is a prefix of $v$, and $v$ is not a subword of $u$.
\item[II.3.] There exists $k \geq 0$ such that for all $z = ar \in \gekl{a,b}^*$ and $\bar{z} = bs \in \gekl{a,b}^*$ with no $v$ as a subword and with $\ell^*(r), \ell^*(s) \geq k$, we must have $zP \cap \bar{z}P \subseteq wP$.
\end{enumerate}
Let $L \defeq \max(\menge{\ell^*(x)}{x \in \gekl{v} \cup \cW}) - 1$. Define a graph as follows: As the vertex set $\cE^0$, take the (finite) collection of all finite subwords of length $L$ in $v, \cW$-free infinite words. For the edges, set $\cE^1 \defeq \menge{(y,x)}
{y, x \in \cE^0; \; \exists \, \sigma, \tau \in \gekl{a,b}: \: y \tau \equiv \sigma x \text{ is } v, \cW \text{-free}}$. Define range and source maps by $r(y,x) = y$ and $s(y,x) = x$. Moreover, define $\bm{\sigma}: \: \cE^1 \to \gekl{a,b}$ by setting $\bm{\sigma}(y,x)$ as the first letter of $y$. Then $(\cE,\bm{\sigma})$ is a graph model for $G \curvearrowright Y$.
\end{enumerate}
\etheo
Here a word (finite or infinite) is called $u,v$-free if it does not contain $u$ or $v$ as a finite subword, $v$-free if it does not contain $v$ as a finite subword, and $v, \cW$-free if it is $v$-free and contains no $x \in \cW$ as a finite subword.
\bproof
In case I., note that we have a bijection $\menge{z \in \gekl{a,b}^{\infty}}{z \text{ is } u,v \text{-free}} \cong Y, \, z \ma \chi_z$ by Lemma~\ref{LEM:w=w}. Moreover, the inverse for the map $\cE^{\infty} \to \menge{z \in \gekl{a,b}^{\infty}}{z \text{ is } u, v \text{-free}}, \, \mu \ma \bm{\sigma}(\mu)$ is given as follows: Given a $u,v$-free $z \in \gekl{a,b}^{\infty}$, say $z = \sigma_1 \sigma_2 \sigma_3 \dotsm$ with $\sigma_i \in \gekl{a,b}$, let $x_i \defeq \sigma_i \dotsm \sigma_{i+L-1}$, set $\mu_i \defeq (x_i,x_{i+1})$, and define $\bm{\mu}(w) \defeq \mu_1 \mu_2 \mu_3 \dotsm$. Hence the map $\varphi$ from Definition~\ref{DEF:GraphModel} -- being the composition of two bijections -- is bijective. It is easy to see that $\varphi$ is a homeomorphism. Finally, the inverse of the map $\Phi$ from Definition~\ref{DEF:GraphModel} is given by $\Psi: \: G \ltimes Y \to \cG_{\cE}, \, (g,\chi_z) \ma (\varphi^{-1}(g.\chi_z), \lim_{n \to \infty} (\ell^*(g z_n) - \ell^*(z_n)), \varphi^{-1}(\chi_z))$, where $z$ is a $u,v$-free word, and $z_n = z_{[1,n]}$ is the word consisting of the first $n$ letters of $z$.
\setlength{\parindent}{0cm} \setlength{\parskip}{0.5cm}

In case II., note that as in the proof of Theorem~\ref{THM:GraphModel-nonOre}~(i), given an infinite word, condition II.2 guarantees that we can always replace $v$ by $u$ to arrive at a $v$-free infinite word.
\setlength{\parindent}{0.5cm} \setlength{\parskip}{0cm}

Condition II.1 implies that for a $v$-free infinite word $z$, $\chi_z(wP) = 0$ if and only if $z$ is $\cW$-free: Suppose that $z$ is not $\cW$-free, i.e., it is of the form $z \equiv \alpha x \zeta$ with $x \in \cW$. As $x$ lies in $wP$, we can write $x = wy$ for some $y \in P$. Hence $\alpha x = \alpha w y \in wP$ by Lemma~\ref{LEM:w}. By construction of $\chi_z$ (see \S~\ref{sss:Char_w}), this implies that $\chi_z(wP) = 1$. Conversely, suppose that $\chi_z(wP) = 1$. Again by construction of $\chi_z$ (see \S~\ref{sss:Char_w}), there must exist a finite word $z'$ with $z' \in wP$ such that $z \equiv z'z''$. Since $z$ is $v$-free, so is $z'$. Hence $z'$ lies in $\cX$, and condition II.1 implies that there exists $x \in \cW$ such that $x$ is a subword of $z'$. Hence $x$ is a subword of $z$, so that $z$ is not $\cW$-free.

Condition II.3. ensures that two distinct $v$-free infinite words $x$ and $y$ determine distinguished characters $\chi_x \neq \chi_y \in Y$: Let $x, y$ be $v$-free infinite words with $x \not\equiv y$, say $x \equiv p a \dotsm$ and $y \equiv p b \dotsm$ for some $p \in \gekl{a,b}^*$. For some $i$ and $j$, we have $x_i \equiv x_{[1,i]} \equiv par$ and $y_j \equiv y_{[1,j]} \equiv pbs$ for some $r, s \in \gekl{a,b}^*$ with $\ell^*(r), \ell^*(s) \geq k$ ($x_{[1,i]}$ and $y_{[1,j]}$ are the words consisting of the first $i$ and $j$ letters of $x$ and $y$, respectively). Now suppose that $\chi_x = \chi_y \eqdef \chi$, and hence that $\chi(x_i P) = 1 = \chi(y_j P)$. Then $\chi((parP) \cap (pbsP)) = 1$. Condition II.3. implies that $(parP) \cap (pbsP) \subseteq pwP \subseteq wP$, where the last inclusion is due to Lemma~\ref{LEM:w}~(i). Hence $\chi(wP) = 1$, so that $\chi \notin Y$. 

Therefore, we have a bijection $\menge{z \in \gekl{a,b}^{\infty}}{z \text{ is } v, \cW \text{-free}} \cong Y, \, z \ma \chi_z$. As in case I., we deduce that $\varphi$ from Definition~\ref{DEF:GraphModel} is a homeomorphism. The proof that the map $\Phi$ from Definition~\ref{DEF:GraphModel} is an isomorphism of topological groupoids is now similar to the one of Theorem~\ref{THM:GraphModel-nonOre}~(i).
\eproof
\setlength{\parindent}{0cm} \setlength{\parskip}{0cm}

\bremark
Once a graph model has been found, it will follow that $(G \ltimes \ti{\Omega}_{\infty}) \, \vert \, Y$ and $G \ltimes \ti{\Omega}_{\infty}$ are amenable, and thus $C^*_r((G \ltimes \ti{\Omega}_{\infty}) \, \vert \, Y)$ and $C^*_r(G \ltimes \ti{\Omega}_{\infty})$ are nuclear (see for instance \cite[Remark~3.6]{Li17}).
\eremark
\setlength{\parindent}{0cm} \setlength{\parskip}{0cm}

\bex
\label{ex:ArtinDihedralTorusKnot}
Here are two classes of examples for case I: Consider $u = aba \dotsm$ and $v = bab \dotsm$ with $\ell^*(u) = \ell^*(v)$, or consider $u = a^p$ and $v = b^q$ for $p, q \geq 2$. The corresponding presentations $(\gekl{a,b}, u=v)$ both define Garside monoids, where $w = u = v$ is a Garside element. Thus the conditions in case I of Theorem~\ref{Thm:OreOneRel:graph} are satisfied. The groups defined by these presentations are the Artin-Tits groups of dihedral type and the torus knot groups. It is conjectured in \cite{Pi} that these are the only Garside groups on two generators (see also \cite[Question~1]{Deh15}).
\eex

\bremark
\label{R:pisC}
For all the presentations in Example~\ref{ex:ArtinDihedralTorusKnot} apart from $(\gekl{a,b}, ab=ba)$ and $(\gekl{a,b}, a^2 = b^2)$, $C^*_r((G \ltimes \ti{\Omega}_{\infty}) \, \vert \, Y)$ and $C^*_r(G \ltimes \ti{\Omega}_{\infty})$ are purely infinite simple. This is straightforward to check using, for instance, \cite[Theorem~3.15]{Sp02}. (Note that the conventions in \cite{Sp02} and \cite{Li17} are different, one has to reverse the arrows in the graphs in one of these papers to translate into the convention of the other.) Therefore, $C^*_r((G \ltimes \ti{\Omega}_{\infty}) \, \vert \, Y)$ and $C^*_r(G \ltimes \ti{\Omega}_{\infty})$ are UCT Kirchberg algebras. While $C^*_r((G \ltimes \ti{\Omega}_{\infty}) \, \vert \, Y)$ is unital, $C^*_r(G \ltimes \ti{\Omega}_{\infty})$ is not unital and hence must be stable by \cite{Zh}. As UCT Kirchberg algebras are classified by their K-theory, it is a natural task to compute K-theory for $C^*_r((G \ltimes \ti{\Omega}_{\infty}) \, \vert \, Y)$ and $C^*_r(G \ltimes \ti{\Omega}_{\infty})$. We will discuss K-theory in \S~\ref{ss:Ex-K_*}.
\eremark

In order to discuss some examples for case II, let us formulate a stronger condition which implies condition II.3. Assume that $u$ starts with $b$ and $v$ starts with $a$, i.e., $u = b \dotsm$ and $v = a \dotsm$. For $1 \leq l \leq \ell^*(v) - 1$, let $v_{[1,l]}$ be the word consisting of the first $l$ letters of $v$, and let $\sigma_{l+1} \in \gekl{a,b}$ be the generator which is not the $(l+1)$-th letter of $v$. Now consider the following condition:
\setlength{\parindent}{0cm} \setlength{\parskip}{0cm}

\begin{enumerate}
\item[II.3'.] For all $1 \leq l \leq \ell^*(v) - 1$, $v_{[1,l]} \sigma_{l+1} P \cap bP \subseteq wP$.
\end{enumerate}
Clearly, condition II.3'. implies II.3. (with $k = \ell^*(v)$).

\bex
\label{ex:u=bjv=...}
Let $u = b^j$ and $v = a b^i a b^i a \dotsm a b^i a$, with $i > 0$ and $\ell_b(u) > \ell_b(v)$. Then we saw in \S~\ref{sss:LCM} that $P = \spkl{a,b \, \vert \, u = v}^+$ is right LCM. It is easy to see that $w = b^{i+j}$ satisfies the conditions in Lemma~\ref{LEM:w}. Moreover, our assumptions on $u$ and $v$ in case II are satisfied: Clearly, condition II.2 is satisfied. 
\setlength{\parindent}{0.5cm} \setlength{\parskip}{0cm}

It is easy to check that $W \defeq \menge{x \in \gekl{a,b}^*}{x = w \ {\rm in} \ P} = \gekl{b^{i+j}} \cup \menge{b^h v b^{i-h}}{0 \leq h \leq i}$. The only $v$-free word in $W$ is given by $b^{i+j}$. Let us verify condition II.1 for $\cW \defeq \menge{x \in \gekl{a,b}^* \ v \text{-free}}{x = w \ {\rm in} \ P} = \gekl{b^{i+j}}$. Take $z \in \cX$. As explained in \S~\ref{sss:MonoidsPres}, $z \in wP$ means that there exist $z_1$, ..., $z_s$ in $\gekl{a,b}^*$ such that $z_1 \equiv z$, $z_{n+1}$ is obtained from $z_n$ by replacing (once) $u$ by $v$ or $v$ by $u$, and $z_s$ contains some $y \in W$ as a subword. Without loss of generality we may assume that $s$ is minimal. Let us show that $z$ must already contain $b^{i+j}$ as a subword. Let $r \in \gekl{1, \dotsc, s}$ be maximal such that for all $1 \leq m \leq r-1$, $z_{m+1}$ is obtained from $z_m$ by replacing (once) $u$ by $v$. Assume first that $r=s$. Condition II.2 implies that we must have $z \equiv \dotsm u \dotsm u \dotsm u \dotsm$, and $z_r \equiv z_s \equiv \dotsm v \dotsm v \dotsm v \dotsm$, i.e., $z_r \equiv z_s$ is obtained from $z$ by replacing finitely many copies of $u$ by $v$, and in each step, changes take place in a different part of $z$ (there can be no overlap between a newly created $v$ and the next $u$ which is replaced because of condition II.2). Assume that $r > 1$. Since $z_r \equiv z_s$ contains $y \in W$ as a subword, and because by minimality of $s$ we know that $z_{r-1}$ cannot contain some word in $W$ as a subword, $y$ must have an overlap with a newly created $v$, i.e., $z_r$ is of the form $\dotsm b^h v b^{i-h} \dotsm$ or $\dotsm a b^i v \dotsm$ or $\dotsm v b^i a \dotsm$, and $z_{r-1}$ is of the form $\dotsm b^h u b^{i-h} \dotsm$ or $\dotsm a b^i u \dotsm$ or $\dotsm u b^i a \dotsm$. In each of these cases, we see that $z_{r-1}$ already contains some word in $W$ as a subword, contradicting minimality of $s$. Hence we conclude $r=1$, i.e., already $z$ contains some word in $W$ as a subword. As $z$ is $v$-free, it follows that $z$ must contain $b^{i+j}$ as a subword, as claimed. Now assume $r < s$. Then $z_{r+1}$ is obtained from $z_r$ by replacing (once) $v$ by $u$. So $z_r$ contains $v$ as a subword. We may as well assume that this instance of $v$ was not created as a replacement $u \to v$ in going from $z_c$ to $z_{c+1}$ for some $c < r$. Since $z$ is $v$-free, part of this copy of $v$ must have been newly created when going from $z_q$ to $z_{q+1}$ by replacing (once) $u$ by $v$, for some $1 \leq q \leq r-1$. So we must have $z_q \equiv \dotsm u \zeta \dotsm$ and $z_{q+1} \equiv \dotsm v \zeta \dotsm \equiv \dotsm v' v'' \zeta \dotsm$ with $v \equiv v'v''$ and $v'' \zeta \equiv v$, where $\zeta$ is a proper subword of $v$. (Another possibility would be $z_q \equiv \dotsm \zeta u \dotsm$, which is treated in an analogous way.) As $\zeta \not\equiv v$, $v''$ cannot be empty, and hence $v''$ must end with $a$ since $v$ ends with $a$. Because $s$ is minimal, $\zeta$ cannot be empty. Therefore, $\zeta$ must start with $b^i$, and hence $b^{i+j}$ is a subword of $u \zeta$, which in turn is a subword of $z_q$, contradicting minimality of $s$. So all in all, we have proven that $z$ must contain $b^{i+j}$ as a subword, as desired. 

To verify condition II.3'., note that $v_{[1,l]} \sigma_{l+1}$ is either of the form $a b^i a b^i a \dotsm a b^{i+1}$ (with less letters $a$ appearing as in $v$) or $a b^i a b^i a \dotsm a b^{h} a$ for $h < i$ (with at most as many letters $a$ appearing as in $v$). In both cases, it is easy to check -- for instance with the help of right reversing -- that $v_{[1,l]} \sigma_{l+1} P \cap bP \subseteq b^{i+j} P$.
\eex
\setlength{\parindent}{0cm} \setlength{\parskip}{0.5cm}

\bremark
As in Remark~\ref{R:pisC}, it is straightforward to check that for the presentations in Example~\ref{ex:u=bjv=...}, $C^*_r((G \ltimes \ti{\Omega}_{\infty}) \, \vert \, Y)$ is a unital UCT Kirchberg algebra and $C^*_r(G \ltimes \ti{\Omega}_{\infty})$ is a stable UCT Kirchberg algebra.
\eremark

\bremark
There is an overlap between the two classes of groups in Example~\ref{ex:ArtinDihedralTorusKnot}, because the presentations $(\gekl{a,b}, aba \dotsm = bab \dotsm)$, where each of the relators has $m$ factors, with $m$ odd, and $(\gekl{a,b}, a^2=b^m)$ define isomorphic groups. Moreover, the groups given by the presentations in Example~\ref{ex:u=bjv=...} are not new; they turn out to be isomorphic to the torus knot groups in Example~\ref{ex:ArtinDihedralTorusKnot}.
\eremark

\subsection{Artin-Tits monoids of finite type}
\label{ss:ATM_FiniteType}

Let us now study Artin-Tits monoids of finite type with more than two generators (see \cite{BS,Del,Sal}). We will see that the results in \S~\ref{ss:OreOneRel} do not carry over. Indeed, $G \ltimes \ti{\Omega}_{\infty}$ is typically not amenable. Moreover, our structural results for $C^*_r(G \ltimes \ti{\Omega}_{\infty})$ and thus for the semigroup C*-algebras of Artin-Tits monoids of finite type complement nicely the results in \cite{CL02,CL07,ELR} on C*-algebras of right-angled Artin monoids.

The following refers to \cite{Mi}.  Thus $G \equiv G_S$ is the group with finite generating set $S$, and presentation $\langle S \mid (s,t)_{m_{st}} = (t,s)_{m_{st}}, \text{ for } s,t \in S \rangle$, where $(s,t)_m = ststs...$ with $m$ factors, and $M = (m_{st})_{s,t \in S}$ is a symmetric positive integer matrix with $m_{ss} = 1$ and $m_{st} > 1$ if $s \not= t$.  Moreover, if $W$ is the group with the same presentation as $G$ but with the additional relations $\{s^2 : s \in S \}$ then $W$ is assumed to be finite.  A labeled graph (called the \emph{Coxeter diagram} of $G$) is associated to $G$ as follows.  $S$ is the set of vertices of the graph, and there is an edge between $s$ and $t$ if $m_{st} > 2$.  The edge is labeled $m_{st}$ if $m_{st} > 3$. $G$ is \emph{irreducible} if this graph is connected.  (If $S = S_1 \sqcup S_2$ such that $m_{s_1s_2} = 2$ whenever $s_1 \in S_1$ and $s_2 \in S_2$ then $G_S \cong G_{S_1} \times G_{S_2}$.)  We let $P \equiv P_S$ be the monoid given by $\langle S \mid (s,t)_{m_{st}} = (t,s)_{m_{st}} \text{ for } s,t \in S \rangle^+$.  As in \cite{Mi} we let $p : G \to W$ be the obvious surjection, and $P_{\text{red}} = \{ x \in P : \ell(x) = \ell(p(x)) \}$, where $\ell$ is the length function.  For $s$, $t \in S$ we write $\Delta_{s,t} = (s,t)_{m_{st}}$. We use the following notation from \cite{Mi}:  for $g$, $h \in P$ we write $g \prec h$ if $h \in gP$ and $h \succ g$ if $h \in Pg$.  By \cite[Proposition 2.6]{Mi}, for any two elements $g$, $h \in P$ there is a unique maximal (for $\prec$) element $g \wedge h \in P$ satisfying $g \wedge h \prec g$, $h$, and as in \cite[Section 3]{Mi} also a unique minimal element $g \vee h \in P$ satisfying $g$, $h \prec g \vee h$.  It follows that $P$ is left reversible and right LCM.  For $g \in P$ we let $\cL(g) = \{s \in S : s \prec g \}$ and $\cR(g) = \{s \in S : g \succ s\}$.  For $g \in P_{\text{red}}$ and $s \in S$, if $s \not\in \cR(g)$ then $gs \in P_{\text{red}}$, and similarly on the left (see \cite[remarks before Proposition 1.5]{Mi}).  A finite sequence $(g_1, \ldots, g_k)$ of elements of $P_{\text{red}} \setminus \{1\}$ is called a \emph{normal form} if $\cR(g_i) \supseteq \cL(g_{i+1})$ for $1 \le i < k$.  It is shown in \cite[Corollaries 4.2, 4.3]{Mi} that for each $g \in P$ there is a unique normal form $(g_1, \ldots, g_k)$ such that $g = g_1 \cdots g_k$.  We use the notation $\nu(g)$ from \cite[p. 367]{Mi} for the \emph{decomposition length} of $g \in P$:  $\nu(g) = n$ if $g$ has normal form $(g_1, \ldots, g_n)$.  If $T \subseteq S$ we write $G_T$, respectively $P_T$, for the Artin-Tits group, respectively monoid, corresponding to $\{ m_{st} : s,t \in T \}$.  It is clear that the inclusion $T \subseteq S$ defines a homomorphism $G_T \to G_S$ carrying $P_T$ to $P_S$.
\blemma
\label{lem:subset of generators}
Let $T \subseteq S$ and let $g \in \langle T \rangle_{G_S}$.  Then $\cL(g)$, $\cR(g) \subseteq T$.
\elemma
\setlength{\parindent}{0cm} \setlength{\parskip}{0cm}
\bproof
Let $s \in S \setminus T$.  Suppose that $s \in \cL(g)$.  Write $g = t_1 t_2 \cdots t_m$ with $t_i \in T$.  Since $s$, $t_1 \in \cL(g)$ it follows from \cite[Proposition 2.5]{Mi} that $\Delta_{\{s,t_1\}} \prec g$.  Then since $t_1 s \prec \Delta_{s,t_1}$ we have $t_1 s \prec g$, and hence $s \prec t_2 \cdots t_m$.  Repeating this process we eventually obtain $s \prec t_m$, hence $s = t_m \in T$, contradicting our supposition.  The argument for $\cR(g)$ is similar.
\eproof
\bprop
\label{prop:subset of generators}
Let $T \subseteq S$.  Then the homomorphism $G_T \to G_S$ is injective.
\eprop
\setlength{\parindent}{0cm} \setlength{\parskip}{0cm}
\bproof
Let $g \in P_T$.  Then $g$ has a normal form $(g_1, \ldots, g_k)$ in $P_T$.  By Lemma \ref{lem:subset of generators} this is also a normal form in $P_S$.  By \cite[Corollary 4.3]{Mi} it follows that $G_T \to G_S$ is injective on $P_T$.  By \cite[Corollary 3.2]{Mi} it follows that $G_T \to G_S$ is injective.
\eproof
As in \cite[Section 3]{Mi}, there is a unique element $\Delta \in P_{\text{red}}$ of greatest length (which plays the role of $w$ in Lemma \ref{LEM:w}).  We will write $P_0 = P_{\text{red}} \setminus \{1,\Delta\}$.  Then for $1 \not= g \in P \setminus \Delta P$, $g$ has normal form $(g_1, \ldots, g_k)$ with $g_i \in P_0$.  We wish to give a normal form for infinite words analogous to the normal form described above for elements of $P$.  For this it is convenient to recall the characterization of $\Omega$ given in \cite[Section 7]{Sp14}.  For $g \in P$ let $[g] = \{ h \in P : h \prec g \}$.  A subset $x \subseteq P$ is \emph{hereditary} if $[g] \subseteq x$ whenever $g \in x$, and is \emph{directed} if for all $g$, $h \in x$, $gP \cap hP \cap x \not= \emptyset$. Then $\Omega$ can be identified with the collection of all nonempty directed hereditary subsets of $P$.  Left reversibility implies that $P$ itself is a directed hereditary subset, corresponding to $\infty \in \Omega$. The finite directed hereditary subsets correspond to the elements of $P$ itself, via $g \in P \leftrightarrow [g] \subseteq P$.  Thus $\ti{\Omega}_{\infty}$ corresponds to $\{x \in \Omega : x \text{ is infinite and } x \not= P \}$.  For $k \ge 0$ let $X_k = \{ x \in \ti{\Omega}_{\infty} : \Delta^k \in x,\ \Delta^{k+1} \not\in x \}$.  Then $\ti{\Omega}_{\infty} = \sqcup_{k=0}^\infty X_k$. 

\blemma
\label{lem:infinite normal form}
$X_0$ is in one-to-one correspondence with infinite sequences $(g_1, g_2, \ldots)$ such that $g_i \in P_0$ and $\cR(g_i) \supseteq \cL(g_{i+1})$ for all $i$.  The correspondence pairs such a sequence with the directed hereditary set $\bigcup_{n=1}^\infty [g_1 \cdots g_n]$.
\elemma
\setlength{\parindent}{0cm} \setlength{\parskip}{0cm}
\bproof
First we claim that if $x \in X_0$ then $x$ contains a unique element of $P_0$ of maximal length.  To see this, first consider two elements $g$, $h \in P_0 \cap x$.  Since $x$ is directed there is $f \in x \cap gP \cap hP$.  But then $g \vee h \prec f$. Since $x$ is hereditary, we have $g \vee h \in x$.  Since $g$, $h \in P_0$ we have $g$, $h \prec \Delta$, and hence $g \vee h \in P_{\text{red}}$.  But since $g \vee h \in x$ and $x \in X_0$ we know that $g \vee h \not= \Delta$.  Therefore $g \vee h \in P_0 \cap x$.  Since $x \cap P_0$ is finite we may apply the preceding argument finitely many times to see that $\bigvee (x \cap P_0)$ is the desired element of maximal length in $x \cap P_0$.  Now let $x \in X_0$ and let $g_1$ be the maximal element of $x \cap P_0$.  We write $\sigma^{g_1}(x) = \{ h \in P : g_1 h \in x \} = g_1^{-1} (x \cap g_1 P)$ (with notation borrowed from \cite{Sp14}).  It is easy to check that $\sigma^{g_1}(x) \in X_0$. Let $g_2$ be maximal in $\sigma^{g_1}(x) \cap P_0$.  Then $\cR(g_1) \supseteq \cL(g_2)$: for, if $s \in \cL(g_2) \setminus \cR(g_1)$ we may write $g_2 = s g_2'$, and also we have $g_1 s \in P_{\text{red}}$. Since $g_1 s \prec g_1 g_2$ we know that $g_1s \in x$, and hence $g_1 s \in P_0$, contradicting the maximality of $g_1$.  Repeating this construction we obtain an infinite sequence $(g_1,g_2, \ldots)$ such that $\cR(g_i) \supseteq \cL(g_{i+1})$ for all $i$.  Since $g_1 \cdots g_n \in x$ for all $n$, we have $\bigcup_{n=1}^\infty [g_1 \cdots g_n] \subseteq x$. Conversely, let $h \in x$. We want to show that $h \in \bigcup_{n=1}^\infty [g_1 \cdots g_n]$. First note that if $p \in P_{\rm red}$ and $q \in P$ has normal form of length $n$, i.e., of the form $(q_1, \dotsc, q_n)$, then $p \vee q$ must have normal form of length $n$. This is straightforward using induction on $n$: For $n=1$, this follows from the fact that $P_{\rm red}$ is closed under $\vee$. Assume that our claim is true for $n-1$. Write $p \vee q_1 = q_1 r$ with $r \in P_{\rm red}$. Then $p \vee q = (p \vee q_1) \vee q = (q_1 r) \vee q = q_1 (r \vee (q_2 \dotsm q_n))$. By the induction hypothesis, $r \vee (q_2 \dotsm q_n)$ has normal form of length $n-1$, so that $q_1 (r \vee (q_2 \dotsm q_n))$ has normal form of length at most $n$. The length must then be equal to $n$ since $q$ has normal form of length $n$. Now consider $h \in x$ as above. Suppose that $h$ has normal form of length $n$. Then we claim that there exists $r \in P$ such that $hr = g_1 \dotsm g_n$. We prove this inductively on $n$. Our claim is obviously true for $n=1$. Assume it is true for $n-1$. By the above observation, we know that $g_1 \vee h$ has normal form of length $n$, say $(g_1, q_2, \dotsc, q_n)$. Hence there exists $s \in P$ such that $hs = g_1 q_2 \dotsm q_n$. Applying the induction hypothesis to $q_2 \dotsm q_n$ and $\sigma^{g_1}(x)$, we obtain $t \in P$ such that $q_2 \dotsm q_n t = g_2 \dotsm g_n$. Then set $r = st$. We have $hr = hst = g_1 q_2 \dotsm q_n t = g_1 g_2 \dotsm g_n$. This proves our claim. Hence it follows that $h \prec g_1 \dotsm g_n$, and thus $h \in \bigcup_{n=1}^\infty [g_1 \cdots g_n]$. Therefore, $x \subseteq \bigcup_{n=1}^\infty [g_1 \cdots g_n]$, as desired.
\quad
We have shown that for each $x \in X_0$ there is an infinite sequence $(g_1, g_2, \ldots)$ such that $\cR(g_i) \supseteq \cL(g_{i+1})$ for all $i$ and $x = \bigcup_{n=1}^\infty [g_1 \cdots g_n]$.  If $(h_1, h_2, \ldots)$ is a sequence in $P_0$ such that $\cR(h_i) \supseteq \cL(h_{i+1})$ for all $i$ and $x = \bigcup_{n=1}^\infty [h_1 \cdots h_n]$, then $g_1 = h_1$ since both equal the unique maximal element of $x \cap P_0$.  Inductively, $g_{i+1} = h_{i+1}$ is the unique maximal element of $\sigma^{g_1 \cdots g_i}(x) \cap P_0$.  Therefore the sequence $(g_1,g_2, \ldots)$ associated to $x$ is unique.  Finally, if $(g_1, g_2, \ldots)$ is any sequence in $P_0$ such that $\cR(g_i) \supseteq \cL(g_{i+1})$ for all $i$ then $\Delta \not\in \bigcup_{n=1}^\infty [g_1 \cdots g_n]$.  Therefore $(g_1, g_2, \ldots)$ is associated to the directed hereditary set $\bigcup_{n=1}^\infty [g_1 \cdots g_n]$ in $X_0$.
\eproof

\bremark
The above proof actually shows that $\ti{\Omega}_{\infty}$ is in one-to-one correspondence with sequences $(g_1,g_2,\ldots)$ such that $g_i \in P_{\text{red}} \setminus\{1\}$, $\cR(g_i) \supseteq \cL(g_{i+1})$ for all $i$, and $g_i \in P_0$ eventually.
\eremark

\bremark
\label{rem:cylinder sets}
We extend $\cL$ to $\ti{\Omega}_{\infty}$ by setting $\cL(x) := \cL(g_1)$ if $x$ has infinite normal form $(g_1, g_2, \ldots)$.  For $g \in P$ we let $Z(g) = \{x \in \ti{\Omega}_{\infty} : g \in x \}$.  We note that for $g$, $h \in P$ we have $Z(g) \cap Z(h) \cap X_0 \not= \emptyset$ if and only if $g \vee h \not\in \Delta P$.  In this case $Z(g) \cap Z(h) \cap X_0 = Z(g \vee h) \cap X_0$.  The collection $\{Z(g) \setminus \bigcup_{j=1}^n Z(h_j) : h_j \in gP \text{ for } 1 \le j \le n \}$ is a base of compact-open sets for the topology of $\ti{\Omega}_{\infty}$. Then $X_0$ is a compact-open subset of $\ti{\Omega}_{\infty}$ that is transversal, and hence the groupoid $G \ltimes \ti{\Omega}_{\infty}$ is equivalent to its restriction to $X_0$.
\eremark

\bprop
\label{prop:nonamenable}
Let $G$ be an Artin-Tits group of finite type with at least three generators, and suppose that there exist $s$, $t$ in $S$ with $2 < m_{st} < \infty$.  Then $G \ltimes \ti{\Omega}_{\infty}$ is not amenable.
\eprop
\setlength{\parindent}{0cm} \setlength{\parskip}{0cm}
\bproof
Let $G_{st}$ be the subgroup generated by $s$ and $t$.  Then $G_{st}$ is nonamenable.  Let $x = (\Delta_{s,t}, \Delta_{s,t}, \ldots) \in X_0$.  Then the isotropy of $G \ltimes \ti{\Omega}_{\infty}$  at $x$ contains $G_{st}$, and hence is nonamenable.  Since there is a point with nonamenable isotropy, the groupoid is nonamenable.
\eproof

We now give several lemmas on the existence of certain normal forms.  Note that since $\Delta$ is the unique element of $P_{\text{red}}$ of maximal length, $s \in \cL(\Delta)$ and $s \in \cR(\Delta)$ for all $s \in S$.  Moreover $\Delta$ is the only element of $P_{\text{red}}$ with this property: if $g \in P_{\text{red}}$ is such that $s \prec g$ for all $s \in S$, then $\ell(p(s)p(g)) = \ell(p(g)) - 1 < \ell(p(g))$ for all $s \in S$, so $p(g)$ must have maximal length in $W$. Therefore $g = \Delta$ (and a similar argument works on the right).  If $T \subseteq S$ we write $\Delta_T$ for the unique element in $P_{T,\text{red}}$ of maximal length ($P_{T,\text{red}} \subseteq P_{\text{red}}$ by Proposition \ref{prop:subset of generators}).

\blemma
\label{lem:delta for subgroup}
Let $T \subseteq S$ and $g \in P_{\text{red}}$ be such that $T \subseteq \cL(g)$.  Then $\Delta_T \prec g$.
\elemma
\bproof
Let $M = \{h \in P : h \in \langle T \rangle^+ \text{ and } h \prec g \}$.  Since $g \in P_{\text{red}}$ we know that $M \subseteq P_{\text{red}}$, and so $M$ is finite.  We verify the two conditions in \cite[Lemma 1.4]{Mi}.  For the first condition, let $h \in M$ and $k \in P$ be such that $k \prec h$.  By Proposition \ref{prop:subset of generators} we have $k \in P_T$, and hence $k \in M$.  For the second condition, let $k \in P$ and $s$, $t \in S$ such that $ks$, $kt \in M$.  By the first condition $k \in M \subseteq P_T$.  Therefore $s$, $t \in G_T$, hence $s$, $t \in T$.  Since $k \in M$ we have $k \prec g$, so $g = k g'$. Then $ks$, $kt \prec kg'$, so $s$, $t \prec g'$.  By \cite[Proposition 1.5]{Mi} we have $\Delta_{s,t} \prec g'$, and hence $k \Delta_{s,t} \prec g$.  Therefore $k \Delta_{s,t} \in M$.  Now by \cite[Lemma 1.4]{Mi} there is $f \in M$ such that $M = \{ h \in P : h \prec f \}$. Then $f \in P_{T,\text{red}}$.  Since $T \subseteq M$ we have $t \prec f$ for all $t \in T$.  Then $tf \not\in P_{T,\text{red}}$ for all $t \in T$, so $f$ is maximal in $P_{T,\text{red}}$.  Therefore $f = \Delta_T$.
\eproof
Note that we obtain an alternative proof for Lemma~\ref{lem:delta for subgroup} by observing that the canonical map $P_T \to P_S$ is $\vee$-preserving. This, in turn, follows for instance from the recipe for computing $\vee$ using right reversing (see \cite[Proposition~6.10]{Deh03}), which does not depend on the ambient monoid.

\blemma
\label{lem:commuting subsets}
Let $T$, $U \subseteq S$ be such that $m_{tu} = 2$ for all $t \in T$ and $u \in U$ (in particular, $T$ and $U$ are disjoint).  Let $g \in P_T$.  Then $\cL(g \Delta_U) = \cL(g) \cup U$ and $\cR(g \Delta_U) = \cR(g) \cup U$.
\elemma

\bproof
The containments $\supseteq$ are clear.  Let $s \in \cL(g \Delta_U)$ and suppose that $s \not\in U$.  Then since $U \subseteq \cL(g \Delta_U)$, Lemma \ref{lem:delta for subgroup} implies that $\Delta_{U \cup \{s\}} \prec g \Delta_U = \Delta_U g$, hence $\Delta_U s \prec \Delta_U g$, and hence $s \prec g$.  Therefore $s \in \cL(g)$ as required.  A similar argument on the right gives the second statement of the lemma.
\eproof

For the next several lemmas we will consider a portion of the Coxeter diagram  that is \emph{linear}, i.e. a tree with no vertex having valence greater than two. 

\blemma
\label{lem:linear L and R}
Let $a_1, \ldots, a_p \in S$ be such that $m_{a_\mu,a_{\mu+1}} > 2$ for $1 \le \mu < p$ and $m_{a_\mu,a_\nu} = 2$ if $|\mu - \nu| \ge 2$ (so that $a_1, \ldots, a_p$ forms a linear subgraph of the Coxeter diagram of $G$).  For $1 \le \mu \le p$ we have
\begin{align*}
(i) && \cL(a_\mu \Delta_{\{a_1,\ldots,a_{\mu-1}\}} \Delta_{\{a_{\mu+1},\ldots, a_p\}}) &= \{a_\nu : \nu \not=  \mu-1,\ \mu+1\} \\
(ii) && \cR(a_\mu \Delta_{\{a_1,\ldots,a_{\mu-1}\}} \Delta_{\{a_{\mu+1},\ldots,a_p\}}) &= \{a_\nu : \nu \not= \mu \} \\
(iii) && \cL(\Delta_{\{a_1,\ldots,a_{\mu-1}\}} \Delta_{\{a_{\mu+1},\ldots, a_p\}} a_\mu) &= \{a_\nu : \nu \not=  \mu \} \\
(iv) && \cR(\Delta_{\{a_1,\ldots,a_{\mu-1}\}} \Delta_{\{a_{\mu+1},\ldots,a_p\}} a_\mu) &= \{a_\nu : \nu \not= \mu-1,\ \mu+1 \} \\
(v) && \cL(a_{\mu-1} \Delta_{\{a_1,\ldots,a_{\mu-2}\}} \Delta_{\{a_{\mu+1},\ldots,a_p\}} a_\mu) &= \{a_\nu : \nu \not= \mu-2,\ \mu\} \\
(vi) && \cR(a_{\mu-1} \Delta_{\{a_1,\ldots,a_{\mu-2}\}} \Delta_{\{a_{\mu+1},\ldots,a_p\}} a_\mu) &= \{a_\nu : \nu \not= \mu-1,\ \mu+1\}.
\end{align*}
Moreover the above elements are in $P_0$.
\elemma

\bproof
It is clear that the elements are in $P_0$, and that in all parts the containments $\supseteq$ hold.  We prove the reverse containments.  

For (i) suppose that $a_{\mu-1} \in \cL(a_\mu \Delta_{\{a_1,\ldots,a_{\mu-1}\}} \Delta_{\{a_{\mu+1},\ldots, a_p\}})$. Then $a_\mu$, $a_{\mu-1} \in \cL(a_\mu \Delta_{\{a_1,\ldots,a_{\mu-1}\}} \Delta_{\{a_{\mu+1},\ldots, a_p\}})$, so by \cite[Proposition 2.5]{Mi} we have $a_\mu a_{\mu-1} a_\mu \prec \Delta_{\{a_{\mu-1},a_\mu\}} \prec a_\mu \Delta_{\{a_1,\ldots,a_{\mu-1}\}} \Delta_{\{a_{\mu+1},\ldots, a_p\}}$.  But then $a_{\mu-1} a_\mu \prec \Delta_{\{a_1,\ldots,a_{\mu-1}\}} \Delta_{\{a_{\mu+1},\ldots, a_p\}}$, which contradicts the fact that $a_\mu \not\in G_{\{a_\nu : \nu \not= \mu \}}$ (by Proposition \ref{prop:subset of generators}).  A similar argument shows that $a_{\mu+1} \not\in \cL(a_\mu \Delta_{\{a_1,\ldots,a_{\mu-1}\}} \Delta_{\{a_{\mu+1},\ldots, a_p\}})$. For (ii) suppose $a_\mu \in \cR(a_\mu \Delta_{\{a_1,\ldots,a_{\mu-1}\}} \Delta_{\{a_{\mu+1},\ldots, a_p\}})$.  Then $a_1$, $\ldots$, $a_p \in \cR(a_\mu \Delta_{\{a_1,\ldots,a_{\mu-1}\}} \Delta_{\{a_{\mu+1},\ldots, a_p\}})$. By Lemma \ref{lem:delta for subgroup} we have $a_\mu \Delta_{\{a_1,\ldots,a_{\mu-1}\}} \Delta_{\{a_{\mu+1},\ldots, a_p\}} \succ \Delta_{\{a_1,\ldots,a_p\}}$.  By \cite[Proposition 3.1]{Mi}, applied to $G_{\{a_1,\ldots,a_p\}}$, we have $\Delta_{\{a_1,\ldots,a_p\}} \prec a_\mu \Delta_{\{a_1,\ldots,a_{\mu-1}\}} \Delta_{\{a_{\mu+1},\ldots, a_p\}}$, contradicting (i). The proof of (iii) is similar to that of (ii), and the proof of (iv) is similar to that of (i).  For (v), first suppose that $a_{\mu-2} \in \cL(a_{\mu-1} \Delta_{\{a_1,\ldots,a_{\mu-2}\}} \Delta_{\{a_{\mu+1},\ldots,a_p\}} a_\mu)$. Then $a_{\mu-2}$, $a_{\mu-1} \in \cL(a_{\mu-1} \Delta_{\{a_1,\ldots,a_{\mu-2}\}} \Delta_{\{a_{\mu+1},\ldots,a_p\}} a_\mu)$.  Now the proof is similar to that of (i).  Next suppose that $a_\mu \in \cL(a_{\mu-1} \Delta_{\{a_1,\ldots,a_{\mu-2}\}} \Delta_{\{a_{\mu+1},\ldots,a_p\}} a_\mu)$. As before we then have $a_{\mu-1} a_\mu a_{\mu-1} \prec a_{\mu-1} \Delta_{\{a_1,\ldots,a_{\mu-2}\}} \Delta_{\{a_{\mu+1},\ldots,a_p\}} a_\mu$, and hence $a_\mu a_{\mu-1} \prec \Delta_{\{a_1,\ldots,a_{\mu-2}\}} \Delta_{\{a_{\mu+1},\ldots,a_p\}} a_\mu$, and hence $a_\mu \in \cL(\Delta_{\{a_1,\ldots,a_{\mu-2}\}} \Delta_{\{a_{\mu+1},\ldots,a_p\}} a_\mu)$.  By Lemma \ref{lem:commuting subsets} we then have $a_\mu \in \cL(\Delta_{\{a_{\mu+1},\ldots,a_p\}} a_\mu)$, contradicting (iii) (applied to the subset $\{a_\mu,\ldots,a_p\}$).  The proof of (vi) is similar to that of (v).  Finally, it is clear that, for example, $\Delta_{\{a_1,\ldots,a_{\mu-1}\}} \Delta_{\{a_{\mu+1},\ldots, a_p\}} \in P_{\{a_1,\ldots,a_p\} \setminus \{a_\mu\},\text{red}}$ and that $a_\mu \not\in \cL(\Delta_{\{a_1,\ldots,a_{\mu-1}\}} \Delta_{\{a_{\mu+1},\ldots, a_p\}})$. Therefore $a_\mu \Delta_{\{a_1,\ldots,a_{\mu-1}\}} \Delta_{\{a_{\mu+1},\ldots, a_p\}} \in P_{\{a_1,\ldots,a_p\},\text{red}}$.  Similar arguments demonstrate the same for parts (ii) - (vi).
\eproof

\bdefin
Let $T \subseteq S$.  For $\emptyset \not= T_1$, $T_2 \subsetneq T$ we write $T_1 \sim_T T_2$ if there is $g \in P_T$ with normal form $(g_1, \ldots, g_j)$ such that $\cL(g_1) = T_1$ and $\cR(g_j) = T_2$.
\edefin

\bremark
It is clear that $\sim_T$ is an equivalence relation on the proper nonempty subsets of $T$. ($T$ cannot be $\sim_T$-equivalent to a proper subset by \cite[Proposition 3.1]{Mi}.)
\eremark

\blemma
\label{lem:singletons}
If $T$ determines a connected portion of the Coxeter diagram of $G$ then $\{s\} \sim_T \{t\}$ for all $s$, $t \in T$.
\elemma

\bproof
Let $s$, $t \in T$. Since $T$ is connected there are $s_1$, $\ldots$, $s_k \in T$ defining a linear subgraph of the Coxeter diagram (as in Lemma \ref{lem:linear L and R}) such that $s = s_1$ and $t = s_k$.  Let $g = s_1 s_2^2 s_3^2 \cdots s_{k-1}^2 s_k$.  Then $g$ has normal form $(s_1 s_2, s_2 s_3, \ldots, s_{k-1} s_k)$, and hence $\cL(g_1) = \cL(s_1 s_2) = s_1$ and $\cR(g_{k-1}) = \cR(s_{k-1}s_k) = s_k$.
\eproof

\blemma
\label{lem:subequivalence}
Let $T \subseteq U \subseteq S$.  Then $\sim_T \subseteq \sim_U$.
\elemma

\bproof
This follows from Lemma \ref{lem:subset of generators}.
\eproof

\bremark
By Lemma \ref{lem:subequivalence}, equivalences relative to a subgraph still hold relative to a larger subgraph.  We will use this as needed without further mention.
\eremark

In the following Lemmas~\ref{lem:first equivalence} -- \ref{lem:sixth equivalence}, let $T = \{s_1,\ldots,s_k\} \subseteq S$ be a linear subgraph of the Coxeter diagram of $G$ (as in Lemma \ref{lem:linear L and R}).
\blemma
\label{lem:first equivalence}
Let $1 < i < k$.  Then $\{s_1,\ldots,s_{i-1}\} \sim_{\{s_1,\ldots,s_{i+1}\}} \{s_1,\ldots,s_i\}$.
\elemma

\bproof
Letting $\{s_1,\ldots,s_i\}$ play the role of $\{a_1,\ldots,a_p\}$ in Lemma \ref{lem:linear L and R}~(iii) and (iv), and letting $\mu = i$, we have $\cL(\Delta_{\{s_1,\ldots,s_{i-1}\}}s_i) = \{s_1,\ldots,s_{i-1}\}$ and $\cR(\Delta_{\{s_1,\ldots,s_{i-1}\}}s_i) = \{s_1,\ldots,s_{i-2},s_i\}$ (or $\{s_i\}$ if $i=2$). Next, letting $\{s_1,\ldots,s_{i+1}\}$ play the role of $\{a_1,\ldots,a_p\}$ in Lemma \ref{lem:linear L and R}~(v) and (vi), and letting $\mu = i+1$, we have $\cL(s_i \Delta_{\{s_1,\ldots,s_{i-1}\}}s_{i+1}) = \{s_1,\ldots,s_{i-2},s_i\}$ (or $\{s_i\}$ if $i=2$) and $\cR(s_i \Delta_{\{s_1,\ldots,s_{i-1}\}}s_{i+1}) = \{s_1,\ldots,s_{i-1},s_{i+1}\}$. Finally, letting $\{s_1,\ldots,s_{i+1}\}$ play the role of $\{a_1,\ldots,a_p\}$ in Lemma \ref{lem:linear L and R}~(i) and (ii), and letting $\mu = i+1$, we have $\cL(s_{i+1} \Delta_{\{s_1,\ldots,s_i\}}) = \{s_1,\ldots,s_{i-1},s_{i+1}\}$ and $\cR(s_{i+1} \Delta_{\{s_1,\ldots,s_i\}}) = \{s_1,\ldots,s_i\}$.  Let $g_1 = \Delta_{\{s_1,\ldots,s_{i-1}\}}s_i$, $g_2 = s_i \Delta_{\{s_1,\ldots,s_{i-1}\}}s_{i+1}$, and $g_3 = s_{i+1} \Delta_{\{s_1,\ldots,s_i\}}$.  Then $g = g_1g_2g_3$ has normal form $(g_1,g_2,g_3)$, and $\cL(g_1) = \{s_1,\ldots,s_{i-1}\}$, $\cR(g_3) = \{s_1,\ldots,s_i\}$.
\eproof

\blemma
\label{lem:second equivalence}
Let $1 < i < k$. Then $\{s_i,\ldots,s_k\} \sim_{\{s_{i-1},\ldots,s_k\}} \{s_{i+1},\ldots,s_k\}$.
\elemma

\bproof
The proof is analogous to the proof of Lemma \ref{lem:first equivalence}.
\eproof

\blemma
\label{lem:third equivalence}
Let $1 < i < k$.  Then $\{s_1,s_{i+1},\ldots,s_k\} \sim_{\{s_1,s_i,\ldots,s_k\}} \{s_1,s_k\}$.
\elemma

\bproof
Applying Lemma \ref{lem:second equivalence} repeatedly we find that $\{s_{i+1},\ldots,s_k\} \sim_{\{s_i,\ldots,s_k\}} \{s_k\}$. Suppose that $i > 2$.  Then in fact we have $\{s_{i+1},\ldots,s_k\} \sim_{\{s_3,\ldots,s_k\}} \{s_k\}$.  By Lemma \ref{lem:commuting subsets} we then have $\{s_1,s_{i+1},\ldots,s_k\} \sim_{\{s_1,s_i,\ldots,s_k\}} \{s_1,s_k\}$.  For the case where $i = 2$, we consider the proof of Lemma \ref{lem:second equivalence} for this case.  For that we start with the normal form $(\Delta_{\{s_3,\ldots,s_k\}}s_2, s_2 \Delta_{\{s_4,\ldots,s_k\}}s_3, s_3 \Delta_{\{s_4,\ldots,s_k\}})$, for which $\cL(\cdot_1) = \{s_3,\ldots,s_k\}$ and $\cR(\cdot_3) = \{s_4,\ldots,s_k\}$.  We modify this to include $s_1$ as follows.  It is only in the first two terms of the normal form that $s_2$ occurs, obstructing the use of Lemma \ref{lem:commuting subsets}.  By Lemma \ref{lem:linear L and R}~(v) and (vi), with $\{s_1,\ldots,s_k\}$ playing the role of $\{a_1,\ldots,a_p\}$ and $\mu = 2$, we have $\cL(s_1 \Delta_{\{s_3,\ldots,s_k\}} s_2) = \{s_1,s_3,\ldots,s_k\}$ and $\cR(s_1 \Delta_{\{s_3,\ldots,s_k\}} s_2) = \{s_2,s_4,\ldots,s_k\}$.  By Lemma \ref{lem:linear L and R}~(v) and (vi), with $\{s_1,\ldots,s_k\}$ playing the role of $\{a_1,\ldots,a_p\}$ and $\mu = 3$, we have $\cL(s_2 s_1 \Delta_{\{s_4,\ldots,s_k\}} s_3) = \{s_2,s_4,\ldots,s_k\}$ and $\cR(s_2 s_1 \Delta_{\{s_4,\ldots,s_k\}} s_3) = \{s_1,s_3,s_5,\ldots,s_k\}$.  By Lemma \ref{lem:commuting subsets}, and Lemma \ref{lem:linear L and R}~(i) and (ii) with $\{s_3,\ldots,s_k\}$ playing the role of $\{a_1,\ldots,a_p\}$ and $\mu = 1$, we have $\cL(s_1 s_3 \Delta_{\{s_4,\ldots,s_k\}}) = \{s_1\} \cup \cL(s_3 \Delta_{\{s_4,\ldots,s_k\}}) = \{s_1,s_3,s_5,\ldots,s_k\}$ and $\cR(s_1 s_3 \Delta_{\{s_4,\ldots,s_k\}}) = \{s_1\} \cup \cR(s_3 \Delta_{\{s_4,\ldots,s_k\}}) = \{s_1,s_4,\ldots,s_k\}$. So $(s_1 \Delta_{\{s_3,\ldots,s_k\}} s_2, s_2 s_1 \Delta_{\{s_4,\ldots,s_k\}} s_3, s_1 s_3 \Delta_{\{s_4,\ldots,s_k\}})$ is a normal form with $\cL(\cdot_1) = \{s_1,s_3,\ldots,s_k\}$ and $\cR(\cdot_3) = \{s_1,s_4,\ldots,s_k\}$.  If we combine this with the equivalence at the beginning of the proof we are finished.
\eproof

\blemma
\label{lem:fourth equivalence}
Let $1 < i < k$. Then $\{s_1,s_{i+1},\ldots,s_k\} \sim_T T \setminus \{s_i\}$.
\elemma

\bproof
By repeated application of Lemma \ref{lem:first equivalence} we have that $\{s_1\} \sim_{\{s_1,\ldots,s_{i-1}\}} \{s_1,\ldots,s_{i-2}\}$. Then by Lemma \ref{lem:commuting subsets}, multiplication by $\Delta_{\{s_{i+1},\ldots,s_k\}}$ gives
\begin{equation}
\label{lem:fourth equivalence eqn 1}
\{s_1,s_{i+1},\ldots,s_k\} \sim_{T \setminus \{s_i\}} \{s_1,\ldots,s_{i-2},s_{i+1},\ldots,s_k\}.
\end{equation}
From Lemma \ref{lem:linear L and R}~(iii) and (iv), with $\{s_1,\ldots,s_{i-1}\}$ playing the role of $\{a_1,\ldots,a_p\}$ and $\mu = i-1$, we find that $\{s_1,\ldots,s_{i-2}\} \sim_{\{s_1,\ldots,s_{i-1}\}} \{s_1,\ldots,s_{i-3},s_{i-1}\}$. Then by Lemma \ref{lem:commuting subsets}, multiplication by $\Delta_{\{s_{i+1},\ldots,s_k\}}$ gives
\begin{equation}
\label{lem:fourth equivalence eqn 2}
\{s_1,\ldots,s_{i-2},s_{i+1},\ldots,s_k\} \sim_{\{s_1,\ldots,s_{i-1},s_{i+1},\ldots,s_k\}} \{s_1,\ldots,s_{i-3},s_{i-1},s_{i+1},\ldots,s_i\}.
\end{equation}
By Lemma \ref{lem:linear L and R}~(v) and (vi), with $\{s_1,\ldots,s_k\}$ playing the role of $\{a_1,\ldots,a_p\}$ and $\mu = i$, we have
\begin{equation}
\label{lem:fourth equivalence eqn 3}
\{s_1,\ldots,s_{i-3},s_{i-1},s_{i+1},\ldots,s_k\} \sim_T \{s_1,\ldots,s_{i-2},s_i,s_{i+2},\ldots,s_k\}.
\end{equation}
By Lemma \ref{lem:linear L and R}~(i) and (ii), with $\{s_1,\ldots,s_k\}$ playing the role of $\{a_1,\ldots,a_p\}$ and $\mu = i$, we have
\begin{equation}
\label{lem:fourth equivalence eqn 4}
\{s_1,\ldots,s_{i-2},s_i,s_{i+2},\ldots,s_k\} \sim_T \{s_1,\ldots,s_{i-1},s_{i+1},\ldots,s_k\}.
\end{equation}
The combination of equivalences \eqref{lem:fourth equivalence eqn 1} - \eqref{lem:fourth equivalence eqn 4} proves the lemma.
\eproof

\blemma
\label{lem:fifth equivalence}
We have $\{s_1\} \sim_T \{s_1,s_k\}$.
\elemma

\bproof
We know that $\cL(s_1s_2) = \{s_1\}$ and $\cR(s_1s_2) = \{s_2\}$.  By Lemma \ref{lem:linear L and R}~(i) and (ii), with $\{s_1,s_2,s_3\}$ playing the role of $\{a_1,\ldots,a_p\}$ and $\mu = 2$, we have $\cL(s_2s_1s_3) = \{s_2\}$ and $\cR(s_2s_1s_3) = \{s_1,s_3\}$.  For $i > 2$ we know that $\cL(s_is_{i+1}) = \{s_i\}$ and $\cR(s_is_{i+1}) = \{s_{i+1}\}$.  Then by Lemma \ref{lem:commuting subsets} we have $\cL(s_1s_is_{i+1}) = \{s_1,s_i\}$ and $\cR(s_1s_is_{i+1}) = \{s_1,s_{i+1}\}$.  Combining these gives the lemma.
\eproof

\blemma
\label{lem:sixth equivalence}
Let $1 \le i, j \le k$. Then $\{s_j\} \sim_T T \setminus \{s_i\}$.
\elemma

\bproof
First suppose $1 < i$. By Lemma \ref{lem:singletons} we have $\{s_j\} \sim_T \{s_1\}$. By Lemma \ref{lem:fifth equivalence} we have $\{s_1\} \sim_T \{s_1,s_k\}$. By Lemma \ref{lem:third equivalence} we have $\{s_1,s_k\} \sim_T \{s_1,s_{i+1},\ldots,s_k\}$. By Lemma \ref{lem:fourth equivalence} we have $\{s_1,s_{i+1},\ldots,s_k\} \sim_T T \setminus \{s_i\}$. Combining these equivalence gives the desired result.  For the case $i = 1$, note that repeated use of Lemma \ref{lem:second equivalence} gives $\{s_2,\ldots,s_k\} \sim_T \{s_k\}$.  Combining this with Lemma \ref{lem:singletons} gives $\{s_j\} \sim_T \{s_2,\ldots,s_k\}$.
\eproof

If $G$ is irreducible and of finite type then the Coxeter diagram must be a tree with at most one vertex having valence greater than two; if there is such a vertex it has valence three.  Moreover, in the case that there is a vertex of valence three, let $v \in S$ be the vertex of valence three, and let the three linear pieces of $S \setminus \{v\}$ have $n_1$, $n_2$, $n_3$ vertices, respectively.  Then we may assume that $n_3 = 1$ and that $n_2 = 1$ or 2  (\cite[Chapter VI, Section 4.1, Th\'eor\`eme 1]{B}).

\bremark
\label{rem:tree notation}
We next will generalize the previous lemmas to the case where the Coxeter diagram is a tree of this type.  Thus in the following Lemmas~\ref{lem:tree equivalence one} -- \ref{lem:tree equivalence seven}, we will assume that $S = \{s_1,\ldots,s_k,u\}$ with $k \ge 3$, that $m_{s_i,s_{i+1}} > 2$ for $1 \le i < k$, that $m_{s_\ell,u} > 2$ for $\ell = 2$ or 3, that $\ell < k$, and that $m_{t,t'} = 2$ for all other pairs of elements of $S$.
\eremark

\blemma
\label{lem:tree equivalence one}
In Lemma \ref{lem:linear L and R} let $\{s_1,\ldots,s_k\}$ play the role of $\{a_1,\ldots,a_p\}$.  The elements of $P_{\text{red}}$ in Lemma \ref{lem:linear L and R} may be modified so that $u$ is included in all of their $\cL$ and $\cR$ sets.
\elemma

\bproof
We indicate the procedure for the first element; the other two are managed similarly.  Let $1 \le i \le k$ and consider the element $s_i \Delta_{\{s_1,\ldots,s_{i-1}\}} \Delta_{\{s_{i+1},\ldots,s_k\}}$.  There are three cases depending on the relation between $i$ and $\ell$.  Recall that $m_{s_j,u} > 2$ if and only if $j = \ell$.

\noindent
\emph{Case 1:} $\ell < i$.  We consider the element $g = s_i \Delta_{\{u,s_1,\ldots,s_{i-1}\}} \Delta_{\{s_{i+1},\ldots,s_k\}}$.

\noindent
\emph{Case 2:} $\ell > i$.  We consider the element $g = s_i \Delta_{\{s_1,\ldots,s_{i-1}\}} \Delta_{\{u,s_{i+1},\ldots,s_k\}}$.

\noindent
\emph{Case 3:} $\ell = i$.  We consider the element $g = \Delta_{\{u,s_i\}} \Delta_{\{s_1,\ldots,s_{i-1}\}} \Delta_{\{s_{i+1},\ldots,s_k\}}$.

In all three cases we claim that $\cL(g) = \{u\} \cup \{s_j : j \not= i-1,\ i+1 \}$ and $\cR(g) = \{u\} \cup \{s_j : j \not= i \}$. The containments $\supseteq$ are clear. The proofs of the containments $\subseteq$ are identical to those in the proof of Lemma \ref{lem:linear L and R}~(i) and (ii).  An analogous argument works for parts (iii)-(vi) of Lemma \ref{lem:linear L and R}.
\eproof

\blemma
\label{lem:tree equivalence two}
For $1 \le i < j \le k$, $\{u,s_i\} \sim_S \{u,s_j\}$.
\elemma

\bproof
If $\ell < i$ or $\ell > j$ we have the normal form $(us_is_{i+1}, us_{i+1}s_{i+2}, \ldots, us_{j-1}s_j)$; the $\cL$ and $\cR$ sets of the terms follow from Lemma \ref{lem:commuting subsets}.  If $i < \ell < j$ we have the normal form $(us_is_{i+1}, \ldots, u s_{\ell-2}s_{\ell-1}, us_{\ell-1}s_\ell, s_\ell s_{\ell+1} u, \ldots, u s_{j-1}s_j)$; the $\cL$ and $\cR$ sets of the terms follow from Lemma \ref{lem:commuting subsets} except for the two terms involving $s_\ell$. For $us_{\ell - 1}s_\ell$ we may use Lemma \ref{lem:linear L and R}~(iii) and (iv) (with $a_1 = u$, $a_2 = s_\ell$, $a_3 = s_{\ell - 1}$, and $\mu = 2$), while for $s_\ell s_{\ell + 1} u$ we may use Lemma \ref{lem:linear L and R}~(iii) and (iv) (with $a_1 = u$, $a_2 = s_\ell$, $a_3 = s_{\ell + 1}$, and $\mu = 1$).  If $\ell = i$ we use the normal form $(\Delta_{ \{ u,s_\ell \} } s_{\ell + 1}, u s_{\ell + 1} s_{\ell + 2}, \cdots, u s_{j-1} s_j)$, where for the first term we use Lemma \ref{lem:linear L and R}~(iii) and (iv) (with $a_1 = u$, $a_2 = s_\ell$, $a_3 = s_{\ell + 1}$, and $\mu = 3$).  If $\ell = j$ we use the normal form $(u s_i s_{i+1}, \ldots, u s_{\ell-2} s_{\ell -1}, s_{\ell - 1} \Delta_{ \{u, s_\ell\} })$, where for the last term we use Lemma \ref{lem:linear L and R}~(i) and (ii) (with $a_1 = s_{\ell - 1}$, $a_2 = s_\ell$, $a_3 = u$, and $\mu = 1$).
\eproof

\blemma
\label{lem:tree equivalence three}
We have $\{u,s_1\} \sim_S \{u,s_1,s_k\}$.
\elemma

\bproof
If $\ell = 2$ we have the normal form $(us_1s_2, s_2s_1s_3u, us_1s_3s_4, \ldots us_1s_{k-1}s_k)$: the $\cL$ and $\cR$ sets of the first term are calculated using Lemma \ref{lem:linear L and R}, for the second term the calculation is similar to the arguments in the proof of that lemma, and for the remaining terms Lemma \ref{lem:commuting subsets} applies.  If $\ell = 3$ we have the normal form $(us_1s_2, us_2s_1s_3, s_1s_3s_4u, us_1s_4s_5, \ldots, us_1s_{k-1}s_k)$.  The arguments here are analogous.  The cases $\ell = 2$, $k = 3$ and $\ell = 3$, $k = 4$ are slightly different.  When $\ell = 2$ and $k = 3$ we use the normal form $(u s_1 s_2, s_2 s_1 s_3 u)$.  The calculation of $\cL(s_2 s_1 s_3 u)$ and $\cR(s_2 s_1 s_3 u)$ are as in the proof of Lemma \ref{lem:linear L and R}.   When $\ell = 3$ and $k = 4$ we use the normal form $(u s_1 s_2, u s_2 s_1 s_3, s_1 s_3 s_4 u)$; again, calculation of the relevant $\cL$ and $\cR$ sets is as in the proof of Lemma \ref{lem:linear L and R}.
\eproof

\blemma
\label{lem:tree equivalence four}
Let $1 < i < k$.  Then $\{u,s_1,\ldots,s_{i-1}\} \sim_{\{u,s_1,\ldots,s_{i+1}\}} \{u,s_1,\ldots,s_i\}$.
\elemma

\bproof
This is the same as Lemma \ref{lem:first equivalence} with $u$ included in all $\cL$ and $\cR$ sets.  Since the proof of Lemma \ref{lem:first equivalence} used only Lemma \ref{lem:linear L and R}, an analogous proof can be given here using Lemma \ref{lem:tree equivalence one}.
\eproof

\blemma
\label{lem:tree equivalence five}
Let $1 < i < k$. Then $\{u,s_i,\ldots,s_k\} \sim_{\{u,s_{i-1},\ldots,s_k\}} \{u,s_{i+1},\ldots,s_k\}$.
\elemma

\bproof
This is the same as Lemma \ref{lem:second equivalence} with $u$ included in all $\cL$ and $\cR$ sets.  The proof is analogous to that of Lemma \ref{lem:tree equivalence four}.
\eproof

\blemma
\label{lem:tree equivalence six}
Let $1 < i < k$. Then $\{u,s_1,s_{i+1},\ldots,s_k\} \sim_S \{u,s_1,s_k\}$.
\elemma

\bproof
Applying Lemma \ref{lem:tree equivalence five} repeatedly we find that $\{u,s_{i+1},\ldots,s_k\} \sim_{\{u,s_i,\ldots,s_k\}} \{u,s_k\}$.  As in the proof of Lemma \ref{lem:third equivalence}, if $i > 2$ then Lemma \ref{lem:commuting subsets} finishes the proof.  Suppose $i = 2$. First consider the case $\ell = 2$.  We start with the normal form from Lemma \ref{lem:third equivalence}: $(\Delta_{\{s_3,\ldots,s_k\}}s_2, s_2 \Delta_{\{s_4,\ldots,s_k\}}s_3, s_3 \Delta_{\{s_4,\ldots,s_k\}})$.  As in Lemma \ref{lem:tree equivalence one} we have the modified version $(\Delta_{\{s_3,\ldots,s_k\}}\Delta_{\{u,s_2\}}, \Delta_{\{u,s_2\}} \Delta_{\{s_4,\ldots,s_k\}}s_3, us_3 \Delta_{\{s_4,\ldots,s_k\}})$.  Now, imitating the proof of Lemma \ref{lem:third equivalence} we claim that there is a normal form $(s_1\Delta_{\{s_3,\ldots,s_k\}}\Delta_{\{u,s_2\}}, \Delta_{\{u,s_2\}} s_1 \Delta_{\{s_4,\ldots,s_k\}}s_3, us_1s_3 \Delta_{\{s_4,\ldots,s_k\}})$. We have
\begin{align*}
\cL(s_1\Delta_{\{s_3,\ldots,s_k\}}\Delta_{\{u,s_2\}}) &= \{u,s_1,s_3,\ldots,s_k\}, \text{ by Lemma \ref{lem:tree equivalence one}} \\
\cR(s_1\Delta_{\{s_3,\ldots,s_k\}}\Delta_{\{u,s_2\}}) &= \{u,s_2,s_4,\ldots,s_k\}, \text{ by Lemma \ref{lem:tree equivalence one}} \\
\cL(\Delta_{\{u,s_2\}} s_1 \Delta_{\{s_4,\ldots,s_k\}}s_3) &= \{u,s_2,s_4,\ldots,s_k\}, \text{ by Lemma \ref{lem:tree equivalence one}} \\
\cR(\Delta_{\{u,s_2\}} s_1 \Delta_{\{s_4,\ldots,s_k\}}s_3) &= \{u,s_1,s_3,s_5,\ldots,s_k\}, \text{ by Lemma \ref{lem:tree equivalence one}} \\
\cL(us_1s_3 \Delta_{\{s_4,\ldots,s_k\}}) &= \{u,s_1,s_3,s_5,\ldots,s_k\} \\
\cR(us_1s_3 \Delta_{\{s_4,\ldots,s_k\}}) &= \{u,s_1,s_4,\ldots,s_k\},\end{align*}
where to justify the last two equalities we apply Lemma \ref{lem:tree equivalence one} to $us_3 \Delta_{\{s_4,\ldots,s_k\}}$, and then Lemma \ref{lem:commuting subsets} to incude $s_1$.  In the case $\ell = 3$ we use the normal form $(s_1 \Delta_{\{u, s_3, \ldots, s_k\}} s_2, s_2 s_1 \Delta_{\{ s_4, \ldots, s_k \}} \Delta_{\{s_3,u\}}, s_1 \Delta_{\{u,s_3\}} \Delta_{\{s_4, \ldots, s_k\}})$.  The $\cL$ and $\cR$ sets of the first term are verified using Lemma \ref{lem:tree equivalence one}, of the third term using Lemmas \ref{lem:tree equivalence one} and \ref{lem:commuting subsets}, and of the second term by arguments analogous to those in the proof of Lemma \ref{lem:linear L and R}.
\eproof

\blemma
\label{lem:tree equivalence seven}
Let $1 < i < k$. Then $\{u,s_1,s_{i+1},\ldots,s_k\} \sim_S S \setminus \{s_i\}$.
\elemma

\bproof
This is the same as Lemma \ref{lem:fourth equivalence} with $u$ included in all $\cL$ and $\cR$ sets. The proof of Lemma \ref{lem:fourth equivalence} used Lemmas \ref{lem:first equivalence} and \ref{lem:linear L and R}, so inclusion of $u$ to the $\cL$ and $\cR$ sets is compatible with those lemmas (by Lemmas \ref{lem:tree equivalence four} and \ref{lem:tree equivalence one}).  The proof of Lemma \ref{lem:fourth equivalence} also used Lemma \ref{lem:commuting subsets}, and this does not immediately allow inclusion of $u$.  We consider the two possibilities for $\ell$.  If $\ell = 2$, then $\ell \le i$.  In this case $s_{i+1}$, $\ldots$, $s_k$ do not have edges connecting to $u$, so the use of Lemma \ref{lem:commuting subsets} to multiply by $\Delta_{\{s_{i+1},\ldots,s_k\}}$ is valid even with $u$ included  in the $\cL$ and $\cR$ sets.  If $\ell = 3$ then it is only if $i=2$ that there is an edge between $\{u\}$ and $\{s_{i+1},\ldots,s_k\}$, namely between $u$ and $s_3$.  But in this case the required statement is $\{s_1,s_3,\ldots,s_k\} \sim_{\{s_1,\ldots,s_k\}} \{s_1,s_3,\ldots,s_k\}$, which is true.
\eproof

\bprop
\label{prop:singleton-subset equivalence}
Let $G$ be an irreducible Artin-Tits group of finite type with generating set $S$ having more than two elements.  Let $s$, $t \in S$. Then $\{t\} \sim_S S \setminus \{s\}$.
\eprop

\bproof
If the Coxeter diagram of $G$ is linear then this follows from Lemma \ref{lem:sixth equivalence}.  Now suppose that the Coxeter diagram of $G$ is a tree having one vertex of valence three.  We adopt the notation of Remark \ref{rem:tree notation}. We first consider the case where $s = s_i$ for some $1 \le i \le k$.  Then
\begin{align*}
\{t\} &\sim_S \{s_\ell\}, \text{ by Lemma \ref{lem:singletons},} \\
&\sim_S \{u,s_\ell\}, \text{ by Lemma \ref{lem:first equivalence},} \\
&\sim_S \{u,s_1\}, \text{ by Lemma \ref{lem:tree equivalence two},} \\
&\sim_S \{u,s_1,s_k\}, \text{ by Lemma \ref{lem:tree equivalence three},} \\
&\sim_S \{u,s_1,s_{i+1},\ldots,s_k\}, \text{ by Lemma \ref{lem:tree equivalence six},} \\
&\sim_S \{u,s_1,\ldots,s_{i-1},s_{i+1},\ldots,s_k\}, \text{ by Lemma \ref{lem:tree equivalence seven},} \\
&= S \setminus \{s_i\}.
\end{align*}
Finally we treat the case where $s=u$. Consider the element $\Delta_{\{s_1,\ldots,s_k\}}us_\ell$.  Since $u \not\in \cR(\Delta_{\{s_1,\ldots,s_k\}})$, we know that $\Delta_{\{s_1,\ldots,s_k\}}u \in P_{\text{red}}$.  Next we note that $s_\ell \not\in \cR(\Delta_{\{s_1,\ldots,s_k\}}u)$, by the same argument used in the proof of Lemma \ref{lem:linear L and R}.  Therefore $\Delta_{\{s_1,\ldots,s_k\}}us_\ell \in P_{\text{red}}$. The arguments used in the proof of Lemma \ref{lem:linear L and R} can be used to show that $\cL(\Delta_{\{s_1,\ldots,s_k\}}us_\ell) = \{s_1,\ldots,s_k\}$ and $\cR(\Delta_{\{s_1,\ldots,s_k\}}us_\ell) = \{s_1,\ldots,s_k\} \setminus \{s_{\ell-1},s_{\ell+1}\}$.  Now we have $\{s_1,\ldots,s_k\} \setminus \{s_{\ell-1},s_{\ell+1}\} \sim_{\{s_1,\ldots,s_k\}} \{s_1,\ldots,s_k\} \setminus \{s_\ell\}$ by Lemma \ref{lem:linear L and R}, and $\{s_1,\ldots,s_k\} \setminus \{s_\ell\} \sim_{\{s_1,\ldots,s_k\}} \{s_1\}$ by Lemma \ref{lem:sixth equivalence}.  Combining these equivalences we have $S \setminus \{u\} = \{s_1,\ldots,s_k\} \sim_S \{s_1,\ldots,s_k\} \setminus \{s_{\ell-1},s_{\ell+1}\} \sim_S \{s_1,\ldots,s_k\} \setminus \{s_\ell\} \sim_S \{s_1\} \sim_S \{t\}$, the last equivalence following from Lemma \ref{lem:singletons}.
\eproof

\btheo
\label{Thm:simplepi}
Let $G$ be an irreducible Artin-Tits group of finite type with at least three generators. Then $C^*_r(G \ltimes \ti{\Omega}_{\infty})$ is simple and purely infinite.
\etheo
\bproof
Since $X_0$ is a compact-open transversal in $\ti{\Omega}_{\infty}$ it is enough to consider the restriction to $X_0$.  Thus we understand that the sets in the proof actually represent their intersections with $X_0$.  We will consider a (nonempty) basic open set $U = Z(g) \setminus \bigcup_{i=1}^n Z(h_i)$ in $X_0$, as described in Remark \ref{rem:cylinder sets} (thus $g$, $h_1$, $\ldots$, $h_n \in P \setminus \Delta P$).  Let $y \in U$. Write the infinite normal form of $y$ as $(y_1,y_2, \ldots)$. Choose $p > \max \{ \nu(g), \nu(h_1), \ldots, \nu(h_n) \}$.  We claim that if $x \in X_0$ has infinite normal form beginning $(y_1, y_2, \ldots, y_p, \ldots)$ then $x \in U$.  To prove this, consider such an element $x$.  Since $g \in y$ we know that $g \prec y_1 \cdots y_r$ for some $r$.  By \cite[Lemma 4.10]{Mi} we have that $\nu(g) \le r$ and $g \prec y_1 \cdots y_{\nu(g)} \prec y_1 \cdots y_p$, and hence that $g \in x$, i.e. that $x \in Z(g)$.  If for some $i$ we have $x \in Z(h_i)$ then the same argument implies that $h_i \prec y_1 \cdots y_p$ and hence that $h_i \in y$, contradicting the fact that $y \not\in Z(h_i)$.  Therefore $x \in Z(g) \setminus \bigcup_{i=1}^n Z(h_i) = U$.  Next we claim that there is an element $\varepsilon \in P \setminus \Delta P$ such that $Z(\varepsilon) \subseteq U$.  To prove this, we first give a preliminary result. Consider an element $c \in P \setminus \Delta P$ with normal form $(c_1, \ldots, c_q)$, and suppose that for some $t \in S$ we have $\cL(c_q) = S \setminus \{ t \}$. We claim that every element of $Z(c)$ has infinite normal form $(c_1, \ldots, c_{q-1}, \ldots)$.  To see this, let $w \in P$ be such that $cw \in P \setminus \Delta P$. By \cite[Lemma 4.6]{Mi} the normal form of $c_q w$ is $(c_q w_1', w_1'' w_2', \ldots)$.  Then $S \supsetneq \cL(c_q w_1') \supseteq \cL(c_q) = S \setminus \{ t \}$, and hence $\cL(c_q w_1') = S \setminus \{ t \}$.  Since $(c_1, \ldots, c_q)$ is a normal form we must have $\cR(c_{q-1}) = S \setminus \{ t \}$, and hence the normal form of $cw$ is $(c_1, \ldots, c_{q-1}, c_q w_1', \ldots)$. Let $x \in Z(c)$ with infinite normal form $(x_1, x_2, \ldots)$. Then for all large enough $r$ we have $x_1 \cdots x_r \in cP \setminus \Delta P$, and hence $x_i = c_i$ for $i < q$. Now we return to the  proof of the existence of the element $\varepsilon$.  By Proposition \ref{prop:singleton-subset equivalence} there is $e \in P$ with normal form $(e_1, \ldots, e_k)$ such that $\cR(y_p) \supseteq \cL(e_1)$ and $\cR(e_k) = S \setminus \{ t \}$ for some choice of $t \in S$.  Let $\varepsilon = y_1 \cdots y_p e \Delta_{S \setminus \{ t \}}$.  Let $x \in Z(\varepsilon)$.  By the previous claim, $x$ has infinite normal form $(y_1, \ldots, y_p, e_1, \ldots, e_k, \ldots)$.  By the first claim above we have $x \in U$.  Thus $Z(\varepsilon) \subseteq U$.  We will write the normal form of $\varepsilon$ as $(\varepsilon_1, \ldots, \varepsilon_m)$. It follows almost immediately that the restriction to $X_0$ is minimal: if $x \in X_0$, then $\varepsilon x \in Z(\varepsilon) \subseteq U$. Now we prove that the restriction to $X_0$ is locally contractive. We must show that $U$ contains an open set $V$ such that $\overline{V} \subseteq U$, and such that there is $b \in P$ with $b \overline{V} \subsetneq V$.  Since $S$ has cardinality at least three we may choose $u \in S$ such that $u \not= t$ and $\cL(\varepsilon_1) \not= \{u\}$.  Then $u \in \cR(\varepsilon_m)$.  By Proposition \ref{prop:singleton-subset equivalence} there is $a \in P$ with normal form $(a_1,\ldots,a_\ell)$ such that $\cL(a_1) = \{u\}$ and $\cR(a_\ell) \supseteq \cL(\varepsilon)$. Then elements of $\varepsilon  u a Z(\varepsilon)$ have infinite normal forms $(\varepsilon_1,\ldots, \varepsilon_m, u, a_1, \ldots, a_\ell, \varepsilon_1, \ldots, \varepsilon_{m-1}, \ldots)$, while elements of $\varepsilon u^2 a Z(\varepsilon)$ have infinite normal forms $(\varepsilon_1,\ldots, \varepsilon_m, u, u, a_1, \ldots, a_\ell, \varepsilon_1, \ldots, \varepsilon_{m-1}, \ldots)$. This follows from the preliminary result above, applied to $c = \varepsilon u a \varepsilon$ and $c = \varepsilon u^2 a \varepsilon$. It follows that $\varepsilon u a Z(\varepsilon)$ and $\varepsilon u^2 a Z(\varepsilon)$ are disjoint subsets of $Z(\varepsilon)$.  Therefore letting $b = \varepsilon u a$ we have that $b Z(\varepsilon) \subsetneq Z(\varepsilon)$.

\quad Now we show that the restriction to $X_0$ is topologically principal.  We must find an element $z \in U$ with trivial isotropy. Let $u \in \cR(y_p)$, and let $t \in S$ with $m_{ut} > 2$.  Let $(i_1,i_2,\ldots)$ be an aperiodic sequence in $\prod_1^\infty \{1,2\}$.  In the following, for $s \in S$ we will write $(s)^i$ for the normal form $(s,s,\ldots,s)$ with $i$ terms. Let $z \in X_0$ have normal form $(y_1, \ldots, y_p, u, ut, (t)^{i_1}, tu, ut, (t)^{i_2}, tu, ut, \ldots)$.  Then $z \in U$, by the first claim at the beginning of the proof of the Theorem. Assume $g \in G$ satisfies $gz = z$. Let $z$ have normal form $(z_1, z_2, z_3, \dotsc)$. $z \in U_{g^{-1}}$ implies that $g = p q^{-1}$ and $q \in [z_1 \dotsm z_l]$ for some $l$, i.e., there must exist $x \in P$ with $qx = z_1 \dotsm z_l$. By replacing $q$ by $qx$ and $p$ by $px$, we may assume without loss of generality that $q = z_1 \dotsm z_l$. Then $q^{-1}z$ lies in $X_0$ and has normal form $(z_{l+1}, z_{l+2}, z_{l+3}, \dotsc) = (z_{l+1}, \dotsc, z_n, ut, (t)^{i_j}, tu, ut, (t)^{i_{j+1}}, \dotsc)$. Let $d$ be such that $d (q^{-1}z) \in X_0$.  Let $d$ have normal form $(d_1, \ldots, d_c)$.  Let us consider the product $d_c z_{l+1} \cdots z_n (ut) t^{i_j} (tu)$.  By \cite[Lemma 4.6]{Mi} we can write $z_m = z_m' z_m''$ so that $d_c z_{l+1} \cdots z_n (ut) t^{i_j} (tu)$ has normal form either $(d_c z_{l+1}', z_{l+1}'' z_{l+2}', \ldots, z_{n-1}'' z_n', z_n'', ut, (t)^{i_j}, tu)$ or $(d_c z_{l+1}', z_{l+1}'' z_{l+2}', \ldots, z_{n-1}'' z_n', z_n'' u, t, (t)^{i_j}, tu)$ or $(d_c z_{l+1}', z_{l+1}'' z_{l+2}', \ldots, z_{n-1}'' z_n', z_n'' ut, (t)^{i_j}, tu)$. Similarly, the product $d_{c-1} d_c z_{l+1} \cdots z_n (ut) t^{i_j} (tu) (ut) t^{i_{j+1}} (tu)$ has normal form $( \cdots, ut, (t)^{i_{j+1}}, tu)$.  Continuing we see that for all large enough $s$ we have that $d z_{l+1} \cdots z_s$ has normal form $( \cdots, ut, (t)^{i_k}, tu, ut, (t)^{i_{k+1}}, tu, \ldots)$. So there exists a positive integer $r$ with $(z_{r+1}, z_{r+2}, z_{r+3}, \ldots) = (ut, (t)^{i_k}, tu, \ldots)$ and such that for all $s \geq r$, the normal form of $d z_{l+1} \dotsm z_s$ is given by the normal form of $d z_{l+1} \dotsm z_r$ followed by $z_{r+1}, z_{r+2}, z_{r+3}, \ldots, z_s$. So the normal form of $dz$ is given by the normal form of $d z_{l+1} \dotsm z_r$ followed by $ut, (t)^{i_k}, tu, ut, (t)^{i_{k+1}}, tu, \ldots$. Now apply this to $d=p$ and $d=q$. Then the normal form of $gz = p (q^{-1}z)$ is given by the normal form of $p z_{l+1} \dotsm z_r$ followed by $ut, (t)^{i_k}, tu, ut, (t)^{i_{k+1}}, tu, \ldots$. Similarly, the normal form of $z = q (q^{-1}z)$ is given by the normal form of $q z_{l+1} \dotsm z_r$ followed by $ut, (t)^{i_k}, tu, ut, (t)^{i_{k+1}}, tu, \ldots$. Since $gz = z$, these two normal forms must coincide. This --- because of aperiodicity of $(i_h)$ --- implies that the normal forms of $p z_{l+1} \dotsm z_r$ and $q z_{l+1} \dotsm z_r$ must coincide. Hence $p=q$, and thus $g=1$.

Therefore $z$ has trivial isotropy.  Now, simplicity follows from \cite[Proposition II.4.6]{Re}, and pure infiniteness from \cite[Proposition 2.4]{An}.
\eproof
\setlength{\parindent}{0cm} \setlength{\parskip}{0.5cm}

\section{Exact sequences in K-theory}
\label{sec:exseq-K}

Let $P = \spkl{a,b \, \vert \, u=v}^+$ be a right LCM monoid as in \S~\ref{sec:Ore}, and assume that we are in case I or case II of Theorem~\ref{Thm:OreOneRel:graph}. Let $G = \spkl{a,b \, \vert \, u=v}$ be the group given by the same presentation as $P$. Suppose $A$ is a C*-algebra and $\gamma: \: G \curvearrowright A$ is a $G$-action on $A$. The partial action $G \curvearrowright \Omega$ described in \S~\ref{SgpC} induces a partial action $G \curvearrowright C(\Omega)$, which in turn gives rise to the diagonal partial action $G \curvearrowright C(\Omega) \otimes A$. Recall that $\ti{\Omega} = \Omega \setminus \gekl{\infty}$ is $G$-invariant, so that we obtain a partial action $G \curvearrowright C_0(\ti{\Omega}) \otimes A$ by restriction. As $G$ is exact \cite{Gue}, we obtain an exact sequence
\begin{equation}
\label{exseq:POinfty}
0 \to (C_0(\ti{\Omega}) \otimes A) \rtimes_r G \to (C(\Omega) \otimes A) \rtimes_r G \to A \rtimes_r G \to 0. 
\end{equation}
In case $P \subseteq G$ is Toeplitz, we know how to compute K-theory for $(C(\Omega) \otimes A) \rtimes_r G$ because of \cite{CEL1,CEL2}. Therefore, our goal now is to deduce a six term exact sequence in K-theory which allows us to compute K-theory for $(C_0(\ti{\Omega}) \otimes A) \rtimes_r G$. This will then enable us to determine the K-theory of $A \rtimes_r G$.

Let us explain our strategy. Recall that we defined $X = \menge{\chi \in \Omega}{\chi(wP) = 0}$, with $w \in P$ as in Lemma~\ref{LEM:w}. As $X$ is a compact open subspace of $\ti{\Omega}$ meeting every $G$-orbit, the projection $\boldsymbol{p} \defeq (1_X \otimes 1)$ (where $1 = 1_{M(A)} \in M(A)$) is a full projection, so that $(C_0(\ti{\Omega}) \otimes A) \rtimes_r G \sim_M \boldsymbol{p} \rukl{(C_0(\ti{\Omega}) \otimes A) \rtimes_r G} \boldsymbol{p}$. Moreover, we have an exact sequence
\begin{equation}
\label{exseq:POtiO}
  0 \to \boldsymbol{p} \rukl{(C_0(P) \otimes A) \rtimes_r G} \boldsymbol{p} \to \boldsymbol{p} \rukl{(C_0(\ti{\Omega}) \otimes A) \rtimes_r G} \boldsymbol{p} \to \boldsymbol{p} \rukl{(C_0(\ti{\Omega}_{\infty}) \otimes A) \rtimes_r G} \boldsymbol{p} \to 0.
\end{equation}
Note that we have an isomorphism
\begin{equation}
\label{KA=P}
  \cK(\ell^2(P \cap X)) \otimes A \cong \boldsymbol{p} \rukl{(C_0(P) \otimes A) \rtimes_r G} \boldsymbol{p}, \, e_{x,y} \otimes a \ma u_x (1_{\gekl{e}} \otimes a) u_{y^{-1}}.
\end{equation}
Here $P \cap X = \menge{x \in P}{x \notin wP}$, $e_{x,y}$ is the rank one operator corresponding to the basis vectors $\delta_x$ and $\delta_y$ in $\ell^2(P \cap X)$, and $u_x, u_y$ are the canonical partial isometries implementing the partial action in the crossed product. Moreover, using the graph models from \S~\ref{sec:Ore}, it is possible to compute K-theory for $\boldsymbol{p} \rukl{(C_0(\ti{\Omega}_{\infty}) \otimes A) \rtimes_r G} \boldsymbol{p}$. Hence the K-theoretic six term exact sequence provided by \eqref{exseq:POtiO} already provides a way to compute K-theory for $\boldsymbol{p} \rukl{(C_0(\ti{\Omega}) \otimes A) \rtimes_r G} \boldsymbol{p}$. However, since the C*-algebra we are interested in appears in the middle of the sequence, we would have to work out boundary maps, which can be complicated. Therefore, we present an alternative approach. The idea is to compare \eqref{exseq:POtiO} with the canonical exact sequence for the twisted graph algebra of a graph model $\cE$, where we twist $C^*(\cE)$ by $A$ using the $G$-action $\gamma$. In the following, let us make this precise, i.e., let us construct twisted graph algebras and their canonical extensions. Our construction generalizes the one in \cite{Cun81}.
\setlength{\parindent}{0cm} \setlength{\parskip}{0cm}

\subsection{Graph algebras twisted by coefficients}
\label{twGraphAlg}

Let $\cE$ be a graph with set of vertices $\cE^0$, set of edges $\cE^1$ and range and source maps $r$, $s$. We assume that $\cE$ is finite without sources because this will be the case in our applications. However, everything in this subsection also works for general graphs (with appropriate modifications).
\setlength{\parindent}{0cm} \setlength{\parskip}{0.5cm}

Let $F$ be the free group on the set of edges $\cE^1$ of $\cE$. Let $S_{\cE}$ be the inverse semigroup attached to $\cE$ as in \cite[\S~3.2]{Li17}, and $E_{\cE}$ its semilattice of idempotents. Set $X_{\cE} \defeq \widehat{E_{\cE}}$. $X_{\cE}$ can be identified with the set of finite and infinite paths in $\cE$. Let $W_{\cE}$ be the subspace of finite paths and $Y_{\cE}$ the subspace of infinite paths. Note that in our case, $Y_{\cE} = \partial \cE$. As explained in \cite[\S~3.2]{Li17}, there is a canonical partial action $F \curvearrowright X_{\cE}$ such that $C(X_{\cE}) \rtimes_r F$ is canonically isomorphic to the Toeplitz algebra of $\cE$ and $C(Y_{\cE}) \rtimes_r F$ is canonically isomorphic to the graph algebra of $\cE$. Moreover, we have an exact sequence $0 \to C_0(W_{\cE}) \rtimes_r F \to C(X_{\cE}) \rtimes_r F \to C(Y_{\cE}) \rtimes_r F \to 0$.

Now assume that $A$ is a C*-algebra with an $F$-action $\gamma$. Then we can form the diagonal partial action $F \curvearrowright C(X_{\cE}) \otimes A$. By exactness of $F$, we obtain an exact sequence
\begin{equation}
\label{exseq:graphA}
0 \to (C_0(W_{\cE}) \otimes A) \rtimes_r F \to (C(X_{\cE}) \otimes A) \rtimes_r F \to (C(Y_{\cE}) \otimes A) \rtimes_r F \to 0.
\end{equation}
Let $p_v$, $v \in \cE^0$ and $s_e$, $e \in \cE^1$, be the canonical generators (projections and partial isometries) of the Toeplitz algebra of $\cE$, so that we may consider the elements $p_v$ and $s_e s_e^*$ in $C(X_{\cE})$. A combination of \cite[Corollary~3.14]{CEL2} and \cite[\S~6.1]{LNII} (see also \cite{LNI}) yields that the homomorphisms
\begin{align}
\label{OplusA-W}
  & \bigoplus_{\cE^0} A \to (C_0(W_{\cE}) \otimes A) \rtimes_r F, \, (a_v)_v \ma \sum_v \Big( p_v - \sum_{e \in r^{-1}(v)} s_e s_e^* \Big) \otimes a_v,\\
\label{OplusA-X}
& \bigoplus_{\cE^0} A \to (C(X_{\cE}) \otimes A) \rtimes_r F, \, (a_v)_v \ma \sum_v p_v \otimes a_v 
\end{align}
induce isomorphisms in K-theory, and that these $K_*$-isomorphisms fit into a commutative diagram
\begin{equation}
\label{CD:KtwGraphAlg}
\begin{tikzcd}
\bigoplus_{\cE^0} K_*(A) \ar["\id - M"]{r} \ar["\cong"']{d} & \bigoplus_{\cE^0} K_*(A) \ar["\cong"]{d} \\
K_*((C_0(W_{\cE}) \otimes A) \rtimes_r F) \ar[r] & K_*((C(X_{\cE}) \otimes A) \rtimes_r F)
\end{tikzcd}
\end{equation}
where the lower horizontal map is the one in \eqref{exseq:graphA} and $M$ is the $\cE^0 \times \cE^0$-matrix $M = (M_{w,v})_{w,v}$ with entries $M_{w,v} \in \End(K_*(A))$ given by $M_{w,v} = \sum_{e \in s^{-1}(w) \cap r^{-1}(v)} (\gamma_e)_*^{-1}$.

In order to compare the exact sequence \eqref{exseq:graphA} with \eqref{exseq:POtiO}, we need to establish a universal property for $(C(X_{\cE}) \otimes A) \rtimes_r F$. Let us assume that $A$ is unital and treat the non-unital case later. Let $\cT_A$ be the universal C*-algebra which comes with a homomorphism $\iota: \: A \to \cT_A$ such that $\cT_A$ is generated by projections $p_v$, $v \in \cE^0$, partial isometries $s_e$, $e \in \cE^1$, and $\iota(A)$, subject to the relations that $p_v$, $v \in \cE^0$, are pairwise orthogonal, $s_e^* s_e = p_{s(e)}$, $p_v \geq \sum_{e \in r^{-1}(v)} s_e s_e^*$ and $s_e \iota(a) = \iota(\gamma_e(a)) s_e$ for all $e \in \cE^1$, $v \in \cE^0$ and $a \in A$. This means that whenever we have a C*-algebra $T$, a homomorphism $i: \: A \to T$, projections $q_v$, $v \in \cE^0$, and partial isometries $t_e$, $e \in \cE^1$, such that $q_v$, $t_e$ and $i(a)$ satisfy analogous relations as the generators of $\cT_A$, then there is a (unique) homomorphism $\cT_A \to T$ sending $p_v$ to $q_v$, $s_e$ to $t_e$ and $\iota(a)$ to $i(a)$. Hence, by universal property of $\cT_A$, there is a homomorphism
\begin{equation}
\label{T_A-->}
  \cT_A \to (C(X_{\cE}) \otimes A) \rtimes_r F, \, p_v \ma 1_v \otimes 1_A \in C(X_{\cE}) \otimes A, \, s_e \ma u_e, \, \iota(a) \ma 1 \otimes a.
\end{equation}
Here $1_v$ is the characteristic function of the subspace of $X_{\cE}$ of all finite and infinite paths with range $v$. We see that $\iota: \: A \to \cT_A$ must be injective, because $A \to (C(X_{\cE}) \otimes A) \rtimes_r F, a \ma 1 \otimes a$ is. So we will from now on identify $A$ with $\iota(A)$ and write $a$ instead of $\iota(a)$. Our goal is to show that the homomorphism in \eqref{T_A-->} is an isomorphism.

For a finite path $\mu \in \cE^*$, where $\mu = \mu_1 \dotsm \mu_m$  with $\mu_i \in \cE^1$, let $s_{\mu} = s_{\mu_1} \dotsm s_{\mu_m}$. First of all, it is easy to see that $\cT_A = \clspan \big( \menge{s_{\mu} s_{\nu}^* a}{a \in A, \, \mu, \nu \in \cE^*, \, s(\mu) = s(\nu)} \big)$. Moreover, by universal property of $\cT_A$, for every $z \in \Tz$ there is an isomorphism $\delta_z: \: \cT_A \cong \cT_A$ given by $\delta_z(a) = a$, $\delta_z(p_v) = p_v$ and $\delta_z(s_e) = z s_e$. The following is easy to see (as in the case of ordinary graph algebras without coefficients):
\blemma
$\theta(x) \defeq \int_{\Tz} \delta_z(x) dz$ defines a faithful conditional expectation $\cT_A \onto \cT_A^{\delta}$ determined by $\theta(s_{\mu} s_{\nu}^* a) = \delta_{l(\mu), \, l(\nu)} s_{\mu} s_{\nu}^* a$.
\elemma
\setlength{\parindent}{0cm} \setlength{\parskip}{0cm}

Now let $\cT$ be the classical Toeplitz algebra of $\cE$, so that $\cT$ is the universal C*-algebra generated by projections $p_v$, $v \in \cE^0$, and partial isometries $s_e$, $e \in \cE^1$, satisfying the same relations as the generators of $\cT_A$. Again, we have a canonical $\Tz$-action $\delta$ on $\cT$, with fix point algebra $\cT^{\delta}$.
\blemma
There is an isomorphism $\cT^{\delta} \otimes A \to \cT_A^{\delta}$ sending $s_{\mu} s_{\nu}^* \otimes a$ to $s_{\mu} a s_{\nu}^*$.
\elemma

\bproof
Let
\begin{align*}
\cT_l^{\delta} &= \clspan \menge{s_{\mu} s_{\nu}^*}{s(\mu) = s(\nu), \, l(\mu) = l(\nu) \leq l},\\
\cT_{A,l}^{\delta} &= \clspan \menge{s_{\mu} s_{\nu}^* a}{a \in A, \, s(\mu) = s(\nu), \, l(\mu) = l(\nu) \leq l}. 
\end{align*}
Then $\cT^{\delta} \otimes A \cong \ilim_l \cT_l^{\delta} \otimes A$ and $\cT_A^{\delta} \cong \ilim_l \cT_{A,l}^{\delta}$. So it suffices to show that for every $l$, we have an isomorphism $\cT_l^{\delta} \otimes A \to \cT_{A,l}^{\delta}, \, s_{\mu} s_{\nu}^* \otimes a \ma s_{\mu} a s_{\nu}^*$. For fixed $l$, define $\epsilon_{\mu, \nu} = s_{\mu} s_{\nu}^* - \sum_{e \in r^{-1}(s(\mu))} s_{\mu e} s_{\nu e}^*$ if $l(\mu) = l(\nu) \leq l-1$ and  $\epsilon_{\mu, \nu} = s_{\mu} s_{\nu}^*$ if $l(\mu) = l(\nu) = l$. Then we have $\epsilon_{\kappa, \lambda}^* = \epsilon_{\lambda, \kappa}$ and $\epsilon_{\kappa, \lambda} \epsilon_{\mu, \nu} = \delta_{\lambda, \mu} \epsilon_{\kappa,\nu}$. Thus $\cT_l^{\delta} \otimes A$ is the universal C*-algebra generated by $\epsilon_{\mu,\nu} \otimes a$ for $\mu, \nu \in \cE^*$ with $s(\mu) = s(\nu)$ and $l(\mu) = l(\nu) \leq l$, and for $a \in A$, subject to the relations $(\epsilon_{\kappa, \lambda} \otimes a)^* = \epsilon_{\lambda, \kappa} \otimes a^*$ and $(\epsilon_{\kappa, \lambda} \otimes a) (\epsilon_{\mu, \nu} \otimes b) = \delta_{\lambda, \mu} \epsilon_{\kappa,\nu} \otimes ab$. This shows existence of the homomorphism $\cT_l^{\delta} \otimes A \to \cT_{A,l}^{\delta}, \, s_{\mu} s_{\nu}^* \otimes a \ma s_{\mu} a s_{\nu}^*$. It is clearly surjective. To prove injectivity, it suffices to show that $A \to \cT_{A,l}^{\delta}, \, a \ma s_{\mu} a s_{\nu}^* - \sum_{e \in r^{-1}(s(\mu))} s_{\mu e} a s_{\nu e}^*$ (for $l(\mu) = l(\nu) \leq l-1$) and $A \to \cT_{A,l}^{\delta}, \, a \ma s_{\mu} a s_{\nu}^*$ (for $l(\mu) = l(\nu) = l$) are injective. But this is easy to see by considering their compositions $A \to \cT_{A,l}^{\delta} \to \cT_A^{\delta} \to \cT_A \to (C(X_{\cE}) \otimes A) \rtimes_r F$.
\eproof

\bcor
The homomorphism in \eqref{T_A-->} is an isomorphism.
\ecor
\bproof
Consider the composition $\Theta: \: \cT_A \overset{\theta}{\lori} \cT_A^{\delta} \cong \cT^{\delta} \otimes A \to C(X_{\cE}) \otimes A$, where the last map is the tensor product of the canonical faithful conditional expectation $\cT^{\delta} \onto C(X_{\cE})$ with $\id_A$. $\Theta$ is a faithful conditional expectation determined by $s_{\mu} a s_{\nu}^* \ma \delta_{\mu,\nu} 1_{\mu} \otimes a$, where $1_{\mu}$ is the characteristic function of the subspace of $X_{\cE}$ of all finite and infinite paths starting with $\mu$. It is now easy to see that we have a commutative diagram
\begin{equation*}
\begin{tikzcd}
\cT_A \ar{r} \ar["\Theta"']{d} & (C(X_{\cE}) \otimes A) \rtimes_r F \ar{d} \\
C(X_{\cE}) \otimes A \ar[equal]{r} & C(X_{\cE}) \otimes A
\end{tikzcd}
\end{equation*}
where the right vertical map is the canonical faithful conditional expectation and the upper horizontal map is the homomorphism in \eqref{T_A-->}. It follows that the map in \eqref{T_A-->}, which is clearly surjective, is also injective.
\eproof

\subsection{Comparing exact sequences in K-theory}

Now we return to the setting described at the beginning of \S~\ref{sec:exseq-K}. Let $\cE$ be a graph model for $(G \ltimes \ti{\Omega}_{\infty}) \, \vert \, Y$, where $Y = \Omega_{\infty} \cap X = \ti{\Omega}_{\infty} \cap X$. The $G$-action $\gamma$ induces an action of the free group $F$ on $\cE^1$ on $A$ by letting $e \in \cE^1$ act via $\gamma_{\boldsymbol{\sigma}(e)}$. 
\setlength{\parindent}{0cm} \setlength{\parskip}{0.5cm}

First let $A$ be unital. By universal property of $(C(X_{\cE}) \otimes A) \rtimes_r F$, as established in \S~\ref{twGraphAlg}, we have a homomorphism
\begin{equation*}
  (C(X_{\cE}) \otimes A) \rtimes_r F \to \boldsymbol{p} \rukl{(C_0(\ti{\Omega}) \otimes A) \rtimes_r G} \boldsymbol{p}, \, p_v \ma 1_X 1_{vP} \otimes 1_A, \, u_e \ma \boldsymbol{p} u_{\boldsymbol{\sigma}(e)} (1_{s(e)P} \otimes 1_A) \boldsymbol{p}, \, a \ma 1_X \otimes a.
\end{equation*}
Note that in general, this homomorphism is neither injective nor surjective. By our choice of $\cE$, this homomorphism fits into a commutative diagram
\begin{equation}
\label{compare:exseq}
\begin{tikzcd}
0 \ar{r} & (C(W_{\cE}) \otimes A) \rtimes_r F \ar{r} \ar{d} & (C(X_{\cE}) \otimes A) \rtimes_r F \ar{r} \ar{d} & (C(Y_{\cE}) \otimes A) \rtimes_r F \ar{r} \ar["\cong"]{d} & 0 \\
0 \ar{r} & \boldsymbol{p} \rukl{(C_0(P) \otimes A) \rtimes_r G} \boldsymbol{p} \ar{r} & \boldsymbol{p} \rukl{(C_0(\ti{\Omega}) \otimes A) \rtimes_r G} \boldsymbol{p} \ar{r} & \boldsymbol{p} \rukl{(C_0(\ti{\Omega}_{\infty}) \otimes A) \rtimes_r G} \boldsymbol{p} \ar{r} & 0 
\end{tikzcd}
\end{equation}
where the third vertical map is an isomorphism because $\cE$ is a graph model for $(G \ltimes \ti{\Omega}_{\infty}) \, \vert \, Y$.

For non-unital $A$, just apply our constructions above to the unitalization $A \ti{\ }$ and then restrict all the maps to the crossed products with $A$ in place of $A \ti{\ }$.

The isomorphism \eqref{KA=P} implies that the homomorphism $A \to \boldsymbol{p} \rukl{(C_0(P) \otimes A) \rtimes_r G} \boldsymbol{p}, \, a \ma 1_{\gekl{e}} \otimes a$ induces an $K_*$-isomorphism, which together with the $K_*$-isomorphism induced by \eqref{OplusA-W} fit into a commutative diagram
\begin{equation}
\label{CD:phi}
\begin{tikzcd}
\bigoplus_{\cE^0} K_*(A) \ar["\cong"]{r} \ar["\phi"']{d} & K_*((C(W_{\cE}) \otimes A) \rtimes_r F) \ar{d} \\
K_*(A) \ar["\cong"]{r} & K_*(\boldsymbol{p} \rukl{(C_0(P) \otimes A) \rtimes_r G} \boldsymbol{p}) 
\end{tikzcd}
\end{equation}
where the right vertical map is induced by the first vertical map in \eqref{compare:exseq}, and $\phi = \sum_{v \in \cE^0} \phi_v$ with
\begin{equation*}
  \phi_v = \sum_{x \in \rukl{vP \, \setminus \, \bigcup_{e \in r^{-1}(v)} \boldsymbol{\sigma}(e) s(e) P} \cap X} (\gamma_x)_*^{-1}.
\end{equation*}

Let us now assume that our graph model $\cE$ is such that the map $\phi$ in \eqref{CD:phi} is surjective. This is the only requirement we need on our graph models. In \S~\ref{ss:KexseqATK}, we will construct concrete graph models with this property.

Here is a general lemma which we will apply in our special situation:
\blemma
\label{LEM:HomAlg}
Suppose we have the following commutative diagram of abelian groups with exact rows:
\begin{equation}
\label{CDexrows}
\begin{tikzcd}
\dotso \ar["p_{i-1}"]{r} & \bar{G}_{i-1} \ar["\partial_{i-1}"]{r} \ar["\cong"', "\psi_{i-1}"]{d} & \check{G}_i \ar["j_i"]{r} \ar[twoheadrightarrow, "\phi_i"]{d} & G_i \ar["p_i"]{r} \ar["\pi_i"]{d} & \bar{G}_i \ar["\partial_i"]{r} \ar["\cong"', "\psi_i"]{d} & \check{G}_{i+1} \ar["j_{i+1}"]{r} \ar[twoheadrightarrow, "\phi_{i+1}"]{d} & \dotso \\
\dotso \ar["q_{i-1}"]{r} & \bar{H}_{i-1} \ar["\epsilon_{i-1}"]{r} & \check{H}_i \ar["k_i"]{r} & H_i \ar["q_i"]{r} & \bar{H}_i \ar["\epsilon_i"]{r} & \check{H}_{i+1} \ar["k_{i+1}"]{r} & \dotso
\end{tikzcd}
\end{equation}
Then we have an exact sequence
\begin{equation}
\label{CD->exseq}
\begin{tikzcd}
\dotso \ar["\partial_{i-1} \circ \psi_{i-1}^{-1} \circ q_{i-1}"]{rr} & & \ker(\phi_i) \ar["j_i"]{r} & G_i \ar["\pi_i"]{r} & H_i \ar["\partial_i \circ \psi_i^{-1} \circ q_i"]{rr} & & \ker(\phi_{i+1}) \ar["j_{i+1}"]{r} & G_{i+1} \ar["\pi_{i+1}"]{r} & \dotso
\end{tikzcd}
\end{equation}
\elemma
\setlength{\parindent}{0cm} \setlength{\parskip}{0cm}
This lemma is most likely well-known, but for the sake of completeness, we include a proof.

\bproof
Since $\phi_{i+1} \circ \partial_i \circ \psi_i^{-1} \circ q_i = \epsilon_i \circ \psi_i \circ \psi_i^{-1} \circ q_i  = \epsilon_i \circ q_i = 0$, we have $\img(\partial_i \circ \psi_i^{-1} \circ q_i) \subseteq \ker(\phi_{i+1})$. Moreover, composition of consecutive homomorphisms gives zero since $\partial_i \circ \psi_i^{-1} \circ q_i \circ \pi_i = \partial_i \circ \psi_i^{-1} \circ \psi_i \circ p_i = \partial_i \circ p_i = 0$. 
\setlength{\parindent}{0.5cm} \setlength{\parskip}{0cm}

Let us verify exactness at $G_i$: Take $z \in G_i$ with $\pi_i(z) = 0$. Then $\psi_i(p_i(z)) = q_i(\pi_i(z)) = 0$, so that $p_i(z) = 0$. So there exists $y \in \check{G}_i$ with $j_i(y) = z$. Since $k_i(\phi_i(y)) = \pi_i(j_i(y)) = 0$, there is $x \in \bar{H}_{i-1}$ with $\epsilon_{i-1}(x) = \phi_i(y)$. Set $w \defeq \psi_{i-1}^{-1}(x)$ and $v \defeq y - \partial_{i-1}(w)$. Then $\phi_i(v) = \phi_i(y) - \phi_i(\partial_{i-1}(w)) = \phi_i(y) - \epsilon_{i-1}(\psi_{i-1}(w)) = \phi_i(y) - \epsilon_{i-1}(x) = 0$ and $j_i(v) = j_i(y) = z$.

Now we show exactness at $H_i$: Let $z \in H_i$ satisfy $\partial_i(\psi_i^{-1}(q_i(z))) = 0$. Then there exists $y \in G_i$ with $p_i(y) = \psi_i^{-1}(q_i(z))$. Since $q_i(\pi_i(y)) = \psi_i(p_i(y)) = q_i(z)$, there is $x \in \check{H}_i$ with $k_i(x) = z - \pi_i(y)$. As $\phi_i$ is surjective, there is $w \in \check{G}_i$ with $\phi_i(w) = x$. Let $v \defeq y + j_i(w)$. Then $\pi_i(v) = \pi_i(y) + \pi_i(j_i(w)) = \pi_i(y) + k_i(\phi_i(w)) = \pi_i(y) + k_i(x) = z$.

Finally, for exactness at $\ker(\phi_{i+1})$, let $z \in \ker(\phi_{i+1})$ satisfy $j_{i+1}(z) = 0$. Then there is $y \in \bar{G}_i$ with $\partial_i(y) = z$. Since $\epsilon_i(\psi_i(y)) = \phi_{i+1}(\partial_i(y)) = \phi_{i+1}(z) = 0$, there is $x \in H_i$ with $q_i(x) = \psi_i(y)$. Then $\partial_i(\psi_i^{-1}(q_i(x))) = \partial_i(\psi_i^{-1}(\psi_i(y))) = \partial_i(y) = z$.
\eproof
\setlength{\parindent}{0cm} \setlength{\parskip}{0cm}

Let us now apply Lemma~\ref{LEM:HomAlg} to the diagram induced by \eqref{compare:exseq} in K-theory, where $\cE$ is a graph model for $(G \ltimes \ti{\Omega}_{\infty}) \, \vert \, Y$ such that the map $\phi$ in \eqref{CD:phi} is surjective. Let us write $I \defeq \boldsymbol{p} \rukl{(C_0(\ti{\Omega}) \otimes A) \rtimes_r G} \boldsymbol{p}$. Plugging \eqref{CD:KtwGraphAlg} and \eqref{CD:phi} into the K-theory diagram for \eqref{compare:exseq}, we obtain the following commutative diagram with exact rows
\begin{equation}
\label{CD:Kcompare}
\begin{tikzcd}
\dotso \ar["\partial_{i-1}"]{r} & \bigoplus_{\cE^0} K_i(A) \ar["j_i"]{r} \ar[twoheadrightarrow, "\phi_i"]{d} & \bigoplus_{\cE^0} K_i(A) \ar{r} \ar["\pi_i"]{d} & K_i((C(Y_{\cE}) \otimes A) \rtimes_r F) \ar["\partial_i"]{r} \ar["\cong"', "\psi_i"]{d} & \dotso \\
\dotso \ar{r} & K_i(A) \ar{r} & K_i(I) \ar["q_i"]{r} & K_i(\boldsymbol{p} \rukl{(C_0(\ti{\Omega}_{\infty}) \otimes A) \rtimes_r G} \boldsymbol{p}) \ar{r} & \dotso
\end{tikzcd}
\end{equation}
Applying Lemma~\ref{LEM:HomAlg}, we obtain an exact sequence
\begin{equation}
\label{exseq-K1}
\begin{tikzcd}
\dotso \ar{r} & \ker(\phi_i) \ar["j_i"]{r} & \bigoplus_{\cE^0} K_i(A) \ar["\pi_i"]{r} & K_i(I) \ar["\partial_i \circ \psi_i^{-1} \circ q_i"]{rr} & & \ker(\phi_{i+1}) \ar{r} & \dotso
\end{tikzcd}
\end{equation}
Furthermore, write $J \defeq (C_0(\ti{\Omega}) \otimes A) \rtimes_r G$. Since $\boldsymbol{p}$ is a full projection in $J$, the canonical inclusion $I \into J$ induces a $K_*$-isomorphism. Let us assume that $P \subseteq G$ is Toeplitz, which in our case amounts to saying that $P \subseteq G$ is quasi-lattice ordered as in Remark~\ref{Rem:Toe=qlo-K}. Then \cite[Corollary~3.14]{CEL2} implies that the homomorphism $A \to (C(\Omega) \otimes A) \rtimes_r G, \, a \ma 1 \otimes a$ induces an isomorphism in K-theory (because $G$ satisfies the Baum-Connes conjecture with coefficients by \cite{BBV}). Hence, applying K-theory to \eqref{exseq:POinfty} and using our observations above, we obtain a second exact sequence
\begin{equation}
\label{exseq-K2}
\begin{tikzcd}
\dotso \ar{r} & K_i(I) \ar["\iota_i"]{r} & K_i(A) \ar{r} & K_i(A \rtimes_r G) \ar{r} & \dotso
\end{tikzcd}
\end{equation}
where $K_i(A) \to K_i(A \rtimes_r G)$ is induced by the canonical map and $\iota_i \circ \pi_i$ is given by
\begin{equation}
\label{iota-pi-formula}
  \iota_i \circ \pi_i = \sum_{v \in \cE^0} \rukl{(\gamma_v)_*^{-1} - (\gamma_{v \vee w})_*^{-1}}: \: \bigoplus_{v \in \cE^0} K_i(A) \to K_i(A),
\end{equation}
where $v \vee w \in P$ is the least common multiple of $v$ and $w$ determined by $(v \vee w) P = vP \cap wP$.

\subsection{K-theory exact sequences for Artin-Tits groups of dihedral type and torus knot groups}
\label{ss:KexseqATK}

Let us now consider two classes of monoids and groups which already appeared in Example~\ref{ex:ArtinDihedralTorusKnot}: Let $P = \spkl{a,b \, \vert \, u=v}^+$ and $G = \spkl{a,b \, \vert \, u=v}$, where $u = aba \dotsm$ and $v = bab \dotsm$ with $\ell^*(u) = \ell^*(v)$, or $u = a^p$ and $v = b^q$ for $p, q \geq 2$. In the first case, $P$ and $G$ are called Artin-Tits monoids and Artin-Tits groups of dihedral type $I_2(m)$ with $m = \ell^*(u) = \ell^*(v)$, and in the second case, $P$ and $G$ are called torus knot monoids or torus knot groups of type $(p,q)$. We refer to the first case as the dihedral Artin-Tits case of type $I_2(m)$ and to the second case as the torus knot case of type $(p,q)$. In both of these cases, our goal is to provide concrete graph models for $(G \ltimes \ti{\Omega}_{\infty}) \, \vert \, Y$ such that $\phi$ in \eqref{CD:phi} is surjective, and then to work out the exact sequences \eqref{exseq-K1} and \eqref{exseq-K2} explicitly and further simplify them. Note that by \cite[\S~5.8]{CELY}, $P \subseteq G$ is Toeplitz because $P$ is right reversible (not only left reversible). 
 We start by constructing graph models.

\blemma
\label{LEM:K_ex--GraphModels}
In the dihedral Artin-Tits case of type $I_2(m)$, with $m$ even, construct a graph by setting $\cE^0$ as the following collection of finite words in $\gekl{a,b}^*$ of length at most $m-1$:
\begin{equation}
\label{E0_I_2even}
aa, abaa, ababaa, \dotsc, \underbrace{ab \dotsm ba}_{m-1}, \underbrace{ab \dotsm bb}_{m-1}, \dotsc, ababb, abb, baa, babaa, \dotsc, \underbrace{ba \dotsm baa}_{m-1}, \underbrace{ba \dotsm bab}_{m-1}, \dotsc, babb, bb.
\end{equation}
Let 
\begin{align*}
  \cE^1 = &\menge{(y,x)}{y, x \in \cE^0; \; \ell^*(y) = \ell^*(x) + 1; \; \exists \, \sigma \in \gekl{a,b}: \: y \equiv \sigma x}\\
  \cup &\menge{(y,x)}{y, x \in \cE^0; \; \ell^*(y) = \ell^*(x); \; \exists \, \sigma, \tau \in \gekl{a,b}: \: y \tau \equiv \sigma x \text{ is } u,v \text{-free}}\\
  \cup &\menge{(aa,x)}{x \in \cE^0; \; x \text{ starts with } ab}
  \cup \menge{(bb,x)}{x \in \cE^0; \; x \text{ starts with } ba}.
\end{align*}
Define range and source maps by $r(y,x) = y$ and $s(y,x) = x$, and set $\bm{\sigma}(y,x)$ as the first letter of $y$. Then $(\cE,\bm{\sigma})$ is a graph model for $G \curvearrowright Y$.
In the dihedral Artin-Tits case of type $I_2(m)$, with $m$ odd, construct a graph by setting $\cE^0$ as the following collection of finite words in $\gekl{a,b}^*$ of length at most $m-1$:
\begin{equation}
\label{E0_I_2odd}
aa, abaa, ababaa, \dotsc, \underbrace{ab \dotsm baa}_{m-1}, \underbrace{ab \dotsm ab}_{m-1}, \dotsc, ababb, abb, baa, babaa, \dotsc, \underbrace{ba \dotsm ba}_{m-1}, \underbrace{ba \dotsm bb}_{m-1}, \dotsc, babb, bb.
\end{equation}
Define $\cE^1$, $r$, $s$ and $\bm{\sigma}$ as in the dihedral Artin-Tits case of type $I_2(m)$, with $m$ even. Then $(\cE,\bm{\sigma})$ is a graph model for $G \curvearrowright Y$.
In the torus knot case of type $(p,q)$, construct a graph by setting $\cE^0$ as the following collection of finite words in $\gekl{a,b}^*$:
\begin{equation}
\label{E0_TorusKnot}
ab, a^2b, \dotsc, a^{p-2}b, a^{p-1}, ba, b^2a, \dotsc, b^{q-1}a.
\end{equation}
Let 
\begin{align*}
  \cE^1 = &\menge{(y,x)}{y, x \in \cE^0; \; \ell^*(y) = \ell^*(x) + 1; \; \exists \, \sigma \in \gekl{a,b}: \: y \equiv \sigma x}\\
  \cup &\menge{(y,x)}{y, x \in \cE^0; \; \ell^*(y) = \ell^*(x); \; \exists \, \sigma, \tau \in \gekl{a,b}: \: y \tau \equiv \sigma x \text{ is } u,v \text{-free}}\\
  \cup &\menge{(ab,x)}{x \in \cE^0; \; x \text{ starts with } b}
  \cup \menge{(ba,x)}{x \in \cE^0; \; x \text{ starts with } a}.
\end{align*}
Define $r$, $s$ and $\bm{\sigma}$ as in the dihedral Artin-Tits case of type $I_2(m)$, with $m$ even. Then $(\cE,\bm{\sigma})$ is a graph model for $G \curvearrowright Y$.
\elemma

\bproof
We first treat the dihedral Artin-Tits case. As in case~I of Theorem~\ref{Thm:OreOneRel:graph}, we have a bijection 
$$
  \menge{z \in \gekl{a,b}^{\infty}}{z \text{ is } u,v \text{-free}}, \, z \ma \chi_z.
$$
Moreover, the inverse of $\cE^{\infty} \to \menge{z \in \gekl{a,b}^{\infty}}{z \text{ is } u,v \text{-free}}, \, \mu \ma \bm{\sigma}(\mu)$ is given as follows: Every $u,v$-free infinite word $z \in \gekl{a,b}^{\infty}$ is of the form $z = z_1 \sigma_1^{e_1} z_2 \sigma_2^{e_2} \dotsm$ where $z_i \in \gekl{a,b}^*$ are $u,v$-free words ending on $aa$ or $bb$ and containing no other $aa$ or $bb$ as subwords, $\sigma_i$ is the last letter of $z_i$ and differs from the first letter of $z_{i+1}$, and $e_i \in \Zz_{\geq 0} \cup \gekl{\infty}$ ($e_i = \infty$ means that $z = \dotsm z_i \sigma_i \sigma_i \sigma_i \dotsm$). For each $z_i$, write $z_i = \tau_1 \dotsm \tau_l$ with $\tau_k \in \gekl{a,b}$, and for $1 \leq j \leq l-1$, set $x_j \defeq \tau_j \dotsm \tau_l$ if $l-j \leq m-2$; if $l-j > m-2$, then we must have $l=m$ and $j=1$, in this case set $x_1 = \tau_1 \dotsm \tau_{l-1}$. Now define $\bm{\mu}(z_i) = (x_1,x_2) \dotsm (x_{l-2},x_{l-1})$. Finally, set 
\begin{align*}
  \bm{\mu}(z) 
  = \ &\bm{\mu}(z_1) (\sigma_1 \sigma_1, \sigma_1 \sigma_1)^{e_1} (\sigma_1 \sigma_1, (\sigma_1 z_2)_{[1,m-1]}) ((\sigma_1 z_2)_{[1,m-1]}, z_2) \cdot \\
  \cdot &\bm{\mu}(z_2) (\sigma_2 \sigma_2, \sigma_2 \sigma_2)^{e_2} (\sigma_2 \sigma_2, (\sigma_2 z_3)_{[1,m-1]}) ((\sigma_2 z_3)_{[1,m-1]}, z_3) \dotsm.
\end{align*}
This defines a map $\bm{\mu}: \: \menge{z \in \gekl{a,b}^{\infty}}{z \text{ is } u,v \text{-free}} \to \cE^{\infty}$ which is the desired inverse. The rest of the argument proceeds exactly as in the proof of case~I of Theorem~\ref{Thm:OreOneRel:graph}.
\setlength{\parindent}{0cm} \setlength{\parskip}{0.5cm}

Next, we treat the torus knot case of type $(p,q)$. As in case~I of Theorem~\ref{Thm:OreOneRel:graph}, we have a bijection 
$$
  \menge{z \in \gekl{a,b}^{\infty}}{z \text{ is } u,v \text{-free}}, \, z \ma \chi_z.
$$
Moreover, the inverse of $\cE^{\infty} \to \menge{z \in \gekl{a,b}^{\infty}}{z \text{ is } u,v \text{-free}}, \, \mu \ma \bm{\sigma}(\mu)$ is given as follows: Every $u,v$-free infinite word $z \in \gekl{a,b}^{\infty}$ is of the form $z = z_1 z_2 \dotsm$ where $z_i$ are powers of $\sigma_i \in \gekl{a,b}$, i.e., $z_i \in \gekl{\ve, a^1, \dotsc, a^{p-1}, b^1, \dotsc, b^{q-1}}$, and we must have $\sigma_i \neq \sigma_{i+1}$. For each $z_i$, write $z_i \sigma_{i+1} = \tau_1 \dotsm \tau_l$ with $\tau_k \in \gekl{a,b}$. If $z_i \sigma_{i+1} \neq a^{p-1}b$, let $x_j \defeq \tau_j \dotsm \tau_l$ for all $1 \leq j \leq l-1$. If $z_i \sigma_{i+1} = a^{p-1}b$, let $x_1 = a^{p-1}$ and $x_j \defeq \tau_j \dotsm \tau_l$ for all $2 \leq j \leq l-1$. Now define $\mu_i = (x_1,x_2) \dotsm (x_{l-2},x_{l-1})$. Finally, set $\bm{\mu}(w) = \mu_1 \, (\sigma_1 \sigma_2, r(\mu_2)) \, \mu_2 \, (\sigma_2 \sigma_3, r(\mu_3)) \dotsm$. This defines a map $\bm{\mu}: \: \menge{z \in \gekl{a,b}^{\infty}}{z \text{ is } u,v \text{-free}} \to \cE^{\infty}$ which is the desired inverse. The rest of the argument proceeds exactly as in the proof of case~I of Theorem~\ref{Thm:OreOneRel:graph}.
\eproof
\setlength{\parindent}{0cm} \setlength{\parskip}{0cm}

Now let $\gamma: \: G \curvearrowright A$ be an action of $G$ on a C*-algebra $A$, and we write $\alpha_i \defeq (\gamma_a)_*^{-1}$ and $\beta_i \defeq (\gamma_b)_*^{-1}$ for the automorphisms $K_i(A) \cong K_i(A)$ induced by $\gamma_a^{-1}$ and $\gamma_b^{-1}$ in K-theory.

\btheo
\label{THM:K_examples}
For the concrete graph models from Lemma~\ref{LEM:K_ex--GraphModels}, $\phi$ in \eqref{CD:phi} is surjective and \eqref{exseq-K1} simplifies to
\begin{equation}
\label{exseq-K1_ex}
\begin{tikzcd}
\dotso \ar{r} & K_i(A) \ar["\ti{j}_i"]{r} & K_i(A) \oplus K_i(A) \ar["\ti{\pi}_i"]{r} & K_i(I) \ar["\ti{\partial}_i"]{r} & \dotso
\end{tikzcd}
\end{equation}
where
\begin{equation}
\label{j_I_2even}
  \ti{j}_i = 
  \begin{pmatrix}
  \id - \underbrace{\beta_i \alpha_i \dotsm \beta_i \alpha_i}_{m} \\
  \id + \beta_i \alpha_i + \dotso + (\beta_i \alpha_i)^{\frac{m-2}{2}}
  - \beta_i \rukl{\id + \alpha_i \beta_i + \dotso + (\alpha_i \beta_i)^{\frac{m-2}{2}}}
  \end{pmatrix}
\end{equation}
in the dihedral Artin-Tits case of type $I_2(m)$ with $m$ even,
\begin{equation}
\label{j_I_2odd}
  \ti{j}_i = 
  \begin{pmatrix}
  \id + \underbrace{\beta_i \alpha_i \dotsm \alpha_i \beta_i}_{m} \\
  \id + \beta_i \alpha_i + \dotso + (\beta_i \alpha_i)^{\frac{m-1}{2}}
  - \beta_i \rukl{\id + \alpha_i \beta_i + \dotso + (\alpha_i \beta_i)^{\frac{m-3}{2}}}
  \end{pmatrix}
\end{equation}
in the dihedral Artin-Tits case of type $I_2(m)$ with $m$ odd, and
\begin{equation}
\label{j_TorusKnot}
  \ti{j}_i = 
  \begin{pmatrix}
  \id + \alpha_i + \dotso + \alpha_i^{p-1} \\
  \id + \beta_i + \dotso + \beta_i^{q-1}
  \end{pmatrix}
\end{equation}
the torus knot case of type $(p,q)$.
\setlength{\parindent}{0cm} \setlength{\parskip}{0.5cm}

In addition, we have the exact sequence
\begin{equation*}
\begin{tikzcd}
\dotso \ar{r} & K_i(I) \ar["\iota_i"]{r} & K_i(A) \ar{r} & K_i(A \rtimes_r G) \ar{r} & \dotso
\end{tikzcd}
\quad {\rm from} \ \eqref{exseq-K2}.
\end{equation*}
\etheo
\setlength{\parindent}{0cm} \setlength{\parskip}{0cm}

Note the asymmetry between $a$ and $b$ in \eqref{E0_TorusKnot}, which is necessary to make sure that $\phi$ in \eqref{CD:phi} is surjective. In the following, we summarize the main steps in the computations leading to Theorem~\ref{THM:K_examples}. We point out that everything is constructive, so that it is possible to keep track of all the identifications made along the way.
\bproof
In the following, we write $\alpha = \alpha_i$ and $\beta = \beta_i$ to simplify notation. We have \eqref{exseq-K2} because in our cases, $P$ is right reversible (not only left reversible), so that $P \subseteq G$ is Toeplitz by \cite[\S~5.8]{CELY}. 
\setlength{\parindent}{0cm} \setlength{\parskip}{0.5cm}

We first treat the dihedral Artin-Tits case. With respect to the ordering of $\cE^0$ as in \eqref{E0_I_2even} and \eqref{E0_I_2odd}, the original homomorphism $j_i$ in \eqref{CD:Kcompare} is given by
\begin{equation}
\label{e:j_i}
  j_i = \id - 
  \rukl{
\resizebox{.4\hsize}{!}{
$
  \begin{array}{ccccccc|ccccccc}
  \alpha & & & & & & & \beta & & & & & & 0\\
  \vdots & & & & & & & & \ddots & & & & & \vdots \\ 
   & & & & & & & & & \beta & & & & \\ \hline
  \alpha & & & & & & & & \dots & 0 & 0 & \dots & \dots & 0 \\ \hline
   & & & & & & & & & & \beta & & & \\
  \vdots & & & & & & & & & & & \ddots & & \vdots\\
  \alpha & & & & & & & & & & & & \beta & 0 \\ \hline
  0 & \alpha & & & & & & & & & & & & \beta \\
  \vdots & & \ddots & & & & & & & & & & & \vdots \\ 
   & & & \alpha & & & & & & & & & & \\ \hline
  0 & \dots & \dots & 0 & 0 & \dots & & & & & & & & \beta \\ \hline
   & & & & \alpha & & & & & & & & & \\
  \vdots & & & & & \ddots & & & & & & & & \vdots \\
  0 & & & & & & \alpha & & & & & & & \beta
  \end{array}
$
}
  }
\end{equation}
where the row
$
  \begin{array}{cc|cccccc}
  \hline
  \alpha & \dots & \dots & 0 & 0 & \dots & \dots & 0 \\ \hline
  \end{array}
$
is at the position corresponding to the finite word $\underbrace{ab \dotsm ba}_{m-1}$ in \eqref{E0_I_2even} if $m$ is even and to $\underbrace{ab \dotsm ab}_{m-1}$ in \eqref{E0_I_2odd} if $m$ is odd, and the row
$
  \begin{array}{cccccc|cc}
  \hline
  0 & \dots & \dots & 0 & 0 & \dots & \dots & \beta \\ \hline
  \end{array}
$
is at the position corresponding to the finite word $\underbrace{ba \dotsm bab}_{m-1}$ in \eqref{E0_I_2even} if $m$ is even and to $\underbrace{ba \dotsm ba}_{m-1}$ in \eqref{E0_I_2odd} if $m$ is odd. By performing elementary row operations, which correspond to post-composition with isomorphisms, $j_i$ is transformed to 
\begin{equation}
\label{e:j'_i}
  j'_i = 
  \rukl{
\resizebox{.4\hsize}{!}{
$
  \begin{array}{c|ccccccccccc|c}
  \boldsymbol{w} & & & & & & & & & & & & \boldsymbol{x} \\ \hline
  ? & 1 & & & & & & & & & & & ? \\ 
  ? & & \ddots & & & & & & & & & & ?\\
  ? & & & 1 & & & & & & & & & ? \\ \hline
  - \alpha & 0 & \dots & 0 & 1 & 0 & \dots & & & & & & 0 \\ \hline
  ? & & & & & 1 & & & & & & & ? \\
  ? & & & & & & \ddots & & & & & & ? \\
  ? & & & & & & & 1 & & & & & ? \\ \hline
  0 & & & & & & \dots & 0 & 1 & 0 & \dots & 0 & - \beta \\ \hline
  ? & & & & & & & & & 1 & & & ? \\
  ? & & & & & & & & & & \ddots & & ? \\
  ? & & & & & & & & & & & 1 & ? \\ \hline
  \boldsymbol{y} & & & & & & & & & & & & \boldsymbol{z}
  \end{array}
$
}
  }
\end{equation}
where we write $1$ for $\id$ and 
\begin{eqnarray*}
  \boldsymbol{w} &=& 1 - \alpha - \beta \alpha^2 - \dotso - \underbrace{\beta \alpha \dotsm \beta \alpha^2}_{m-1}, \ \boldsymbol{x} = - \beta^2 - \beta \alpha \beta^2 - \dotso - \underbrace{\beta \alpha \dotsm \beta \alpha \beta^2}_{m-2},\\
  \boldsymbol{y} &=& - \alpha^2 - \alpha \beta \alpha^2 - \dotso - \underbrace{\alpha \beta \dotsm \alpha \beta \alpha^2}_{m-2}, \
  \boldsymbol{z} = 1 - \beta - \alpha \beta^2 - \dotso - \underbrace{\alpha \beta \dotsm \alpha \beta^2}_{m-1}
  \quad {\rm if} \ m \ \text{is even, and}\\
  \boldsymbol{w} &=& 1 - \alpha - \beta \alpha^2 - \dotso - \underbrace{\beta \alpha \dotsm \beta \alpha^2}_{m-2}, \
  \boldsymbol{x} = - \beta^2 - \beta \alpha \beta^2 - \dotso - \underbrace{\beta \alpha \dotsm \beta \alpha \beta^2}_{m-1},\\
  \boldsymbol{y} &=& - \alpha^2 - \alpha \beta \alpha^2 - \dotso - \underbrace{\alpha \beta \dotsm \alpha \beta \alpha^2}_{m-1}, \
  \boldsymbol{z} = 1 - \beta - \alpha \beta^2 - \dotso - \underbrace{\alpha \beta \dotsm \alpha \beta^2}_{m-2}
  \quad {\rm if} \ m \ \text{is odd}.
\end{eqnarray*}
The map $\phi_i$ in \eqref{CD:Kcompare} is given by $\phi_i = (\boldsymbol{\chi}, 0, \dotsc, 0, \bar{\omega}, 0, \dotsc, 0, \omega, 0, \dotsc, 0, \boldsymbol{\psi})$, where 
\begin{eqnarray*}
  \boldsymbol{\chi} &=& \alpha^2 + \beta \alpha^2 + \alpha \beta \alpha^2 + \dotso + \underbrace{\dotsm \beta \alpha \beta \alpha^2}_{m-1}, \
  \boldsymbol{\psi} = \beta^2 + \alpha \beta^2 + \beta \alpha \beta ^2 + \dotso + \underbrace{\dotsm \alpha \beta \alpha \beta ^2}_{m-1},\\
  \bar{\omega} &=& \underbrace{\dotsm \beta \alpha \beta \alpha}_{m-1} \ \text{is at the position corresponding to the finite word} \ \underbrace{ab \dotsm ba}_{m-1} \ {\rm in} \ \eqref{E0_I_2even} \ {\rm if} \ m \ \text{is even} \\ 
  && \text{and to} \ \underbrace{ab \dotsm ab}_{m-1} \ {\rm in} \ \eqref{E0_I_2odd} \ {\rm if} \ m \ \text{is odd, and} \\
  \omega &=& \underbrace{\dotsm \alpha \beta \alpha \beta}_{m-1} \ \text{is at the position corresponding to the finite word} \ \underbrace{ba \dotsm bab}_{m-1}\ {\rm in} \ \eqref{E0_I_2even}\  {\rm if} \ m \ \text{is even} \\
  && \text{and to} \ \underbrace{ba \dotsm ba}_{m-1} \ {\rm in} \ \eqref{E0_I_2odd} \ {\rm if} \ m \ \text{is odd}.
\end{eqnarray*}
Now the following map $\bigoplus_{\cE^0 \setminus \gekl{v}} K_i(A) \to \bigoplus_{\cE^0} K_i(A)$ induces an isomorphism $\bigoplus_{\cE^0 \setminus \gekl{v}} K_i(A) \cong \ker(\phi_i)$, where $v$ is the finite word $\underbrace{ba \dotsm bab}_{m-1}$ in \eqref{E0_I_2even} if $m$ is even and to $\underbrace{ba \dotsm ba}_{m-1}$ in \eqref{E0_I_2odd} if $m$ is odd:
\begin{equation*}
  \rukl{
\resizebox{.5\hsize}{!}{
$
  \begin{array}{c|ccc|c|ccc||ccc|c}
  1 & & & & & & & & & & & \\ \hline
  & 1 & & & & & & & & & & \\
  & & \ddots& & & & & & & & & \\
  & & & 1 & & & & & & & & \\ \hline
  & & & & 1 & & & & & & & \\ \hline
  & & & & & 1 & & & & & & \\
  & & & & & & \ddots & & & & & \\
  & & & & & & & 1 & & & & \\ \hline
  - \omega^{-1} \boldsymbol{\chi} & 0 & \dots & 0 & - \omega^{-1} \bar{\omega} & 0 & \dots & 0 & 0 & \dots & 0 & - \omega^{-1} \boldsymbol{\psi} \\ \hline
  & & & & & & & & 1 & & & \\
  & & & & & & & & & \ddots & & \\
  & & & & & & & & & & 1 & \\ \hline
  & & & & & & & & & & & 1
  \end{array}
$
}
  }
\end{equation*}
where the double vertical line indicates the position of $v$. Pre-composing $j'_i$ by this isomorphism yields the homomorphism $\ti{j'_i}: \: \bigoplus_{\cE^0 \setminus \gekl{v}} K_i(A) \cong \ker(\phi_i) \to \bigoplus_{\cE^0} K_i(A)$ given by

\begin{equation*}
  \ti{j'_i} = 
  \rukl{
\resizebox{.5\hsize}{!}{
$
  \begin{array}{c|ccc|c|ccc||ccc|c}
  \boldsymbol{w} & & & & & & & & & & & \boldsymbol{x} \\ \hline
  ? & 1 & & & & & & & & & & ? \\ 
  ? & & \ddots & & & & & & & & & ?\\
  ? & & & 1 & & & & & & & & ? \\ \hline
  - \alpha & 0 & \dots & 0 & 1 & 0 & \dots & & & & & 0 \\ \hline
  ? & & & & & 1 & & & & & & ? \\
  ? & & & & & & \ddots & & & & & ? \\
  ? & & & & & & & 1 & & & & ? \\ \hline
  - \omega^{-1} \boldsymbol{\chi} & & & & - \omega^{-1} \bar{\omega} & & \dots & 0 & 0 & \dots & 0 & - \beta - \omega^{-1} \boldsymbol{\psi} \\ \hline
  ? & & & & & & & & 1 & & & ? \\
  ? & & & & & & & & & \ddots & & ? \\
  ? & & & & & & & & & & 1 & ? \\ \hline
  \boldsymbol{y} & & & & & & & & & & & \boldsymbol{z}
  \end{array}
$
}
  }
\end{equation*}
Since we have many columns with only one non-zero entry given by $1$, we can first perform elementary column operations and then permutations so that $\ti{j'_i}$ is transformed to a matrix whose upper left corner consists of the identity matrix and the lower right corner is of the form
\begin{equation*}
  \rukl{
  \begin{array}{ccc}
  \boldsymbol{w} & 0 & \boldsymbol{x} \\ 
  -\alpha & 1 & 0 \\
  - \omega^{-1} \boldsymbol{\chi} & - \omega^{-1} \bar{\omega} & - \beta - \omega^{-1} \boldsymbol{\psi} \\
  \boldsymbol{y} & 0 & \boldsymbol{z}
  \end{array}
  }
\end{equation*}
Let us keep the upper left corner and further transform the lower right corner. Multiplying the third row from the left by $\omega$ and then adding the second column, multiplied by $\alpha$ from the right, to the first column, we obtain in the lower right corner
\begin{equation*}
  \rukl{
  \begin{array}{ccc}
  \boldsymbol{w} & 0 & \boldsymbol{x} \\ 
  0 & 1 & 0 \\
  - \boldsymbol{\chi} - \bar{\omega} \alpha & - \bar{\omega} & - \omega \beta - \boldsymbol{\psi} \\
  \boldsymbol{y} & 0 & \boldsymbol{z}
  \end{array}
  }
\end{equation*}
Eliminating the entry $- \bar{\omega}$ and permuting, we obtain
\begin{equation*}
  \rukl{
  \begin{array}{c|cc}
  1 & & \\ \hline
  &\boldsymbol{w} & \boldsymbol{x} \\ 
  & - \boldsymbol{\chi} - \bar{\omega} \alpha & - \omega \beta - \boldsymbol{\psi} \\
  & \boldsymbol{y} & \boldsymbol{z}
  \end{array}
  }
\end{equation*}
We further transform the lower right $3 \times 2$ matrix. Multiplying its first row by $\alpha$ from the left, subtracting the second row and adding the third row, multiplied by $\beta$ from the left, the lower right $3 \times 2$ matrix becomes
\begin{equation*}
  \rukl{
  \begin{array}{cc}
  \alpha & \beta \\ 
  - \boldsymbol{\chi} - \bar{\omega} \alpha & - \omega \beta - \boldsymbol{\psi} \\
  \boldsymbol{y} & \boldsymbol{z}
  \end{array}
  },
\end{equation*}
where we used the equations $\alpha \boldsymbol{x} + (\omega \beta + \boldsymbol{\psi}) + \beta \boldsymbol{z} = \beta$ and $\alpha \boldsymbol{w} + (\boldsymbol{\chi} + \bar{\omega} \alpha) + \beta \boldsymbol{y} = \alpha$. Subtracting from the second column the first column, multiplied from the right by $\alpha^{-1} \beta$, we obtain
\begin{equation*}
  \rukl{
  \begin{array}{cc}
  \alpha & 0 \\ 
  - \boldsymbol{\chi} - \bar{\omega} \alpha & - \omega \beta - \boldsymbol{\psi} - (- \boldsymbol{\chi} - \bar{\omega} \alpha) \alpha^{-1} \beta \\
  \boldsymbol{y} & \boldsymbol{z} - \boldsymbol{y} \alpha^{-1} \beta
  \end{array}
  }
\end{equation*}
Now eliminate the entries in the first column below $\alpha$, multiply the first column by $\alpha^{-1}$, conjugate the second column by $\beta$ and multiply the second row by $-1$ to obtain
\begin{equation*}
  \rukl{
  \begin{array}{cc}
  1 & 0 \\ 
  0 & \beta \omega + \beta \boldsymbol{\psi} \beta^{-1} + \beta (- \boldsymbol{\chi} - \bar{\omega} \alpha) \alpha^{-1} \\
  0 & \beta \boldsymbol{z} \beta^{-1} - \beta \boldsymbol{y} \alpha^{-1}
  \end{array}
  }
\end{equation*}
Adding the third row, multiplied by $1 + \beta$ from the left, to the second row yields 
\begin{equation*}
  \rukl{
  \begin{array}{c|c}
  1 & \\ \hline 
  & \ti{j}_i
  \end{array}
  }
\end{equation*}
with $\ti{j}_i$ as in \eqref{j_I_2even} or \eqref{j_I_2odd}.

All in all, we have transformed $\ti{j'_i}$ by performing elementary column and row operations (which correspond to pre- and post-composition with isomorphisms) to $\ti{j''_i}: \: \bigoplus_{\cE^0 \setminus \gekl{v}} K_i(A) \cong \ker(\phi_i) \to \bigoplus_{\cE^0} K_i(A)$ given by
\begin{equation*}
  \ti{j''_i} = 
  \rukl{
  \begin{array}{ccc|c}
  1 & & & \\ 
   & \ddots & & \\
   & & 1 & \\ \hline
   & & & \ti{j}_i
  \end{array}
  }
\end{equation*}
with $\ti{j}_i$ as in \eqref{j_I_2even} or \eqref{j_I_2odd}. This completes the proof in the dihedral Artin-Tits case.

Let us now treat the torus knot case. With respect to the ordering of $\cE^0$ as in \eqref{E0_TorusKnot}, the original homomorphism $j_i$ in \eqref{CD:Kcompare} is given by
\begin{equation}
\label{e:j_i_Torus}
  j_i = \id - 
  \rukl{
  \begin{array}{c|ccc|c|ccc}
  0 & \alpha & & 0 & \beta & & & \\
  \vdots & & \ddots & & \vdots & & 0 & \\
  0 & 0 & & \alpha & & & & \\ \hline
  0 & 0 & \dots & 0 & \beta & 0 & \dots & 0 \\ \hline
  \alpha & & & & 0 & \beta & & 0 \\
  \vdots & & 0 & & \vdots & & \ddots & \\
  & & & & 0 & 0 & & \beta \\
  \alpha & & & & 0 & 0 & & 0  
  \end{array}
  }
\end{equation}
where the row
$$
  \begin{array}{c|ccc|c|ccc} \hline
  0 & 0 & \dots & 0 & \beta & 0 & \dots & 0 \\ \hline
  \end{array}
$$
is at the position corresponding to the finite word $a^{p-1}$ in \eqref{E0_TorusKnot}, and the column
$$
  \rukl{
  \begin{array}{ccc|c|cccc}
  \hline
  \beta & \dots & & \beta & 0 & \dots & 0 & 0
  \\ \hline
  \end{array}
  }^t
$$
is at the position corresponding to the finite word $ba$ in \eqref{E0_TorusKnot}. By performing elementary row operations, which correspond to post-composition with isomorphisms, $j_i$ is transformed to 
\begin{equation}
\label{e:j'_i_Torus}
  j'_i = 
  \rukl{
  \begin{array}{ccc|c|ccc}
  1 & & 0 & \boldsymbol{w} & & & \\
   & \ddots & & ? & & 0 & \\ \hline
   0 & & 1 & - \beta & & & \\ \hline
  0 & \dots & 0 & \boldsymbol{x} & 0 & \dots & 0 \\
  & & & ? & 1 & & 0 \\
   & 0 & & \vdots & & \ddots & \\
  & & & ? & 0 & & 1  
  \end{array}
  }
\end{equation}
where 
\begin{eqnarray*}
  \boldsymbol{w} &=& - \beta - \alpha \beta - \alpha^2 \beta - \dotso - \alpha^{p-2} \beta, \ 
  \boldsymbol{x} = 1 - (\alpha + \beta \alpha + \beta^2 \alpha + \dotso + \beta^{q-2} \alpha)(\beta + \alpha \beta + \alpha^2 \beta + \dotso + \alpha^{p-2} \beta).
\end{eqnarray*}

The map $\phi_i$ in \eqref{CD:Kcompare} is given by $\phi = (\boldsymbol{y}, 0, \dotsc, 0, \omega, \boldsymbol{z}, 0, \dotsc, 0)$, where 
\begin{eqnarray*}
  \boldsymbol{y} &=& \beta \alpha + \beta^2 \alpha + \dotso + \beta^{q-1} \alpha, \ 
  \omega = \alpha^{p-1}, \
  \boldsymbol{z} = \alpha \beta + \alpha^2 \beta + \dotso + \alpha^{p-2} \beta,
\end{eqnarray*}
where $\omega$ is at the position corresponding to the finite word $a^{p-1}$ in \eqref{E0_TorusKnot} and $\boldsymbol{z}$ is at the position corresponding to the finite word $ba$ in \eqref{E0_TorusKnot}.

Now the following map $\bigoplus_{\cE^0 \setminus \gekl{a^{p-1}}} K_i(A) \to \bigoplus_{\cE^0} K_i(A)$ induces an isomorphism $\bigoplus_{\cE^0 \setminus \gekl{a^{p-1}}} K_i(A) \cong \ker(\phi_i)$:
\begin{equation*}
  \rukl{
  \begin{array}{cccc||c|ccc}
  1 & 0 & \dots & 0 & 0 & 0 & \dots & 0 \\
  0 & 1 & & 0 &  & & & \\
  \vdots & & \ddots & & \vdots & & 0 & \\
  0 & 0 & & 1 & 0 & & & \\ \hline
  - \omega^{-1} \boldsymbol{y} & 0 & \dots & 0 & - \omega^{-1} \boldsymbol{z} & 0 & \dots & 0 \\ \hline
  & & & & 1 & 0 & \dots & 0 \\
  & & & & 0 & 1 & & 0 \\
  & 0 & & & \vdots &  & \ddots & \\
  & & & & 0 & 0 & & 1
  \end{array}
  }
\end{equation*}
where the double vertical line indicates the position of $a^{p-1}$. Pre-composing $j'_i$ by this isomorphism yields the homomorphism $\ti{j'_i}: \: \bigoplus_{\cE^0 \setminus \gekl{a^{p-1}}} K_i(A) \cong \ker(\phi_i) \to \bigoplus_{\cE^0} K_i(A)$ given by
\begin{equation*}
  \ti{j'_i} = 
  \rukl{
  \begin{array}{cccc||c|ccc}
  1 & 0 & \dots & 0 & \boldsymbol{w} & 0 & \dots & 0 \\
  0 & 1 & & 0 & ? & & & \\
  \vdots & & \ddots & & \vdots & & 0 & \\
  0 & 0 & & 1 & ? & & & \\ \hline
  - \omega^{-1} \boldsymbol{y} & 0 & \dots & 0 & - \beta - \omega^{-1} \boldsymbol{z} & 0 & \dots & 0 \\ \hline
  & & & & \boldsymbol{x} & 0 & \dots & 0 \\
  & & & & ? & 1 & & 0 \\
  & 0 & & & \vdots &  & \ddots & \\
  & & & & ? & 0 & & 1
  \end{array}
  }
\end{equation*}
Finally, by performing elementary column and row operations corresponding to pre- and post-composition with isomorphisms, $\ti{j'_i}$ is transformed to $\ti{j''_i}: \: \bigoplus_{\cE^0 \setminus \gekl{a^{p-1}}} K_i(A) \cong \ker(\phi_i) \to \bigoplus_{\cE^0} K_i(A)$ given by
\begin{equation*}
  \ti{j''_i} = 
  \rukl{
  \begin{array}{ccc|c}
  1 & & & \\ 
   & \ddots & & \\
   & & 1 & \\ \hline
   & & & \ti{j}_i
  \end{array}
  }
\end{equation*}
with $\ti{j}_i$ as in \eqref{j_TorusKnot}. This completes the proof in the torus knot case.
\eproof

\bremark
\label{R:j--j}
The proof shows that $\ti{j}_i$, $\ti{\pi}_i$ and $\ti{\partial}_i$ are related to the maps $j_i$, $\pi_i$ and $\partial_i \circ \psi_i^{-1} \circ q_i$ from \eqref{exseq-K1} in the following way: We have isomorphisms $\ker(\phi_i) \cong K_i \oplus K_i(A)$ and $\bigoplus_{\cE^0} K_i(A) \cong K_i \oplus (K_i(A) \oplus K_i(A))$ for some abelian group $K_i$, which fit into a commutative diagram
$$
\begin{tikzcd}
\ker(\phi_i) \ar["j_i"]{r} \ar["\cong"']{d} & \bigoplus_{\cE^0} K_i(A) \ar["\cong"]{d} \\
K_i \oplus K_i(A) \ar[r] & K_i \oplus (K_i(A) \oplus K_i(A)) 
\end{tikzcd}
$$
where the lower horizontal arrow is given by $\rukl{\begin{smallmatrix} * & 0 \\ 0 & \ti{j}_i \end{smallmatrix}}$ for some automorphism $*: \: K_i \cong K_i$. Then $\ti{\pi}_i$ is given by the composition
$$
  K_i(A) \oplus K_i(A) \into K_i \oplus (K_i(A) \oplus K_i(A)) \cong \bigoplus_{\cE^0} K_i(A) \overset{\pi_i}{\lori} K_i(I),
$$
where the first map is the canonical inclusion, and $\ti{\partial}_i$ is given by the composition
$$
  K_i(I) \overset{\partial_i \circ \psi_i^{-1} \circ q_i}{\lori} \ker(\phi_i) \cong K_i \oplus K_i(A) \onto K_i(A),
$$
where the last map is the canonical projection. In particular, we have $\img(\ti{\pi}_i) = \img(\pi_i)$.
\eremark
\setlength{\parindent}{0cm} \setlength{\parskip}{0.5cm}

\subsection{Examples}
\label{ss:Ex-K_*}

We will use the following fact several times.
\bremark
\label{abelianization}
The commutator (or derived) subgroup $G'$ of a group $G$ is the normal subgroup generated by all elements of the form $ghg^{-1}h^{-1}$ for $g,h\in G$.
The quotient $G_{\textup{ab}}=G/G'$ is called the abelianization of $G$ and we let $f_{\textup{ab}}\colon G\to G_{\textup{ab}}$ denote the canonical quotient map.
The following result can now be deduced from \cite{MO} (see Theorem~1.2 and the following paragraph in that paper):
\newline
Suppose that $G$ is a torsion-free group satisfying the Baum-Connes conjecture.
Then there is a well-defined split-injective group homomorphism $G_{\textup{ab}} \to K_1(C^*_r(G))$ defined by $f_{\textup{ab}}(g)\mapsto [\lambda(g)]_1$.
\eremark

\begin{example}
Let $A=\Cz$, so of course, the action is trivial, $K_0(\Cz)=\Zz$, and $K_1(\Cz)=\{0\}$, and \eqref{exseq-K1_ex} gives the exact sequence
\[
0 \longrightarrow K_1(I) \overset{\tilde{\partial}_1}{\longrightarrow} \Zz \overset{\tilde{j}_0}{\longrightarrow} \Zz\oplus\Zz \overset{\tilde{\pi}_0}{\longrightarrow} K_0(I) \longrightarrow 0.
\]
Consider the case $I_2(m)$ with $m\geq 2$ even.
Then \eqref{j_I_2even} gives that $\tilde{j}_0$ is the zero map, and consequently $K_1(I)\cong\Zz$ and $K_0(I)\cong\Zz^2$.
Moreover, \eqref{exseq-K2} now gives the exact sequence
\[
0 \longrightarrow K_1(C^*_r(I_2(m))) \longrightarrow K_0(I) \overset{\iota_0}{\longrightarrow} K_0(A) \longrightarrow K_0(C^*_{\lambda}(I_2(m))) \longrightarrow K_1(I) \longrightarrow 0.
\]
Since $\iota_0\circ\pi_0$ is the zero map (using the formula \eqref{iota-pi-formula}), and $\pi_0$ is surjective, we get that $\iota_0$ is the zero map.
Hence, $K_0(C^*_{\lambda}(I_2(m)))\cong K_0(A) \oplus K_1(I) \cong \Zz^2$ and $K_1(C^*_{\lambda}(I_2(m)))\cong K_0(I) \cong \Zz^2$. It follows from Remark~\ref{abelianization} that the $K_1$-group is generated by $[\lambda(a)]$ and $[\lambda(b)]$. Furthermore, one copy of $\Zz$ in $K_0$ is generated by $[1]$.
\setlength{\parindent}{0.5cm} \setlength{\parskip}{0cm}

Moreover, having computed the maps $j_i$ and $j'_i$ in \eqref{e:j_i} and \eqref{e:j'_i}, it is easy to determine the K-theory of $Q \defeq C^*_r(G \ltimes \ti{\Omega}_{\infty})$: We get $K_0(Q) \cong \Zz \oplus (\Zz / \tfrac{m-2}{2} \Zz)$ and $K_1(Q) \cong \Zz$. In case $m \geq 4$, Remark~\ref{R:pisC} tells us that this determines the stable UCT Kirchberg algebra $Q$ up to isomorphism. Furthermore, we have exact sequences $0 \to J \to C^*_{\lambda}(P) \to C^*_{\lambda}(I_2(m)) \to 0$ and $0 \to \cK(\ell^2 P) \to J \to Q \to 0$.
\setlength{\parindent}{0cm} \setlength{\parskip}{0.5cm}

Assume next that $m\geq 3$ is odd.
Then \eqref{j_I_2odd} gives that $\tilde{j}_0$ is the map $\left(\begin{smallmatrix} 2 \\ 1 \end{smallmatrix}\right)$, which is injective with cokernel isomorphic to $\Zz$.
Therefore, $K_0(I)\cong\Zz$ and $\tilde{\partial}_1$ is the zero map, so $K_1(I)=\{0\}$.
From \eqref{exseq-K2} we get the exact sequence
\[
0 \longrightarrow K_1(C^*_{\lambda}(I_2(m))) \longrightarrow K_0(I) \overset{\iota_0}{\longrightarrow} K_0(A) \longrightarrow K_0(C^*_{\lambda}(I_2(m))) \longrightarrow 0.
\]
Since $\iota_0\circ\pi_0$ is the zero map, and $\pi_0$ is surjective, then $\iota_0$ is the zero map.
Hence, $K_0(C^*_{\lambda}(I_2(m)))\cong K_0(A)\cong \Zz$ and thus $K_1(C^*_{\lambda}(I_2(m)))\cong K_0(I) \cong \Zz$. The $K_0$-group is generated by the identity, and by Remark~\ref{abelianization}, the $K_1$-group is generated by $[\lambda(a)]=[\lambda(b)]$.
\setlength{\parindent}{0.5cm} \setlength{\parskip}{0cm}

As above, using \eqref{e:j_i} and \eqref{e:j'_i}, it is easy to determine the K-theory of $Q \defeq C^*_r(G \ltimes \ti{\Omega}_{\infty})$: We get $K_0(Q) \cong \Zz / (m-2) \Zz$ and $K_1(Q) \cong \gekl{0}$. Remark~\ref{R:pisC} tells us that this determines the stable UCT Kirchberg algebra $Q$ up to isomorphism, and we get $Q \cong \cK \otimes \cO_{m-1}$ (it is possible to write down an explicit isomorphism). So we have exact sequences $0 \to J \to C^*_{\lambda}(P) \to C^*_{\lambda}(I_2(m)) \to 0$ and $0 \to \cK(\ell^2 P) \to J \to \cK \otimes \cO_{m-1} \to 0$.
\setlength{\parindent}{0cm} \setlength{\parskip}{0.5cm}

Finally, consider the torus knot group $T(p,q)$ of type $(p,q)$ for $p,q\geq 2$, and set $g=\gcd(p,q)$.
Then \eqref{j_TorusKnot} gives that $\tilde{j}_0$ is the map $\left(\begin{smallmatrix} p \\ q \end{smallmatrix}\right)$, which is injective with cokernel isomorphic to $\Zz\oplus(\Zz / g \Zz)$.
Therefore, $K_0(I)\cong\Zz\oplus(\Zz / g\Zz)$ and $\tilde{\partial}_1$ is the zero map, so $K_1(I)=\{0\}$.
From \eqref{exseq-K2} we get the exact sequence
\[
0 \longrightarrow K_1(C^*_{\lambda}(T(p,q))) \longrightarrow K_0(I) \overset{\iota_0}{\longrightarrow} K_0(A) \longrightarrow K_0(C^*_{\lambda}(T(p,q))) \longrightarrow 0.
\]
Since $\iota_0\circ\pi_0$ is the zero map, and $\pi_0$ is surjective, then $\iota_0$ is the zero map.
Hence, $K_0(C^*_{\lambda}(T(p,q)))\cong K_0(A)\cong \Zz$ and $K_1(C^*_{\lambda}(T(p,q)))\cong K_0(I)\cong\Zz\oplus(\Zz / g\Zz)$. The $K_0$-group is generated by the identity. To write down generators for $K_1$, let $x$ and $y$ be integers such that $\tfrac{p}{g} \cdot x - \tfrac{q}{g} \cdot y = 1$. Then $[\lambda(a^y b^{-x})]$ generates the direct summand $\Zz$ and $[\lambda(a^{p/g} b^{-q/g})]$ generates the direct summand $\Zz / g \Zz$ of $K_1$ (see Remark~\ref{abelianization}).
\setlength{\parindent}{0.5cm} \setlength{\parskip}{0cm}

Having computed the maps $j_i$ and $j'_i$ in \eqref{e:j_i_Torus} and \eqref{e:j'_i_Torus}, it is easy to determine the K-theory of $Q \defeq C^*_r(G \ltimes \ti{\Omega}_{\infty})$: We get $K_0(Q) \cong \Zz$, $K_1(Q) \cong \Zz$ if $p = q = 2$, and $K_0(Q) \cong \Zz / ((p-1)(q-1) - 1) \Zz$, $K_1(Q) \cong \gekl{0}$ otherwise. In case $(p,q) \neq (2,2)$, Remark~\ref{R:pisC} tells us that this determines the stable UCT Kirchberg algebra $Q$ up to isomorphism, and we get $Q \cong \cK \otimes \cO_{(p-1)(q-1)}$. Thus we obtain exact sequences $0 \to J \to C^*_{\lambda}(P) \to C^*_{\lambda}(T(p,q)) \to 0$ and $0 \to \cK(\ell^2 P) \to J \to \cK \otimes \cO_{(p-1)(q-1)} \to 0$.
\setlength{\parindent}{0cm} \setlength{\parskip}{0.5cm}

\end{example}

\begin{example}
The braid group $B_4$ is given by the presentation
\[
\langle a,b,c \, \vert \, aba=bab, bcb=cbc, ac=ca\rangle.
\]
Consider the surjective homomorphism $\varphi\colon B_4\to B_3$ given by $a\mapsto a$, $b\mapsto b$, $c\mapsto a$.
Then it follows from \cite[Theorem~2.1]{GL} that $\ker\varphi$ is isomorphic to the free group $\Fz_2$ with generators
\[
x_1=ca^{-1} \quad\text{and}\quad x_2=bca^{-1}b^{-1}.
\]
Moreover, the obvious inclusion $B_3\into B_4$ defines an action $\gamma$ of $B_3$ on $\Fz_2$ given by
\[
\begin{split}
\gamma_a(x_1)=ax_1a^{-1}=x_1, &\quad \gamma_b(x_1)=bx_1b^{-1}=x_2, \\
\gamma_a(x_2)=ax_2a^{-1}=x_1^{-1}x_2, &\quad \gamma_b(x_2)=bx_2b^{-1}=x_2x_1^{-1}x_2,
\end{split}
\]
and $B_4\cong\Fz_2\rtimes_\gamma B_3$.

Set $A=C^*_{\lambda}(\Fz_2)$, and note that $C^*_{\lambda}(B_4)\cong A\rtimes_r B_3$, where the action is induced from the one above.
It is shown in \cite[Corollary~3.2]{PV} that the $K$-theory of $C^*_{\lambda}(\Fz_2)$ is given by $K_0(C^*_{\lambda}(\Fz_2))=\Zz[1]$ and $K_1(C^*_{\lambda}(\Fz_2))=\Zz[\lambda(x_1)] \oplus \Zz[\lambda(x_2)]$,
where $\lambda$ is the left regular representation
Hence, it follows that $\alpha_0=\beta_0=\textup{id}$ on $K_0$, while on $K_1$ we have
\[
\alpha_1=\begin{pmatrix} 1 & -1 \\ 0 & 1 \end{pmatrix}^{-1}=\begin{pmatrix} 1 & 1 \\ 0 & 1 \end{pmatrix}
\quad\text{and}\quad
\beta_1=\begin{pmatrix} 0 & -1 \\ 1 & 2 \end{pmatrix}^{-1}=\begin{pmatrix} 2 & 1 \\ -1 & 0 \end{pmatrix}.
\]
Thus, using \eqref{j_I_2odd}, we compute that
\[
\tilde{j}_0=\begin{pmatrix} 2 \\ 1 \end{pmatrix}
\quad\text{and}\quad
\tilde{j}_1=\begin{pmatrix} 2 & 2 \\ -1 & 0 \\ 1 & 2 \\ 0 & 0 \end{pmatrix}.
\]
Both maps are injective, with cokernels isomorphic to $\Zz$ and $\Zz^2\oplus(\Zz / 2\Zz)$, respectively.
This means that the index maps are trivial, so \eqref{exseq-K1_ex} gives short exact sequences
\[
0 \longrightarrow K_i(A) \overset{\tilde{j}_i}{\longrightarrow} K_i(A) \oplus K_i(A) \overset{\tilde{\pi}_i}{\longrightarrow} K_i(I) \longrightarrow 0
\]
for $i=0,1$.
Thus, $K_0(I)=\Zz$ and $K_1(I)=\Zz^2\oplus(\Zz / 2\Zz)$.

Moreover, using the formula \eqref{iota-pi-formula}, we see that $\iota_i\circ\pi_i\colon K_i(A) \oplus K_i(A) \oplus K_i(A) \oplus K_i(A) \to K_i(A)$ is the zero map for $i=0$ and
\begin{equation*}
  \rukl{
  \begin{array}{cc|cc|cc|cc}
  0 & -1 & 1 & 1 & 0 & -1 & 3 & 1 \\
  1 & 3 & 0 & 0 & 0 & 1 & -1 & 0  
  \end{array}
  }
\end{equation*}
for $i=1$ (with respect to $\cE^0 = \gekl{aa,ab,ba,bb}$ with  this ordering). In both cases, $\pi_i$ is surjective (as $\img(\pi_i) = \img(\ti{\pi}_i)$ by Remark~\ref{R:j--j}), so $\iota_0=0$, and $\iota_1$ maps $\Zz^2\oplus(\Zz / 2\Zz)$ onto $\Zz^2$.
Therefore, \eqref{exseq-K2} gives an exact sequence
\[
0 \longrightarrow K_0(A) \longrightarrow K_0(A\rtimes_r B_3) \longrightarrow K_1(I) \overset{\iota_1}{\longrightarrow} K_1(A) \longrightarrow K_1(A\rtimes_r B_3) \longrightarrow K_0(I) \longrightarrow 0,
\]
that is,
\[
0 \longrightarrow \Zz \longrightarrow K_0(A\rtimes_r B_3) \longrightarrow \Zz^2\oplus(\Zz / 2\Zz) \overset{\iota_1}{\longrightarrow} \Zz^2 \longrightarrow K_1(A\rtimes_r B_3) \longrightarrow \Zz \longrightarrow 0.
\]
Since $\iota_1$ is surjective, $K_1(C^*_{\lambda}(B_4))\cong K_0(I) \cong \Zz$, while $K_0(C^*_{\lambda}(B_4))$ is either $\Zz$ or $\Zz\oplus(\Zz / 2\Zz)$.
It is known that $[1]$ must generate a copy of $\Zz$ as a direct summand of the $K_0$-group (see for instance \cite[Proof of Lemma~4.4]{Li14}),
so therefore it follows that $K_0(C^*_{\lambda}(B_4))\cong\Zz\oplus(\Zz / 2\Zz)$. Moreover, it follows from Remark~\ref{abelianization} that $K_1$ is generated by $[\lambda(a)]=[\lambda(b)]=[\lambda(c)]$. 
\end{example}

\begin{example}
Artin's representation of braid groups (see \cite{Om} for a discussion) is defined for $n\geq 3$ as the action $\gamma$ of $B_n$ on $\Fz_n$, with canonical generators $\sigma_i$ and $x_i$, respectively, given by
\[
\gamma(\sigma_i)(x_j)=
\begin{cases}
x_{i+1} & \text{if } j=i, \\
x_{i+1}^{-1}x_i^{\vphantom{-1}}x_{i+1}^{\vphantom{-1}} & \text{if } j=i+1, \\
x_j & \text{else}.
\end{cases}
\]
We consider the case $n=3$, and want to compute K-theory for the reduced group C*-algebra of $G = \Fz_3 \rtimes B_3$. Note that by \cite[Proposition~2.1]{CP2}, $G$ is the irreducible Artin-Tits group given by the presentation 
$$
  \spkl{a, b, c \, \vert \, aba=bab, bcbc=cbcb, ac=ca}.
$$ 
Let us compute the induced action (denoted by the same symbol) $\gamma$ of $B_3$, now generated by $a,b$, on $A = C^*_{\lambda}(\Fz_3)$ in K-theory. We get $\alpha_0=\beta_0=\textup{id}$ on $K_0$, while on $K_1$ we have
\[
\alpha_1=\begin{pmatrix} 0 & 1 & 0 \\ 1 & 0 & 0 \\ 0 & 0 & 1 \end{pmatrix},
\quad
\beta_1=\begin{pmatrix} 1 & 0 & 0 \\ 0 & 0 & 1 \\ 0 & 1 & 0 \end{pmatrix},
\quad\text{and}\quad
(\gamma_{aba=bab})^{-1}_*=\begin{pmatrix} 0 & 0 & 1 \\ 0 & 1 & 0 \\ 1 & 0 & 0 \end{pmatrix}.
\]
Thus, using \eqref{j_I_2odd}, we compute that
\[
\tilde{j}_0=\begin{pmatrix} 2 \\ 1 \end{pmatrix}
\quad\text{and}\quad
\tilde{j}_1=\begin{pmatrix} 1 & 0 & 1 \\ 0 & 2 & 0 \\ 1 & 0 & 1 \\ 0 & 1 & 0 \\ 0 & 1 & 0 \\ 1 & -1 & 1 \end{pmatrix}.
\]
We note that $\ker\tilde{j}_1\cong\Zz$, $\operatorname{im}\tilde{j}_1\cong\Zz^2$, and $\operatorname{coker}\tilde{j}_1\cong\Zz^4$.
Since $\tilde{j}_0$ is injective, we get the exact sequence
\[
0 \longrightarrow K_0(A) \overset{\tilde{j}_0}{\longrightarrow} K_0(A)\oplus K_0(A) \overset{\tilde{\pi}_0}{\longrightarrow} K_0(I) \overset{\tilde{\partial}_0}{\longrightarrow} K_1(A) \overset{\tilde{j}_1}{\longrightarrow} K_1(A)\oplus K_1(A) \overset{\tilde{\pi}_1}{\longrightarrow} K_1(I) \longrightarrow 0,
\]
which gives
\[
0 \longrightarrow \Zz \overset{\tilde{j}_0}{\longrightarrow} \Zz^2 \overset{\tilde{\pi}_0}{\longrightarrow} K_0(I) \overset{\tilde{\partial}_0}{\longrightarrow} \Zz^3 \overset{\tilde{j}_1}{\longrightarrow} \Zz^6 \overset{\tilde{\pi}_1}{\longrightarrow} K_1(I) \longrightarrow 0.
\]
Here, $K_0(I)\cong\operatorname{coker}\tilde{j}_0\oplus\ker\tilde{j}_1\cong\Zz^2$ and $K_1(I)\cong\operatorname{coker}\tilde{j}_1\cong\Zz^4$.
Moreover, using the formula \eqref{iota-pi-formula}, we see that $\iota_i\circ\pi_i\colon  K_i(A) \oplus K_i(A) \oplus K_i(A) \oplus K_i(A) \to K_i(A)$ is the zero map for $i=0$ and
\begin{equation*}
  \rukl{
  \begin{array}{ccc|ccc|ccc|ccc}
  1 & 0 & -1 & 0 & 1 & -1 & 0 & 0 & 0 & 1 & -1 & 0 \\ 
  -1 & 1 & 0 & 0 & -1 & 1 & 1 & -1 & 0 & 0 & 1 & -1 \\ 
  0 & -1 & 1 & 0 & 0 & 0 & -1 & 1 & 0 & -1 & 0 & 1  
  \end{array}
  }
\end{equation*}
for $i=1$ (with respect to $\cE^0 = \gekl{aa,ab,ba,bb}$ with  this ordering). Since $\pi_1$ is surjective ($\img(\pi_1) = \img(\ti{\pi}_1)$ by Remark~\ref{R:j--j}), we compute that for $\iota_1\colon\Zz^4\to\Zz^3$ that $\ker\iota_1\cong\Zz^2$, $\operatorname{im}\iota_1\cong\Zz^2$, and $\operatorname{coker}\iota_1\cong\Zz$.
On the other hand, the image of $\pi_0$ is $\Zz$ (as $\img(\pi_0) = \img(\ti{\pi}_0)$ by Remark~\ref{R:j--j}), so the kernel of $\iota_0$ is either $\Zz$ or $\Zz^2$.
If the former holds, then $\iota_0$ is not zero, so that $[1] \in K_0(A \rtimes_r B_3)$ must be a torsion element, which is not possible, because the canonical trace does not vanish on $[1]$.
Thus, we must have that the kernel is $\Zz^2$, which means that $\iota_0$ is the zero map, so we get
\[
0 \longrightarrow K_0(A) \longrightarrow K_0(A\rtimes_r B_3) \longrightarrow K_1(I) \overset{\iota_1}{\longrightarrow} K_1(A) \longrightarrow K_1(A\rtimes_r B_3) \longrightarrow K_0(I) \longrightarrow 0,
\]
which gives
\[
0 \longrightarrow \Zz \longrightarrow K_0(A\rtimes_r B_3) \longrightarrow \Zz^4 \overset{\iota_1}{\longrightarrow} \Zz^3 \longrightarrow K_1(A\rtimes_r B_3) \longrightarrow \Zz^2 \longrightarrow 0.
\]
Hence, $K_0(C^*_{\lambda}(G)) \cong K_0(A\rtimes_r B_3)\cong\Zz^3$ and $K_1(C^*_{\lambda}(G)) \cong K_1(A\rtimes_r B_3)\cong\Zz^3$.
The identity generates one copy of $\Zz$ in $K_0$, while $[\lambda(a)]=[\lambda(b)]$ and $[\lambda(c)]$ generate two copies of $\Zz$ in $K_1$ by Remark~\ref{abelianization}.
\end{example}

\bremark
In all our examples, we have $K_0(C^*_r(G))\cong H_0G\oplus H_2G$ and $K_1(C^*_r(G))\cong H_1G\oplus H_3G$. The reason is the following: As explained in the introduction of \cite{DH}, there is a manifold model for the classifying space $BG$ with dimension at most $3$ in all our examples, because we are considering Artin-Tits groups of finite type with rank at most $3$. Thus, by \cite[Proposition~2.1.~(ii)]{Mat}, we have $K_*(C^*_{\lambda}(G)) \cong K_*^G(\underline{E}G) = K_*^G(EG) \cong K_*(BG) \cong H_*(G) \oplus H_{*+2}(G)$. For the first identification, we used that $G$ satisfies the Baum-Connes conjecture, while for the second equality, we used that $G$ is torsion-free.
\eremark

\end{document}